\documentclass[11pt]{amsart}

\usepackage{latexsym}
\usepackage{amsfonts}
\usepackage{amssymb}
\usepackage{amsmath}
\usepackage{mathrsfs}
\usepackage{bbm}
\usepackage{dsfont}
\usepackage{tikz-cd}
\usepackage{tikz}
\allowdisplaybreaks
\usepackage[all]{xy}
\usepackage[foot]{amsaddr}
\usepackage{hyperref}
\usepackage[foot]{amsaddr}

\newcommand{\be}{\begin{equation}}
\newcommand{\ee}{\end{equation}}
\newcommand{\bea}{\begin{eqnarray}}
\newcommand{\eea}{\end{eqnarray}}
\newcommand{\bean}{\begin{eqnarray*}}
\newcommand{\eean}{\end{eqnarray*}}
\newcommand{\brray}{\begin{array}}
\newcommand{\erray}{\end{array}}

\newtheorem{dfn}{Definition}[section]
\newtheorem{thm}[dfn]{Theorem}
\newtheorem{lmma}[dfn]{Lemma}
\newtheorem{ppsn}[dfn]{Proposition}
\newtheorem{crlre}[dfn]{Corollary}
\newtheorem{xmpl}[dfn]{Example}
\newtheorem{rmrk}[dfn]{Remark}

\newcommand{\bdfn}{\begin{dfn}\rm}
\newcommand{\bthm}{\begin{thm}}
\newcommand{\blmma}{\begin{lmma}}
\newcommand{\bppsn}{\begin{ppsn}}
\newcommand{\bcrlre}{\begin{crlre}}
\newcommand{\bxmpl}{\begin{xmpl}}
\newcommand{\brmrk}{\begin{rmrk}\rm}

\newcommand{\edfn}{\end{dfn}}
\newcommand{\ethm}{\end{thm}}
\newcommand{\elmma}{\end{lmma}}
\newcommand{\eppsn}{\end{ppsn}}
\newcommand{\ecrlre}{\end{crlre}}
\newcommand{\exmpl}{\end{xmpl}}
\newcommand{\ermrk}{\end{rmrk}}

\newcommand{\bbc}{\mathbb{C}}
\newcommand{\bbz}{\mathbb{Z}}

\newcommand{\bbn}{\mathbb{N}}
\newcommand{\bbr}{\mathbb{R}}

\newcommand{\bbt}{\mathbb{T}}

\newcommand{\cla}{\mathcal{A}}
\newcommand{\clb}{\mathcal{B}}
\newcommand{\cld}{\mathcal{D}}
\newcommand{\cle}{\mathcal{E}}
\newcommand{\clf}{\mathcal{F}}

\newcommand{\clk}{\mathcal{K}}
\newcommand{\cll}{\mathcal{L}}

\newcommand{\clg}{\mathcal{G}}

\author{Md Amir Hossain and S. Sundar}
\title[ Reduced $C^{*}$-algebras of product systems]{  Reduced $C^{*}$-algebras of  product systems---an $E_0$-semigroup and a groupoid perspective}

\address{The Institute of Mathematical Sciences, A CI of Homi Bhabha National Institute, 4th cross street, CIT Campus, Taramani, Chennai, India, 600113}
\email{mdamirhossain18@gmail.com, sundarsobers@gmail.com}

\begin{document}
\maketitle

\begin{abstract}
    For Ore semigroups $P$ with an order unit, we prove that there is a bijection between $E_0$-semigroups over $P$ and product systems of $C^{*}$-correspondences over $P^{op}$. We exploit this bijection and show that the reduced $C^{*}$-algebra of a proper product system is Morita equivalent to the reduced crossed product of the associated semigroup dynamical system given by the corresponding $E_0$-semigroup. We appeal to the groupoid picture of the reduced crossed product of a semigroup dynamical system derived in~\cite{Sundar_Khoshkam} to prove that, under good conditions, the reduced $C^{*}$-algebra of a proper product system is nuclear/exact if and only if the coefficient algebra is nuclear/exact. We also discuss the invariance of $K$-theory under homotopy of product systems. 
\end{abstract}
\noindent {\bf AMS Classification No. :} {Primary 46L55; Secondary 22A22.}  
\tableofcontents

\section{Introduction}
In this paper, we make a contribution to the study of $C^{*}$-algebras associated with product systems. 
Product systems can be studied from various perspectives. They
\begin{enumerate}
    \item are the classifying objects for $E_0$-semigroups in Arveson's classification programme, 
    \item are indispensable in the dilation theory of $CP$-semigroups, and
    \item provide a unifying framework to construct interesting $C^{*}$-algebras which encompasses higher-rank graph $C^{*}$-algebras, semigroup crossed products, etc...
    \end{enumerate}
The purpose of this paper is to demonstrate the utility of viewing product systems as those arising from $E_0$-semigroups in the analysis of the reduced $C^{*}$-algebras associated with them. For, this allows us to write, up to a Morita equivalence, the reduced $C^{*}$-algebra of a product system as a semigroup crossed product for which a nice groupoid description was derived in~\cite{Sundar_Khoshkam}, which makes it possible to apply the well-developed groupoid machinery to answer questions concerning its nuclearity, exactness and K-theory.
The motivation, the context, and the results obtained in this paper are explained in more detail below. 

  The term \emph{product system}  originated in  the seminal work of Arveson~(\cite{Arv_Fock}), and at the end of  a series of fundamental papers (\cite{Arv_Fock}, \cite{Arv_Fock2}, \cite{Arv_Fock3}, \cite{Arv_Fock4}), he established the following  equality
   \begin{equation}
 \label{Arveson}
 \Big\{\textrm{Product systems of Hilbert spaces}\Big\}=\Big\{\textrm{$E_0$-semigroups on $B(H)$}\Big \}
 \end{equation}
  in the one-parameter case.
The above equality states that the problem of classifying $1$-parameter $E_0$-semigroups is equivalent to the problem of classifying product systems over $(0,\infty)$.
 While a product system of Hilbert spaces, over $(0,\infty)$,  is a measurable field of Hilbert spaces over the base space $(0,\infty)$  with an associative multiplication, an $E_0$-semigroup is a $1$-parameter semigroup of unital $*$-endomorphisms of $B(H)$ which can be interpreted as an action of the semigroup $\bbr_{+}=(0,\infty)$ on the algebra $B(H)$. The primary problem in Powers' and Arveson's theory of $E_0$-semigroups is to classify $E_0$-semigroups up to cocycle conjugacy, and the validity of Eq.~\ref{Arveson} was regarded as a seminal result. 
 This reduces the problem to the study of  product systems, which are usually easier to construct than $E_0$-semigroups.
 For example, many interesting classes of product systems can be constructed naturally using probabilistic techniques (\cite{Tsi}, \cite{Tsirelson}, \cite{Liebscher}). 
  Both the notion of product systems and $E_0$-semigroups generalise to the Hilbert module setting and also over general semigroups. We replace Hilbert spaces by Hilbert $C^{*}$-modules, $B(H)$ by the algebra of adjointable operators on a Hilbert module, and $(0,\infty)$ by a general semigroup~$P$. 

On the $C^{*}$-algebra side, 
the seed for studying product systems in the Hilbert module setting was sown by Pimsner who, in his seminal paper (\cite{Pimsner}), associated two $C^{*}$-algebras $\mathcal{T}_E$ and $\mathcal{O}_E$ for a single $C^{*}$-correspondence $E$ (a product system over $\bbn$). The algebra $\mathcal{T}_E$ is called the Toeplitz algebra, and $\mathcal{O}_E$ is called the Cuntz-Pimsner algebra of  $E$. He further demonstrated that the Cuntz algebra $\mathcal{O}_n$,  Cuntz-Krieger algebras $\mathcal{O}_A$ and crossed products by $\bbz$ fall under this class. Moreover, he derived a six-term exact sequence (the Pimsner-Voiculescu six-term exact sequence is a special case of it), a powerful tool, to compute the $K$-theory of the associated algebras. The need to consider product systems over general discrete semigroups became apparent towards the end of the 1990s, and Fowler (\cite{Fowler_2002}) first formally defined the notion of product systems of Hilbert modules over general discrete semigroups. Since then, the study of the associated $C^{*}$-algebras, which include semigroup crossed products and higher-rank graph $C^{*}$-algebras, has remained one of the active areas of research pursued by many operator algebraists. The literature is ever-growing, and the papers listed in the bibliography form a small sample list.

Despite the consensus that semigroup crossed products fall under the broader framework of product systems, and the fact that product systems have their origins in the classification theory of $E_0$-semigroups, a systematic use of $E_0$-semigroups in the analysis of $C^{*}$-algebras associated with product systems over general semigroups is hard to find in the literature. This could partly be explained by the fact that it is only in the recent years (\cite{Sundar_Existence}, \cite{Murugan_Sundar_discrete_mathsci}, \cite{Murugan_Sundar_continuous}) that the bijection between $E_0$-semigroups and product systems is extended beyond the $1$-parameter case to a reasonable class of semigroups even in the Hilbert space setting. The primary motivation of this paper is to show that viewing product systems as those associated with $E_0$-semigroups is advantageous. This allows us to reverse the viewpoint and express the reduced $C^{*}$-algebra of a product system, up to a Morita equivalence, as the reduced $C^{*}$-algebra of a semigroup dynamical system. This, for the case of a single correspondence, was alluded to in the paper of Khoshkam and Skandalis (\cite{KS97}), a major source of inspiration for this work, where they mention at the end of the introductory section that, although their construction falls within the framework of Pimsner's (\cite{Pimsner}), their results would lead to a better understanding of Pimsner's algebras. We share this viewpoint, and we believe that our paper will shed some new light on the study of $C^{*}$-algebras of product systems. We also hope our work brings back into focus the problem of determining the class of semigroups for which Eq.~\ref{Arveson} and its module version hold. However, one drawback with our results is that they are only applicable to product systems that are  \emph{proper}, i.e. when the left action of the coefficient algebra on each fibre is by compact operators.

\subsection{Main results}

Next, we give an overview of the main results obtained. The main results can be divided into two parts. The first part is concerned with expressing (up to a Morita equivalence) the reduced $C^{*}$-algebra of a proper product system as the reduced crossed product of a semigroup dynamical system. The second part deals with applications where we exploit the groupoid picture of the reduced crossed product of a semigroup dynamical system to deduce a few nice consequences.

 Throughout this introduction,  $P$ denotes a discrete, countable, cancellative semigroup with identity $e$, $B$ denotes a separable $C^{*}$-algebra, and $\cle$ stands for a countably generated, full Hilbert $B$-module. The opposite of the semigroup $P$ is denoted $P^{op}$. A product system over~$P$ of $B$-$B$ $C^{*}$-correspondences, or simply a product system with coefficient algebra $B$ is a semigroup~$X$ equipped with a surjective homomorphism $p:X \to P$ such that if we set $X_{s}:=p^{-1}(s)$ for $s \in P$, the fibre $X_s$ is a Hilbert $B$-module equipped with a left action of $B$ such that the map 
 \[X_{s}\otimes_B X_t \ni x \otimes y \to xy \in X_{st}\]
 is a unitary which is also a bimodule map. We assume that $X_{e}=B$, and the `multiplication' map $X_{e} \otimes_{B} X_s \to X_s$ `coincides' with the left action of $B$ on $X_s$ and the multiplication map $X_{s} \otimes_B X_e \to X_s $ is the right action of $B$ on $X_s$. 
 The reduced $C^{*}$-algebra of $X$, denoted $C_{red}^{*}(X)$,  is defined as the $C^{*}$-algebra generated by the left creation operators $\{\phi(x):x \in X_s, s \in P\}$  on the Fock space $\displaystyle \bigoplus_{s \in P}X_s$, where for $x \in X_s$, the operator $\phi(x)$ is defined by 
 \[
 \phi(x)(y \otimes \delta_t):=xy \otimes \delta_{st}.
 \]
 
    A semigroup of unital endomorphisms $\alpha:=\{\alpha_s\}_{s \in P}$ of $\cll_B(\cle)$ is called an $E_0$-semigroup over $P$ on $\cll_B(\cle)$ if for each $s \in P$, $\alpha_s$ is locally strictly continuous \footnote{For a Hilbert $B$-module $\cle$ and a Hilbert $C$-module $\mathcal{F}$, a homomorphism $\pi:\cll_B(\cle) \to \cll_C(\clf)$ is said to be locally strictly continuous if $\pi$ restricted to each norm bounded set is strictly continuous.}, or equivalently, for each $s \in P$, $\alpha_{s}|_{\clk_B(\cle)}$ is non-degenerate. We demand that $\alpha_e$ is the identity morphism. If $\alpha_s$ leaves $\clk_B(\cle)$ invariant for every $s \in P$, we say that $\alpha$ is of compact type. If $A$ is a $C^{*}$-algebra, and if  $\alpha=\{\alpha_s\}_{s \in P}$ is a semigroup of endomorphisms of $A$ such that $\overline{\alpha_s(A)A}=A$ (non-degeneracy condition),  then we call the triple $(A,P,\alpha)$ a semigroup dynamical system. Thanks to the non-degeneracy condition, a semigroup dynamical system $(A,P,\alpha)$ can be viewed as an $E_0$-semigroup on $M(A)$ which is of compact type. Conversely, an $E_0$-semigroup $\alpha$ on $\cll_B(\cle)$ which is of compact type gives rise to the semigroup dynamical system $(\clk_B(\cle),P,\alpha)$. Just like the reduced $C^{*}$-algebra of the product system, we can define the reduced crossed product of a semigroup dynamical system $(A,P,\alpha)$ as follows: for $x \in A$ and $s \in P$, let $\pi(x)$ and $V_s$ be the operators on the external tensor product $\ell^2(P) \otimes A$   defined by 
  \[
  \pi(x)(\delta_t \otimes y):=\delta_{t} \otimes\alpha_t(x) y~; \quad V_{s}(\delta_t \otimes y):=\delta_{ts} \otimes y.\]
 The reduced crossed product is defined as the $C^{*}$-algebra generated  by $\{V_s\pi(x): x \in A,s \in P\}$ and denoted by $A \rtimes_{red} P$.
  
Next, we associate a product system to an $E_0$-semigroup. Arveson's construction of the product system given an $E_0$-semigroup, in the Hilbert space setting, is via the space of intertwiners. However, it does not adapt quite well to the Hilbert module setting. The various subtleties that arise in the Hilbert module setting, and the appropriate way to attach a product system in the module setting were explained in great detail by Skeide in many of his papers~(\cite{Skeide_modulekey}, \cite{Skeide_module1}, \cite{skeide},
\cite{Skeide_module3}). The construction is as follows.

 Suppose $\alpha:=\{\alpha_s\}_{s \in P}$ is an $E_0$-semigroup over $P$ on $\cll_B(\cle)$. For $s \in P$, set \begin{equation}
     \label{fibre}
     X_s:=\cle^{*} \otimes_{\clk_B(\cle)}\cle,
 \end{equation} where in the above internal tensor product, the left action of $\clk_B(\cle)$ on $\cle$ is via the homomorphism $\alpha_s|_{\clk_B(\cle)}:\clk_B(\cle) \to \cll_B(\cle)$. Note that $X_s$ is a Hilbert $B$-module that carries a natural left action of $B$ given by the left action of $B$ on $\cle^{*}$.  Set $\displaystyle X:=\coprod_{s \in P}X_s$. Then, $X$ is a product system over $P^{op}$ with the multiplication rule given by the map 
 \begin{equation}
 \label{multiplication}
\displaystyle X_s \otimes_B X_t \ni (x^* \otimes y) \otimes_B (z^* \otimes u) \to x^* \otimes \alpha_{t} (yz^*)u \in X_{ts}.
 \end{equation}
 We call $X$ the product system associated with $\alpha$. Note that while $\alpha$ is an $E_0$-semigroup over $P$, $X$ is a product system over $P^{op}$. Also, $X$ is proper if and only if $\alpha$ is of compact type. If $B=\bbc$, then $H:=\cle$ is a Hilbert space, and in this case, we can identify the fibre $X_s$ with the space of intertwiners, i.e. 
 \[
 X_s:\cong \{T \in B(H): \alpha_s(A)T=TA \textrm{~for all $A \in B(H)$}\}.
 \]
 However, the multiplication rule then becomes the opposite of the usual composition.  In Arveson's theory and in the Hilbert space setting, $X^{op}$, which is a product system over $P$, is called the product system associated with $\alpha$. In this paper, as we work in the Hilbert module setting, we do not use Arveson's picture. For us, the product system associated with $\alpha$ is $X$, which is a product system over $P^{op}$ and whose fibres are given by Eq.~\ref{fibre} with multiplication given by Eq.~\ref{multiplication}.

The following is our first main result, analogous to a result of Muhly and Solel (\cite{Muhly_Solel_2000}) for a single correspondence. 

\begin{thm}
\label{main_Muhly}
    Let $\alpha$ be an $E_0$-semigroup over $P$ on $\cll_B(\cle)$, and let $X$ be the associated product system over $P^{op}$. Suppose that $\alpha$ is of compact type. Then,  $C_{red}^{*}(X)$ and $\clk_B(\cle)\rtimes_{red} P$ are Morita equivalent.  
\end{thm}
The above theorem is applicable to the study of the reduced $C^{*}$-algebra of a (proper) product system, provided every product system comes from an $E_0$-semigroup. However, this is not true if the semigroup $P$ does not embed in a group (see Section 5 of \cite{Sundar_Existence}). However, it is shown to be true recently by the second author (\cite{Sundar_Existence}) in the Hilbert space setting for a large class of subsemigroups of groups. The techniques used in \cite{Sundar_Existence}, which in turn are inspired by Skeide's ideas (\cite{skeide}) in the $1$-parameter case,   easily adapt themselves to the Hilbert module setting, and we prove the following theorem\footnote{The proof for the base $P=\bbn$ can be found in~\cite{Skeide_modulekey} and~\cite{Skeide_non_unital}.}.

\begin{thm}
    \label{main_existence}
   Let $B$ be a separable  $C^{*}$-algebra, and let $P$ be a subsemigroup  of a group $G$. Assume that $P$ is right Ore, i.e. $PP^{-1}=G$, and that $P$ has an order unit, i.e. there exists $a \in P$ such that $\bigcup_{n=1}^{\infty}Pa^{-n}=G$. Suppose $X$ is a product system over $P^{op}$ with coefficient algebra $B$ such that $X_s$ is full for every $s \in P$. Suppose that $X_s$ is countably generated for every $s \in P$. Then, there exists a countably generated full Hilbert $B$-module $\cle$ and an $E_0$-semigroup $\alpha$ over $P$ on $\cll_B(\cle)$ such that $X$  is isomorphic to the product system associated with $\alpha$.
\end{thm}

Thanks to the above theorem, in the `proper' case and for Ore semigroups with an order unit, the study of the reduced $C^{*}$-algebra of a product system (up to Morita equivalence) boils down to the study of the reduced crossed product of a semigroup dynamical system. 

Next, we turn to the second half of our main results regarding the reduced crossed product of a semigroup dynamical system. Using Thm.~\ref{main_Muhly} and Thm.~\ref{main_existence}, we can easily translate them to the product system context when the product system is proper.  A groupoid picture of the reduced crossed product of a semigroup dynamical system was obtained by the second author in~\cite{Sundar_Khoshkam}. We must also mention here that in the quasi-lattice ordered case and for compactly aligned product systems, the Nica-Toeplitz $C^{*}$-algebra of a product system was expressed as a Fell bundle $C^{*}$-algebra over a groupoid in~\cite{Rennie_Sims}.
The groupoid involved in both~\cite{Sundar_Khoshkam} and~\cite{Rennie_Sims} is the Wiener-Hopf groupoid, which first appeared in the work of Muhly and Renault (\cite{Renault_Muhly}) in their analysis of the Wiener-Hopf algebra associated with a convex cone. The groupoid perspective of the Wiener-Hopf algebra was further developed in~\cite{Nica_WienerHopf}, \cite{Hilgert_Neeb}, and in~\cite{Jean_Sundar}.  

The construction of the Wiener-Hopf groupoid is as follows: let $G$ be a discrete, countable group, and let $P$ be a subsemigroup of~$G$. Let $\mathcal{P}(G)=\{0,1\}^{G}$ be the power set of $G$ endowed with the product topology. Let 
\begin{align*}
\Omega &:=\overline{\{P^{-1}a: a \in P\}}\\\
\widetilde{\Omega} &:=\bigcup_{g \in G}\Omega g.
\end{align*}
The Wiener-Hopf groupoid $\mathcal{G}$ is defined to be the reduction of the transformation groupoid $\mathcal{P}(G) \rtimes G$ onto the clopen set $\Omega$. Moreover, $\mathcal{G}$ and $\widetilde{\Omega} \rtimes G$ are equivalent as groupoids. 
It was proved in \cite{Sundar_Khoshkam} that if $(A,P,\alpha)$ is a semigroup dynamical system, then $A \rtimes_{red} P$ is isomorphic to the reduced crossed product $\mathcal{D} \rtimes_{red} \mathcal{G}$ of a groupoid dynamical system $(\mathcal{D},\clg)$ provided the pair $(P,G)$ satsifies the Toeplitz condition, a technical condition due to Li, which holds if either $(P,G)$ is quasi-lattice ordered or if $PP^{-1}=G$.  Since $\clg$ is equivalent to a transformation groupoid, we can also write $A \rtimes_{red} P$, up to a Morita equivalence, as a reduced crossed product of a group dynamical system. In the product system context, this means that the reduced $C^{*}$-algebra of a proper product system over $P^{op}$ is Morita equivalent to the reduced crossed product of a groupoid dynamical system which in turn is Morita equivalent to a group crossed product when the semigroup $P$ is right Ore and has an order unit. 

We prove two results by applying the groupoid picture of the reduced crossed product of a semigroup dynamical system. We need the following notation to state them: for $x,y \in G$, we say $x \leq y$ if $yx^{-1} \in P$. 
\begin{thm}
\label{main_nuclearity}
Let $P$ be a subsemigroup of a group $G$ such that $(P,G)$ satisfies the Toeplitz condition. Assume that the groupoid $\clg$ (described above) is amenable. 
\begin{enumerate}
    \item Let  $(A,P,\alpha)$ be a semigroup dynamical system with $A$ separable. Then,  we have the following.
\begin{enumerate}
    \item[(i)] If $P^{-1}P=G$ and if $\alpha_s$ is injective for every $s \in P$, then $A \rtimes_{red} P$ is exact if and only if $A$ is exact.
    \item[(ii)] Suppose that for every $F \in \widetilde{\Omega}$, $F$ is directed with respect to the pre-order $\leq$. Then, $A \rtimes_{red} P$ is nuclear if and only if $A$ is nuclear. 
\end{enumerate}
\item Suppose that $P$ has an order unit, i.e. there exists $a \in P$ such that $\bigcup_{n=1}^{\infty}Pa^{-n}=G$. Let $X$ be a proper product system over $P^{op}$ with coefficient algebra $B$. Assume that $B$ is separable and \(X_s\) is full for every $s \in P^{op}$. Then, we have the following. 
\begin{enumerate}
    \item[(a)] If $P^{-1}P=G$ and if the left action of $B$ on $X_s$ is injective for each $s \in P$, then $C_{red}^{*}(X)$ is exact if and only if $B$ is exact. 
    \item[(b)] Suppose that for every $F \in \widetilde{\Omega}$, $F$ is directed. Then, $C_{red}^{*}(X)$ is nuclear if and only if $B$ is nuclear. 
\end{enumerate}
\end{enumerate}
\end{thm}

A similar nuclearity result in the quasi-lattice ordered case was derived in~\cite{Rennie_Sims} using their groupoid Fell bundle picture. The advantage of our result is that it is applicable to a broader class of semigroups, which includes finitely generated subsemigroups of abelian groups.  We also note here that nuclearity and exactness results for product systems over finitely generated subsemigroups of abelian groups were also obtained in~\cite{Laca_et_al} under different hypotheses and by using different methods. 

 As yet another application, we prove the invariance of $K$-theory of the reduced crossed product under homotopy. We consider the $C(Z)$-version of semigroup dynamical systems and product systems. Let $Z$ be a compact metric space. By a $C(Z)$-semigroup dynamical system, we mean a dynamical system $(A,P,\alpha)$, where $A$ is a $C(Z)$-algebra, and $\alpha_s$ is a $C(Z)$-morphism for every $s \in P$, i.e. for $s \in P$, $f \in C(Z)$ and $a \in A$, 
\[
\alpha_s(f\cdot a)=f\cdot \alpha_s(a).
\]
Thanks to the above condition, for every $z \in Z$, we get an action $\alpha^z$ of $P$ on the fibre $A^z$ by endomorphisms resulting in a semigroup dynamical system $(A^z,P,\alpha^z)$.  We define a similar notion of $C(Z)$-product systems. Here, we restrict ourselves to the unital case. A $C(Z)$-product system over $P^{op}$ is a product system $X$ over $P^{op}$ whose coefficient algebra $B$  is a $C(Z)$-algebra such that 
for $f \in C(Z)$, $s \in P$ and $x \in X_s$, 
\[
f\cdot x=x\cdot f.
\]
The above conditions ensure that for every $z \in P$, we get a product system $X^{z}$ with coefficient algebra $B^z$. Here, for $z \in Z$, $B^z$ is the fibre of $B$ over $z$. The following result, the product system version of it and their corollaries form the final main results of our paper. 
\begin{thm}
\label{main_invariance}
     Let $(A,P,\alpha)$ be a $C[0,1]$-semigroup dynamical system. Assume that $A$ is separable. Suppose that the following conditions hold:
     \begin{enumerate}
        \item[(1)] For every $z \in [0,1]$, the evaluation map $K_*(ev_z):K_*(A) \to K_*(A^z)$ is an isomorphism, where $ev_z:A \to A^z$ is the evaluation map.
         \item[(2)]  Every element of $\widetilde{\Omega}$ is directed. 
         \item[(3)] The group $G$ is torsion-free and satisfies the Baum-Connes conjecture with coefficients.
         \item[(4)] The pair $(P,G)$ satisfies the Toeplitz condition.
     \end{enumerate}
     Then, the $K$-groups $K_*(A \rtimes_{red} P)$ and $K_*(A^z \rtimes_{red} P)$ are isomorphic for every $z \in [0,1]$.
\end{thm}
The product system version of the above result holds under similar hypotheses when the semigroup is right Ore with an order unit; in particular, when the semigroup is a finitely generated subsemigroup of an abelian group. In~\cite{Gillaspy}, the notion of homotopy of product systems was considered, and the invariance of $K$-theory under homotopy was proved when the semigroup is $\bbn^k$.  The authors applied the iteration procedure (\cite{Fletcher}, \cite{Deaconu_iterate}) along with the six-term exact sequence in $K$-theory to deduce invariance. But this cannot be used for more general semigroups. However, we could apply our groupoid/crossed product presentation, which makes it possible to make use of the Baum-Connes machinery along with the `descent principle' (\cite{Phillips_Ecterhoff}, \cite{Li-Cuntz-Echterhoff}) in $K$-theory to prove the invariance of $K$-theory.

We end this introduction by indicating the organisation of this paper.  In Section~\ref{sec-preli}, we collect the basic definitions concerning representations of product systems and semigroup dynamical systems. We also prove a few basic facts concerning the reduced crossed product of a semigroup dynamical system. For ease of reference,  we again explain the procedure of attaching a product system to an $E_0$-semigroup. Section~\ref{sec-morita-equi} is dedicated to the proof of Thm.~\ref{main_Muhly}. In Section~\ref{sec-exist-E0-semigroup}, we adapt the techniques of~\cite{Sundar_Existence} and prove Thm.~\ref{main_existence}. The groupoid presentation of a semigroup dynamical system is recalled in Section~\ref{sec-gpd-presentation}. The last two sections are concerned with applications. In Section~\ref{sec-nuc-exact}, we discuss the nuclearity and the exactness of the reduced crossed product of a semigroup dynamical system. Thm.~\ref{main_nuclearity} is proved in this section. In Section~\ref{sec-homo-K-theo}, we prove Thm.~\ref{main_invariance}, which establishes the invariance of $K$-theory under `homotopy'. We end our paper by giving a few examples of product systems that are homotopic.  

The notation used throughout the paper and the standing assumptions are listed below.
 \begin{itemize}
\item We let $\bbn_0=\{0,1,2,\cdots\}$ and $\bbn=\{1,2,\cdots\}$. 
\item For a Hilbert $B$-module $\cle$, the algebra of adjointable operators on $\cle$ is denoted $\cll_B(\cle)$, and we will denote the algebra of compact operators on $\cle$ by $\clk_B(\cle)$. For $x,y \in \cle$, the operator 
\[
\cle \ni z \mapsto x\langle y | z \rangle\ \in \cle\]
is denoted $xy^{*}$. 
\item For a full Hilbert $B$-module $\cle$, let $\cle^{*}:=\{x^{*}:x \in \cle\}$. Then, $\cle^{*}$ is a Hilbert $\clk_B(\cle)$-module with the inner product given by 
\[
\langle x^* | y^{*}\rangle_{\clk_B(\cle)}=xy^{*},\]
and the right action given by 
\[
x^{*}\cdot T=(T^*x)^{*}\]
for $x^{*} \in \cle^{*}$ and $T \in \clk_B(\cle)$. Moreover, $\cle^{*}$ is a $B$-$\clk_B(\cle)$ imprimitivity bimodule, where the left action of $B$ is given by $b\cdot x^*=(xb^*)^{*}$. 

\item In a direct sum $\displaystyle V:=\bigoplus_{i \in I}V_i$, for $i \in I$ and $u \in V_i$, $u \otimes \delta_i$ stands for the element in $V$ which vanishes at all expect at the $i^{th}$-coordinate, where it takes the value $u$. 

\item We only consider separable semigroup dynamical systems, i.e. if $(A,P,\alpha)$ is a semigroup dynamical system, then $A$ is assumed to be separable unless we mention otherwise.

\item We only consider countably generated Hilbert $B$-modules, and we further assume  $B$ is separable. So, we assume that the fibres of a product system are countably generated.

\item If $\cla$ is an upper semi-continuous bundle of $C^{*}$-algebras over a locally compact space~$X$, then the  algebra of continuous sections of $\cla$  is denoted $C_0(X,\cla)$ and $C(X,\cla)$ if $X$ is compact. 
\item The semigroups considered are countable, cancellative, discrete and have an identity element $e$.

\end{itemize}

\section{Preliminaries}
\label{sec-preli}

In this section, we collect the basic definitions concerning representations of semigroup dynamical systems and product systems. We also prove a few results concerning the reduced crossed product of a semigroup dynamical system that we will use later. Till the end of Section~\ref{sec-morita-equi}, unless otherwise mentioned, the letter $P$ stands for a discrete, countable, cancellative semigroup with an identity $e$.

\subsection{Reduced crossed product of a semigroup dynamical system}
\begin{dfn}\label{def-cov-rep-dyn-syt}
Let $(A,P,\alpha)$ be a semigroup dynamical system. Let $\cle$ be a Hilbert $B$-module, let $\pi:A \to \cll_B(\cle)$ be a homomorphism, and let $V:=\{V_s\}_{s \in P}$ be a family of adjointable isometries on $\cle$. We say that $(\pi,V)$ is a covariant representation of $(A,P,\alpha)$ on $\cle$ if 
 \begin{enumerate}
 \item $\pi$ extends to a locally strictly continuous representation of $M(A)$,
 \item for $s \in P$ and $x \in A$, $\pi(x)V_s=V_s\pi(\alpha_s(x))$, and
 \item for $s,t \in P$, $V_{s}V_{t}=V_{ts}$. 
 \end{enumerate}
 The representation $\pi$ is said to be non-degenerate if $\overline{\pi(A)\cle}=\cle$. 
\end{dfn}

Let $(A,P,\alpha)$ be a semigroup dynamical system. 
Here is an example of a covariant representation.  
 Consider the Hilbert $A$-module $\cle:=\ell^2(P) \otimes A$. Here, $\otimes$ denotes the external tensor product. Let $\{\delta_s:s \in P\}$ be the standard orthonormal basis of $\ell^2(P)$. For $s \in P$, let $v_s:\ell^2(P) \to \ell^2(P)$ be the isometry defined by \begin{equation}
 \label{isometries}
 v_s(\delta_t)=\delta_{ts}.\end{equation} Set \[V_s:=v_s \otimes 1 \] for $s \in P$. 
For $x \in A$, let $\pi(x):\cle \to \cle$ be the adjointable operator defined by 
\begin{equation}
\label{regular_representation}
\pi(x)( \delta_s \otimes x):=\delta_s \otimes \alpha_s(x)y .\end{equation}
Note that for each $s \in P$, $V_s$ is an isometry, and $V_sV_t=V_{ts}$ for $s,t \in P$. Moreover, the pair $(\pi,V)$ is a covariant representation. 
The pair $(\pi,V)$ is called the \emph{regular representation} of the dynamical system $(A, P,\alpha)$.

\begin{dfn}[\cite{Sundar_Khoshkam}]
\label{first_Wiener}
Let $(A,P,\alpha)$ be a semigroup dynamical system. Let $(\pi,V)$ be the regular representation of 
$(A,P,\alpha)$. The reduced crossed product $A\rtimes_{red} P$ is defined as the $C^{*}$-algebra generated by $\{V_s\pi(x):s \in P,x \in A\}$.

    Suppose $P$ is a subsemigroup of a group $G$. Let $\rho:=\{\rho_g\}_{g \in G}$ be the right regular representation of $G$ on $\ell^2(G)$. For $g \in G$, let $w_g$ be the compression of $\rho_g$ onto $\ell^2(P)$. Set $W_g:=w_g \otimes 1$. The \emph{Wiener-Hopf algebra} is defined as the $C^{*}$-algebra generated by $\{\pi(x)W_g: x \in A, g \in G\}$ and denoted by $\mathcal{W}(A,P,G,\alpha)$.
\end{dfn}

\begin{rmrk}
In \cite{Sundar_Khoshkam}, the reduced crossed product $A \rtimes_{red} P$ was defined as the $C^{*}$-algebra generated by $\{\pi(x)V_{s_1}^{*}V_{t_1}V_{s_2}^{*}V_{t_2}\cdots V_{s_n}^{*}V_{t_n}: n \in \bbn_0, x \in A, s_i, t_i \in P\}$. Using the fact that $(\pi,V)$ is covariant and $\alpha_s$ is non-degenerate for every $s \in P$, it is easily verifiable that  that the $C^{*}$-algebra generated by   $\{\pi(x)V_{s_1}^{*}V_{t_1}V_{s_2}^{*}V_{t_2}\cdots V_{s_n}^{*}V_{t_n}: n \in \bbn_0, x \in A, s_i, t_i \in P\}$, and the $C^{*}$-algebra generated by $\{V_s\pi(x):s \in P,x \in A\}$ are equal. 
  The case when $P=\bbn$ was first considered by Khohskam and Skandalis (\cite{KS97}). 
  
  If $A=\bbc$, we denote  $\bbc \rtimes_{red} P$ by $C_{red}^{*}(P)$. The $C^{*}$-algebra $C_{red}^{*}(P)$ is called  the reduced $C^{*}$-algebra of the semigroup $P$.

  \end{rmrk}

\begin{ppsn}
\label{non-degeneracy-faithful-reduced-basic}
    Let $(A,P,\alpha)$ be a semigroup dynamical system, and let $(\pi,V)$ be the regular representation. Then, $\pi$ is faithful, and $\overline{\pi(A)(A \rtimes_{red} P)}=A \rtimes_{red} P$. 
\end{ppsn}
\textit{Proof.} Since we have assumed that $P$ contains the identity element and $\alpha_e$ is the identity map, it is clear that $\pi$ is one-one. Let $(e_\lambda)_{\lambda}$ be an approximate unit of $A$. Then, for $s \in P$, $x \in A$, and since $\alpha_s$ is locally strictly continuous,
\[
\lim_{\lambda}\pi(e_\lambda)V_s\pi(x)=\lim_{\lambda}V_s\pi(\alpha_s(e_\lambda))\pi(x)=\lim_{\lambda}V_{s}\pi(\alpha_s(e_\lambda) x)=V_s\pi(x). 
\]
Also, $\lim_{\lambda}V_s\pi(x)\pi(e_\lambda)=\lim_{\lambda}V_{s}\pi(xe_\lambda)=V_{s}\pi(x)$. Since $\{V_{s}\pi(x):s \in P,x \in A\}$ generates $A \rtimes_{red} P$, $\pi(A)(A \rtimes_{red} P)$ is total in $A \rtimes_{red} P$. \hfill $\Box$

\begin{rmrk}
\label{Toeplitz condition}
    If $P \subset G$ is a subsemigroup, we say that $(P,G)$ satisfies the Toeplitz condition if for each $g$, the partial isometry $w_g$, as defined in Defn.~\ref{first_Wiener}, lies in the semigroup generated by $\{v_s,v_t^*:s,t \in P\} \cup \{0\}$, where $\{v_s\}_{s \in P}$ is as defined in Eq. \ref{isometries}. If $(P,G)$ satisfies the  Toeplitz condition, then for a semigroup dynamical system $(A,P,\alpha)$, the reduced crossed product $A \rtimes_{red} P$ and $\mathcal{W}(A,P,G,\alpha)$ coincide. The Toeplitz condition is satisfied if $PP^{-1}=G$, or if $(P,G)$ is quasi-lattice ordered (see Section 2 of \cite{Sundar_Khoshkam}). 
\end{rmrk}

\begin{ppsn}
Let $(A,P,\alpha)$ be a semigroup dynamical system. Let $(\pi,V)$ be the regular representation of $(A,P,\alpha)$.
Then,  $\pi(x) \in M(A \rtimes_{red} P)$ and $V_s \in M(A\rtimes_{red} P)$ for $x \in A$ and $s \in P$.
Suppose $D$ is a $C^{*}$-algebra, and $L$ is a Hilbert $D$-module. Let  $\eta:A \rtimes_{red} P \to \cll_D(L)$ be a non-degenerate representation. Then, there exists a covariant representation $(\widetilde{\pi},\widetilde{V})$ of $(A,P,\alpha)$ such that 
\[
\eta(V_s\pi(x))=\widetilde{V}_s\widetilde{\pi}(x)\]
for $x \in A$ and $s \in P$.
\end{ppsn}
\textit{Proof.}  By definition, $\pi(A) \subset A \rtimes_{red} P \subset M(A \rtimes_{red} P)$. Let
\[
\mathcal{S}:=\{V_s\pi(x):s \in P, x \in A\}.\] 
Let $t \in P$. We claim the following:
\begin{enumerate}
    \item $V_t\mathcal{S} \subset A \rtimes_{red} P$,
    \item $\mathcal{S}V_t \subset A \rtimes_{red} P$,
    \item $V_{t}^{*}\mathcal{S} \subset A \rtimes_{red} P$, and 
    \item $\mathcal{S}V_{t}^{*} \subset A \rtimes_{red} P$. 
\end{enumerate}
Note that $(1)$ and $(2)$ follows from the fact that $(\pi,V)$ is a covariant representation. More precisely, $(1)$ follows from the third condition of Defn.~\ref{def-cov-rep-dyn-syt} and $(2)$ follows from the second condition of Defn.~\ref{def-cov-rep-dyn-syt}. 

Let $d \in \mathcal{S}$ be given, and write $d=V_s\pi(x)$ for some $s \in P$ and $x \in A$. Let $(e_\lambda)_{\lambda \in \Lambda}$ be an approximate identity of $A$. Note that 
\begin{align*}
    V_t^{*}d &=V_{t}^{*}(V_s\pi(x))\\
    &=\lim_{\lambda \in \Lambda}V_{t}^{*}(V_s\pi(\alpha_{s}(e_{\lambda})x) \quad  (\textrm{since $\alpha_s$ and $\pi$ are locally strictly continuous})
\\
&=\lim_{\lambda \in \Lambda}V_{t}^{*}\pi(e_\lambda)V_s\pi(x)\\
&=\lim_{\lambda \in \Lambda}(\pi(e_\lambda)V_t)^{*}V_s\pi(x)\\
&=\lim_{\lambda \in \Lambda}(V_t\pi(\alpha_t(e_\lambda))^{*}V_s\pi(x) \in A\rtimes_{red} P.
\end{align*}
The proof of $\mathcal{S}V_t^{*} \subset A \rtimes_{red} P$ is similar. This proves the claim. 

Since $\mathcal{S}$ generates $A \rtimes_{red} P$, it follows that $A \rtimes_{red} P$ is left invariant when multiplied by $\{V_t,V_{t}^{*}\}$ both on the right as well as on the left. Hence, $V_{t} \in M(A \rtimes_{red} P)$. The proof of the first assertion is over. The second assertion follows from the first. \hfill $\Box$

We next derive a few good functorial properties of the reduced crossed product. The proofs are analogous to the case when $P$ is a group.  Let $(A,P,\alpha)$ be a dynamical system.  Let $(\pi,V)$ be the regular representation of $(A,P,\alpha)$ on the Hilbert $A$-module $\ell^2(P)\otimes A$. Let $\rho$ be a faithful representation of $A$ on a Hilbert space $H$. For $s \in P$ and $x \in A$, let $\widetilde{\pi}(x)$ and $\widetilde{V}_s$ be the operators on $\ell^2(P)\otimes H$ defined by 
\[
\widetilde{V}_{s}(\delta_t \otimes \eta)=\delta_{ts}\otimes \eta~; \quad \widetilde{\pi}(x)(\delta_t \otimes \eta)=\delta_t \otimes \rho(\alpha_t(x))\eta.
\]
Then, $(\widetilde{\pi},\widetilde{V})$ is a covariant representation. 

\begin{ppsn}
\label{alternate_description_reduced}
With the foregoing notation, there exists a unique injective  $*$-homomorphism $\mu:A \rtimes_{red} P \to \cll(\ell^2(P) \otimes H)$ such that 
\[
\mu(V_s\pi(x))=\widetilde{V}_s\widetilde{\pi}(x)
\]
for $s \in P$ and $x \in A$.  Hence, the $C^{*}$-algebra $A \rtimes_{red} P$ is isomorphic to  the $C^{*}$-subalgebra of \(\mathcal{L}(\ell^2(P)\otimes H)\) generated by $\{\widetilde{V}_s\widetilde{\pi}(x):x \in A, s \in P\}$.
    \end{ppsn}
    \textit{Proof.} 
    Let $U:(\ell^2(P) \otimes A) \otimes_{\rho} H \to \ell^2(P)\otimes H$ be the isometry given by the equation     \[
    U((\delta_s \otimes a) \otimes_{\rho}\xi)=\delta_s \otimes \rho(a)\xi.
    \]
    for $s \in P$, $a \in A$ and $\xi \in H$. 
    It follows from a routine computation that for $x \in A$, $s \in P$, 
    \[
    U(V_{s}\pi(x)\otimes 1)U^{*}=\widetilde{V}_{s}\widetilde{\pi}(x).
    \]
    Let $\Phi:\cll_{A}(\ell^2(P) \otimes A) \to \cll((\ell^2(P) \otimes A)\otimes_{\rho}H)$ be the map defined by \[\Phi(T):=T \otimes 1.\] Since $\rho$ is faithful, $\Phi$ is faithful. The map $\mu(\cdot):=U\Phi(\cdot)U^*$ is the desired homomorphism. \hfill $\Box$

The following two corollaries are immediate. 

\begin{crlre}
\label{functorial_one}
Let $(A,P,\alpha)$ and $(B,P,\beta)$ be two semigroup dynamical systems. Let $\epsilon:B \to A$ be an injective homomorphism such that $\epsilon$ is $P$-equivariant, i.e. $\epsilon(\beta_s(b))=\alpha_s(\epsilon(b))$ for $b \in B$ and $s \in P$. Let $(\pi,V)$ be the regular representation of $(A,P,\alpha)$ and $(\widehat{\pi},\widehat{V})$ be the regular representation of $(B,P,\beta)$. Then, there exits a unique $*$-homomorphism $\widetilde{\epsilon}:B \rtimes_{red} P \to A \rtimes_{red} P$ such that $\widetilde{\epsilon}(\widehat{V}_s\widehat{\pi}(b))=V_s\pi(\epsilon(b))$ for every $s \in P$ and $b \in B$.
\end{crlre}

\begin{crlre}
\label{functorial_ideal}
Let $(A,P,\alpha)$ be a semigroup dynamical system, and let $I \subset A$ be an ideal which is $P$-invariant, i.e. $\alpha_s(I) \subset I$ for $s \in P$. Suppose that  $\overline{\alpha_s(I)I}=I$ for every $s \in P$. Let $(\pi,V)$ be the regular representation of $(A,P,\alpha)$ and $(\widehat{\pi},\widehat{V})$ be the regular representation of $(I,P,\alpha)$. Then, there exists a unique $*$-homomorphism $\mu: I \rtimes_{red} P \to A \rtimes_{red} P$ such that $\mu(\widehat{V}_s\widehat{\pi}(x))=V_s\pi(x)$ for $x \in I$ and $s \in P$. Moreover, $\mu$ is injective. In short, $I \rtimes_{red} P \subset A \rtimes_{red} P$ and is also an ideal.    
\end{crlre}

\begin{ppsn}
\label{equivariance implies homomorphism}
    Let $(A,P,\alpha)$ and $(B,P,\beta)$ be semigroup dynamical systems. Let $(\pi,V)$ be the regular representation of $(A,P,\alpha)$, and let $(\widetilde{\pi},\widetilde{V})$ be the regular representation of $(B,P,\beta)$.  Let $\phi:A \to B$ be a homomorphism which is $P$-equivariant, i.e. $\phi \circ \alpha_s=\beta_s\circ \phi$ for every $s \in P$.  Suppose that $\overline{\phi(A)B}=B$. Then, there exists a unique $*$-homomorphism $\widetilde{\phi}:A \rtimes_{red} P \to B \rtimes_{red} P$ such that $\widetilde{\phi}(V_s\pi(x))=\widetilde{V}_s\widetilde{\pi}(\phi(x))$ for $x \in A$ and $s \in P$. 
\end{ppsn}
\textit{Proof.} Let $U:(\ell^2(P) \otimes A)\otimes_{\phi}B \to \ell^2(P) \otimes B$ be the unitary defined by 
\[
U((\delta_s \otimes a) \otimes_{\phi} b)=\delta_{s}\otimes \phi(a)b.
\]
for $s \in P$, $a \in A$ and $b \in B$. 
Let $\Phi:\cll_{A}(\ell^2(P) \otimes A) \to \cll_{B}((\ell^2(P) \otimes A)\otimes_{\phi} B)$ be defined by $\Phi(T):=T \otimes 1$. Observe that, for $x \in A$ and $s \in P$, 
\[
U\Phi(V_s\pi(x))U^{*}=\widetilde{V}_s\widetilde{\pi}(x).
\]
Then, $\widetilde{\phi}(\cdot)=U\Phi(\cdot)U^*$, restricted to $A \rtimes_{red} P$, is the required map. \hfill $\Box$

\subsection{Product systems} 
Let $A$ and $B$ be $C^{*}$-algebras. An $A$-$B$ $C^{*}$-correspondence  $X$  is a Hilbert $B$-module together with a left action of $A$ given by a homomorphism $\pi:A \to \cll_B(X)$. The $C^{*}$-correspondence $X$ is said to \emph{proper} if $\pi(A) \subset \clk_B(X)$ and is said to be injective if $\pi$ is faithful. It is said to be regular if it is both proper and injective. We say that $X$ is non-degenerate if $\pi$ is non-degenerate. 
The homomorphism $\pi$ is usually suppressed, and we denote $\pi(b)x$ simply by $b\cdot x$ or $bx$ for $b \in A$ and $x \in X$. We also call an $A$-$B$ $C^{*}$-correspondence a $C^{*}$-correspondence from $A$ to $B$. 

Let $X$ be a product system over $P$ with coefficient algebra $B$. For $s \in P$, denote the fibre over $s $ by $X_s$, which is a $B$-$B$-correspondence. We also call $X$ a product system of $B$-$B$-correspondences. We often abuse  notation, and write $X:=\{X_s\}_{s \in P}$ instead of $X:=\displaystyle \coprod_{s \in P}X_s$. For $s,t \in P$, denote the map \[
X_{s} \otimes_B X_t \ni u \otimes v \mapsto uv \in X_{st}\] by $U_{s,t}$ which is a bimodule map. Recall that we have assumed $X_e=B$. This implies in particular that $X_s$ is non-degenerate for every $s\in P$. The associativity of the multiplication on $X$ is equivalent to the fact that for $r,s,t \in P$, 
\[
U_{r,st}(1 \otimes U_{s,t})=U_{rs,t}(U_{r,s} \otimes 1).
\]
Note that $1 \otimes U_{s,t}$ makes sense as $U_{s,t}$ is a bimodule map. 

\begin{rmrk}
 We mention here that there are important/natural examples where non-degeneracy fails. We do not strive for generality, and we demand non-degeneracy in this paper.  We refer the reader to Katsura's papers (\cite{Katsura1}, \cite{Katsura2}, \cite{Katsura3}) for more on these issues.
\end{rmrk}

\begin{dfn}\label{def-rep-prod-mod}
Let $X:=\{X_s\}_{s \in P}$ be a product system of $B$-$B$-correspondences over $P$. Let~$\cle$ be a Hilbert $D$-module, where $D$ is a $C^{*}$-algebra. Let $\displaystyle \phi:\coprod_{s \in P}X_s \to \cll_D(\cle)$ be a map whose restriction to the fibre $X_s$ is denoted by $\phi_s$. We call $\phi=\{\phi_s\}_{s \in P}$ a representation of $X$ on $\cle$ if 
\begin{enumerate}
\item for $s \in P$, $u,v \in X_s$, $\phi_s(u)^{*}\phi_s(v)=\phi_{e}(\langle u|v \rangle)$,
\item for $s \in P$, there exists an adjointable isometry $T_s:X_{s} \otimes_{B} \cle \to \cle$, necessarily unique,  such that 
\[
T_s(u \otimes \xi)=\phi_s(u)\xi\]
for $u \in X_s$ and $\xi \in \cle$,
\item for $s,t \in P$, $u \in X_s$ and $v \in X_t$, $\phi_{st}(uv)=\phi_s(u)\phi_t(v)$, and 
\item the representation $\phi_e$ extends to a locally strictly continuous representation of $M(B)$. 
\end{enumerate}
We call $\phi$ non-degenerate if $\phi_e$ is non-degenerate.
\end{dfn}
Every product system carries a natural representation called the Fock representation, which is defined below. Let $X:=\{X_s\}_{s \in P}$ be a product system of $B$-$B$ correspondences. Let 
\[
H:=\bigoplus_{s \in P}X_s
\]
be the full Fock module. For $s \in P$ and $u \in X_s$, define the left creation operator  $\phi_s(u)$ by setting
\[
\phi_s(u)(v \otimes \delta_t):=uv \otimes \delta_{st}.
\]
Then, $\phi:=\{\phi_s\}_{s \in P}$ is a representation of $X$ on the Hilbert $B$-module $H$. We call $\phi$ the \emph{Fock representation} or the \emph{regular representation} of the product system $X$. Note that the Fock representation is non-degenerate as we have assumed that the left action on $X_s$ is non-degenerate for every $s \in P$. 
Recall that the \emph{reduced $C^{*}$-algebra of $X$}, denoted by $C_{red}^{*}(X)$, is defined to be the $C^{*}$-subalgebra of $\cll_B(H)$ generated by $\{\phi_s(u): u \in X_s, s \in P\}$.

\subsection{From $E_0$-semigroups to product systems}

Although we explained in the introduction how to associate a product system to an $E_0$-semigroup, for ease of reference, we recall once again. Arveson's space of intertwiners needs to be appropriately replaced, and the idea for the modification required could be traced back to Rieffel's work (\cite{Rieffel}).

Let $B$ be a $C^{*}$-algebra,  and let $\cle$ be a full Hilbert $B$-module. Suppose $\alpha:=\{\alpha_s\}_{s \in P}$ is an $E_0$-semigroup over $P$ on $\cll_B(\cle)$. Fix $s \in P$. By restricting $\alpha_s$ to $\clk_B(\cle)$, we can view $\alpha_s$ as a non-degenerate representation of $\clk_B(\cle)$ on $\cle$. Since $\clk_B(\cle)$ and $B$ are Morita-equivalent with $\cle$ being an imprimitivity bimodule, it follows that the representation $\alpha_s$ arises out of a representation of $B$ on another Hilbert $B$-module call it $X_s$ via Rieffel's induction (see ~\cite[Thm. 5.3]{Rieffel} and~\cite[Thm. 1.4]{Skeide_Muhly_Solel} for more details). A more precise description is given below. 

Let $s \in P$. Let $\cle_s:=\cle$ be the  $\clk_B(\cle)$-$B$-correspondence, where the left action of $\clk_B(\cle)$ on $\cle$ is given by the homomorphism $\alpha_s:\clk_B(\cle) \to \cll_B(\cle)$. 
Define \[X_s:=\cle^{*} \otimes_{\clk_B(\cle)} \cle_s.\]
Then, $X_s$ is a $B$-$B$-correspondence as $\cle^{*}$ is a $B$-$\clk_B(\cle)$-correspondence. The $B$-valued inner product on $X_s$ is given by the formula 
\[
\langle x_1^{*} \otimes y_1 | x_2^{*} \otimes y_2 \rangle_B=\langle y_1|\alpha_s(x_1x_2^{*})y_2 \rangle.\]
The left action of $B$ on $X_s$ is given by 
\[
b\cdot (x^{*} \otimes y)=(xb^*)^*\otimes y.\]
We also write $\cle^{*}\otimes_{\alpha_s} \cle$ instead of $\cle^{*} \otimes_{\clk_B(\cle)}\cle_s$ if we wish to stress the left action on the second factor.

 It is not difficult to prove using the fact that $\alpha_s \circ \alpha_t=\alpha_{st}$ for  $s,t \in P$,  and the fact that $\alpha_s$ is non-degenerate for each $s \in P$ that given $s,t \in P$, there exists a unitary $U_{s,t}:X_s \otimes_B X_t \to X_{ts}$ such that for 
\[
U_{s,t}((y^* \otimes z) \otimes_B (x^* \otimes u)) = y^* \otimes \alpha_{t}(zx^*)u.\]
Note that $U_{s,t}$ is a bimodule map. 
The product on the disjoint union $\displaystyle X:=\coprod_{s}X_s$ is defined by setting 
\[
u\cdot v=U_{s,t}(u \otimes v)\]
for $u \in X_s$ and $v \in X_t$.   The  product is associative~(\cite[Thm. 1.14]{Skeide_Muhly_Solel}) and makes $X:=\{X_s\}_{s \in P^{op}}$ a product system over $P^{op}$. 
We call $X$ \emph{the product system of the $E_0$-semigroup~$\alpha$.}

\begin{rmrk}
\label{absorption}
Observe that for $s \in P$, the map $\sigma_s:\cle \otimes_{B} X_s \to X_s$ defined by 
\[
\sigma_s(x \otimes (y^* \otimes z))=\alpha_s(xy^*)z\]
is a unitary. Thus, $\cle \otimes_{B} X_s \cong \cle$. Moreover,   $\alpha_s(T)=\sigma_s(T \otimes 1)\sigma_s^*$ for every $s \in P$ and $T \in \clk_B(\cle)$. We refer to the isomorphism $\cle \otimes_B X_s \cong X_s$ given by $\sigma_s$ as the `absorption property' of $\cle$. 
\end{rmrk}

	\begin{rmrk}\label{lem:comp-faithful}
			Let \(B\) and \(C\) be two \(C^*\)-algebras. Consider a Hilbert \(B\)-module \(Y\), and let $Z$ be a  \(C^*\)-correspondence from $B$ to $C$. Let \(\Phi: \mathcal{L}_B(Y) \to \mathcal{L}_C(Y \otimes_B Z)\) be the homomorphism defined by \(\Phi(T):= T \otimes 1\).
	\begin{enumerate}
		\item  If the left action of \(B\) on \(Z\) is faithful, then \(\Phi\) is also faithful.
		
	\item If $Z$ is proper, then \(\Phi(T)=T \otimes 1 \in \mathcal{K}_{C}(Y\otimes_{B}Z)\) for \(T\in \mathcal{K}_{B}(Y)\) (see~\cite[Lemma 3.2]{Kwasniewski-Topo-aperiod-prod-syst}).
	\end{enumerate}
	\end{rmrk}

\begin{ppsn}\label{prop-inj-comp-action}
Let $\cle$ be a full Hilbert $B$-module, and let  $\alpha$ be an $E_0$-semigroup over $P$ on $\cll_B(\cle)$. Denote the associated product system over $P^{op}$ by $X$. We have the following. 
\begin{enumerate}
\item The $E_0$-semigroup $\alpha$ is of compact type if and only $X$ is proper. 
\item For $s \in P$, the homomorphism $\alpha_s:\cll_B(\cle) \to \cll_B(\cle)$ is injective if and only if the left action of $B$ on $X_s$ is injective. 
\end{enumerate}
\end{ppsn}
\textit{Proof.}
We use the notation of Remark \ref{absorption}.

\noindent (1).  Suppose the left action of the algebra \(B\) on \(X_s\) is by compact operators for every \(s \in P\). 
	Let \(T\) be an element of \(\mathcal{K}_B(\mathcal{E})\). Then, by  Remark~\ref{lem:comp-faithful}(2), \( T \otimes I \in \mathcal{K}_B(\mathcal{E} \otimes_B X_s)\). It follows from Remark \ref{absorption} that  \(\alpha\) is of compact type. 	Conversely, suppose \(\alpha\) is of compact type. By Remark~\ref{lem:comp-faithful}(2) we conclude that the action of \(B\)  on \(\mathcal{E}^*\otimes_{\mathcal{K}_B(\mathcal{E})} \mathcal{E}=\cle^{*}\otimes_{\alpha_s}\cle\) is by compact operators as \(B\) acts on \(\mathcal{E}^*\) by compact operators. Thus, the left action of \(B\) on \(X_s\) is by compact operators for every $s \in P$.

	\noindent (2).  Suppose the left action of \(B\) on \(X_s\) is injective for $s \in P$. Let \(T\in \mathcal{K}_B(\mathcal{E})\) be such that  \(\alpha_s(T) = 0\). Then,  \(\alpha_s(T)=\sigma_s(T\otimes I)\sigma_s^{*} =0 \implies T \otimes 1=0\). By Remark \ref{lem:comp-faithful}(1), $T=0$. 
	Conversely, assume that \(\alpha_s\) is injective for $s \in P$. Then, by Remark~\ref{lem:comp-faithful}(1), the left action of \(B\) on \(X_s:=\cle^{*}\otimes_{\alpha_s}\cle\) is faithful.	
\hfill $\Box$

\section{A Morita equivalence result}
\label{sec-morita-equi}

In this section, we prove Thm.~\ref{main_Muhly}, which is our first main result. The proof will be presented after several propositions. We denote the opposite semigroup $P^{op}$ by $Q$. We start with the following remark.

\begin{rmrk}
Thm.~\ref{main_Muhly} is the product system version of a result of  Muhly and Solel ([\cite{Muhly_Solel_2000}, Corollary 2.11]). In ~\cite{Muhly_Solel_2000}, Muhly and Solel defined a notion of Morita equivalence of $C^{*}$-correspondences (product systems over $\bbn$) and showed that the  Toeplitz algebra and the Cuntz-Pimsner algebra associated with $C^{*}$-correspondences are Morita equivalent if the underlying $C^{*}$-correspondences are Morita equivalent. It is quite probable that the product system version of this result is known to experts, and Thm.~\ref{main_Muhly} might follow from such a result.

 For, if 
$(A,P,\alpha)$ is a semigroup dynamical system, and if we view it as an $E_0$-semigroup over $P$ on $\cll_A(A)$, then the construction explained in Section~\ref{sec-preli} produces a product system $Y:=\{Y_s\}_{s \in P}$ of $A$-$A$ $C^{*}$-correspondences over $Q=P^{op}$. It follows from definition that $C_{red}^{*}(Y) =A \rtimes_{red} P$. Now, let $\alpha$ be a compact type $E_0$-semigroup over $P$ on $\cll_B(\cle)$. Then,  this construction applied to the dynamical system $(\clk_B(\cle),P,\alpha)$ yields a product system $Y$ of $\clk_B(\cle)$-$\clk_B(\cle)$-correspondences over $P^{op}$. Denote by $X$ the product system of $B$-$B$-correspondences over $P^{op}$ defined as in the previous section. 
The fibres $Y_s$ and $X_s$ are related by the equation
\begin{equation}
\label{Muhly_Solel_equation}
\cle \otimes_{B} X_s\cong Y_s \otimes_{\clk_B(\cle)} \cle,   
\end{equation}
where the isomorphism is given by the map \[\cle \otimes_B X_s \ni x \otimes (y^* \otimes z) \mapsto \alpha_s(xy^*)\otimes z \in Y_s \otimes_{\clk_B(\cle)}\cle. \] Eq.~\ref{Muhly_Solel_equation}  is exactly the equation considered by Muhly and Solel~(\cite{Muhly_Solel_2000}). 
However, we do not strive here to define the notion of Morita equivalence of two product systems, and prove the Morita equivalence of the associated reduced $C^{*}$-algebras, in full generality, working completely in the language of product systems. 

We will be content with proving Thm.~\ref{main_Muhly}. For, we wish to stress the dynamics governed by $\alpha$ and to think of $(\clk_B(\cle),P,\alpha)$ as a genuine semigroup dynamical system and not just as a part of the theory of product systems. In fact, the objective of our paper is to establish that there are advantages in taking the opposite viewpoint, which is that product systems come from semigroup dynamical systems.
\end{rmrk}

For the rest of this section, let $B$ be a separable $C^{*}$-algebra, $\cle$ a countably generated full Hilbert $B$-module, and let $\alpha$ be an $E_0$-semigroup over $P$ on $\cll_B(\cle)$. We denote the product system over $Q:=P^{op}$ associated with $\alpha$ by $X$. We assume that $\alpha$ is of compact type, i.e. $\alpha_s(\clk_B(\cle)) \subseteq \clk_B(\cle)$ for every $s \in P$.
  We first discuss a way to pass from a covariant representation of $(\clk_B(\cle),P,\alpha)$ to a representation of  $X$ and vice versa. 
  
\begin{ppsn}
    \label{equivalence of representation theory1}
    Let $D$ be a $C^{*}$-algebra, and let $H$ be a Hilbert $D$-module. 
    Let $(\pi,V)$ be a non-degenerate covariant representation of the dynamical system $(\clk_B(\cle),P,\alpha)$ on $H$. Let \[L:=\cle^{*}\otimes_{\clk_B(\cle)}H=\cle^{*}\otimes_\pi H.\]  Then, there exists a representation $\phi^{(\pi,V)}:=\phi=\{\phi_s\}_{s \in Q}$  of $X$ on $L$ such that for $s \in P$, \[
    \phi_{s}(x^*\otimes y)(z^* \otimes h)=\phi_s(x^*\otimes_{\alpha_s} y)(z^* \otimes h)=x^* \otimes V_s\pi(yz^*)h.
    \]
    for $x,y,z \in \cle$ and $h \in H$. 
\end{ppsn}
\textit{Proof.} 
	Let \(s\in P\).
		Define a unitary map \(U_s\colon X_s\otimes_{B}(\mathcal{E}^* \otimes_{\pi} H) \to \mathcal{E}^*\otimes_{\pi \circ \alpha_s} H\)
		by 
		\[
		U_s\big((x^*\otimes_{\alpha_s} y)\otimes(z^*\otimes h)\big) = x^*\otimes \pi(yz^*)h
		\]
		where \(x,y,z\in \mathcal{E}\) and \(h\in H\). 
        Define an operator \(T_s\colon X_s\otimes_{B}(\mathcal{E}^* \otimes_{\mathcal{K}_{B}(\mathcal{E})} H) \to \mathcal{E}^*\otimes_{\mathcal{K}_{B}(\mathcal{E})} H\)
		by \(T_s = (1\otimes V_s)U_s\).
		The operator \(T_s\) is well-defined because of the covariance condition \(V_s\pi(\alpha(T)) = \pi(T)V_s\) for \(T\in \mathcal{K}_B(\mathcal{E})\). Since \(V_s\) is an adjointable isometry, \(T_s\) is so. For \(u \in X_s\), we set 
		\[\
				\phi_{s}(u)(z^* \otimes h)=T_{s}(u \otimes (z^* \otimes h)).\]
                Then, for $x,y,z \in \cle$ and $h \in H$, we have
                \begin{align*}
				\phi_s((x^* \otimes y) \otimes (z^* \otimes h)&:= T_s\big((x^*\otimes y)\otimes (z^*\otimes h) \big)\\
                &= (1\otimes V_s)U_s\big((x^*\otimes y)\otimes (z^*\otimes h)\big)\\
                &= (1\otimes V_s)(x^*\otimes \pi(yz^*)h\\
				&=x^* \otimes V_s\pi(yz^*)h.
		\end{align*}

	For \(y_1,y_2,z_1,z_2\in X_s, e_1,e_2\in \mathcal{E}\) and \(h_1,h_2\in H\), we have 
	\begin{align}\label{equ-cov-phi-cond-1}
		\bigl \langle e_1^*\otimes h_1|\phi_s(y_1^*\otimes z_1)^*\phi_s(y_2^*\otimes z_2)(e^*_2\otimes h_2 \bigr \rangle &=	\bigl \langle \phi_s(y_1^*\otimes z_1)(e_1^*\otimes h_1) | \phi_s(y_2^*\otimes z_2)(e_2^*\otimes h_2) \bigr \rangle \nonumber\\ &= \bigl \langle y_1^*\otimes V_s\pi(z_1e_1^*)h_1 | y_2^*\otimes V_s\pi(z_2e_2^*) h_2\bigr \rangle \nonumber\\
		&= \bigl \langle V_s\pi(z_1e_1^*)h_1 | \pi(y_1y_2^*)V_s\pi(z_2e^*_2)h_2 \bigr \rangle \nonumber\\
		&= \bigl \langle V_s\pi(z_1e_1^*)h_1 |  V_s\pi(\alpha_s(y_1y^*_2))\pi(z_2e^*_2)h_2 \bigr \rangle \nonumber\\
		& =  \bigl \langle \pi(z_1e_1^*)h_1 |  \pi(\alpha_s(y_1y^*_2))\pi(z_2e^*_2)h_2 \bigr \rangle.
	\end{align}
	 The fourth equality above follows from the covariance condition of \((\pi, V)\), and the fifth one follows as \(V_s\) is an isometry. Again we have
	\begin{align}\label{equ-cov-phi-cond-2}
		\bigl \langle e_1^*\otimes h_1|\phi_e\big(\langle y_1^*\otimes z_1| y_2^*\otimes z_2\rangle \big) (e^*_2\otimes h_2) \bigr \rangle &= 	\bigl \langle e_1^*\otimes h_1|\phi_e\big(\langle  z_1| \alpha_s(y_1y_2^*)z_2\rangle \big) (e^*_2\otimes h_2) \bigr \rangle  \nonumber\\
		&= 	\bigl \langle e_1^*\otimes h_1| z_1^*\otimes \pi (\alpha_s(y_1y^*_2)z_2e_2^*)h_2
		\bigr \rangle  \nonumber\\
		&= \bigl \langle h_1| \pi(e_1z_1^*) \pi (\alpha_s(y_1y^*_2)z_2e_2^*)h_2 \bigr \rangle  \nonumber\\
		&= \bigl \langle \pi(z_1e_1^*)h_1 |  \pi(\alpha_s(y_1y^*_2))\pi(z_2e^*_2)h_2 \bigr \rangle .
	\end{align}
	Eq.~\ref{equ-cov-phi-cond-1} and Eq.~\ref{equ-cov-phi-cond-2} give us \(\phi_s(y_1^*\otimes z_1)^*\phi_s(y_2^*\otimes z_2) = \phi_e\big(\langle y_1^*\otimes z_1| y_2^*\otimes z_2\rangle \big)\).

	Let \(s,t\in P\). For \(y_1,y_2,z_1,z_2, \in \mathcal{E}\), 
    we have 
	\begin{align*}
		\phi_t(y_1^*\otimes z_1)\phi_s(y_2^*\otimes z_2)(e^*\otimes h) &= \phi_t(y_1^*\otimes z_1)(y^*_2 \otimes V_s\pi(z_2e^*)h) \\
		&= y^*_1\otimes V_t\pi(z_1y_2^*)V_s\pi(z_2e^*)h\\
		&= y_1^*\otimes V_tV_s\pi(\alpha_s(z_1y_2^*))\pi(z_2e^*)h\\
		&= y_1^*\otimes V_{st}\pi(\alpha_s(z_1y_2^*)z_2e^*)h\\
		&= \phi_{st}(y_1^*\otimes \alpha_s(z_1y^*_2)z_2)(e^*\otimes h) \\
		&= \phi_{st}\big((y_1^*\otimes z_1)\cdot (y_2^*\otimes z_2) \big)(e^*\otimes h)
	\end{align*}
	for all \(e\in \mathcal{E}\) and \(h\in H\). Therefore,  \(\phi^{(\pi,V)}:=\phi=\{\phi_s\}_{s\in Q}\) is a representation of \(X\) on \(\mathcal{E}^*\otimes_{\mathcal{K}_B(\mathcal{E})} H\).
  \hfill $\Box$

Strictly speaking, we do not need the following result later.  However, we have included it as the following result along with Prop. \ref{equivalence of representation theory1} give us the first indication that $C_{red}^{*}(X)$ and $\clk_B(\cle)\rtimes_{red} P$ are Morita equivalence. 
\begin{ppsn}
\label{equivalence of representation theory}
Let $D$ be a $C^{*}$-algebra, and let $L$ be a Hilbert $D$-module. 
    Let $\phi=\{\phi_s\}_{s \in Q}$ be a non-degenerate representation of the product system $X$ on  $L$. Then, $L$ carries a left $B$-action via $\phi_e$. Set 
\[
H:=\cle \otimes_B  L=\cle \otimes_{\phi_e} L.\]
 For $T \in \clk_B(\cle)$, let \[\pi(T):=T \otimes 1.\] For every $s \in P$, there exists a unique adjointable isometry $V_s$ on $H$ such that 
\begin{equation}
\label{defn of Vs}
V_s(\alpha_s(xy^{*})z \otimes h)=x \otimes \phi_s(y^{*} \otimes z)h\end{equation}
for $x,y,z \in \cle$ and $h \in L$. Moreover, the pair $(\pi^{\phi},V^{\phi}):=(\pi,V)$ is a non-degenerate covariant representation of the semigroup dynamical system $(\clk_B(\cle),P,\alpha)$. 
\end{ppsn}
\textit{Proof.} 
Note that $\pi$ extends to a locally strictly continuous representation of the multiplier algebra  $M(\clk_B(\cle))=\cll_B(\cle)$ given by the map $\cll_B(\cle) \ni T \mapsto T \otimes 1 \in \cll_D(\cle \otimes_{B} L)$. It is clear that $\pi$ is non-degenerate.

Let $s \in P$. Let $T_{s}:X_s \otimes_B L \to L$ be defined by \[
	T_s(u \otimes \xi)=\phi_s(u)\xi.
	\]  Since \(T_s\) is an adjointable map and commutes with the left action of \(B\), it follows that the operator
	$
	1\otimes T_s\colon \mathcal{E}\otimes_B(X_s\otimes_B L) \to \mathcal{E}\otimes_{B} L
	$
	is a well-defined adjointable operator, and is also an isometry. Let $\sigma_s:\mathcal{E}\otimes_B X_s \to \cle$ be the unitary map defined by  \[\sigma_s(x\otimes (y^*\otimes z))= \alpha_s(xy^*)z\]
    for $x,y,z \in \cle$ and $h \in L$ (see Remark~\ref{absorption}). 
	 Set $V_{s}:=(1 \otimes T_s)(\sigma_s \otimes 1)^{*}$. Then, $V_s$ is an adjointable isometry. Note that for $x,y,z \in \cle$ and $h \in L$,
	\begin{align*}
     V_s(\alpha_s(xy^*)z\otimes h )&=
	(1\otimes T_s)(\sigma_s \otimes 1)^*(\alpha_s(xy^*)z\otimes h )\\
    & = (1\otimes T_s)\bigl(x\otimes(y^*\otimes z) \otimes h\bigr)\\
    &= x \otimes \phi_s(y^*\otimes z)h.
    	\end{align*}
    Thus, $V_{s}$ satisfies Eq.~\ref{defn of Vs}. Since $\{\alpha_s(xy^*)z \otimes h: x,y,z \in \cle, h \in L\}$ is total, it follows that $V_s$ is uniquely determined by Eq.~\ref{defn of Vs}. 

    Let $s,t \in P$. 
	To prove \(V_sV_t =V_{ts}\), let $x,y,z \in \cle$ and $h \in L$ be given. 
		Since $\alpha_s$ is non-degenerate, we can choose a net \((\alpha_s(z_\beta y_\beta^{*\prime}) z_\beta^{\prime})_{\beta}\) such that \( \alpha_s(z_\beta y_\beta^{\prime *}) z_\beta^{\prime} \to  z\). Now, 
	\begin{align*}
    V_{st}(\alpha_{st}(xy^*)z \otimes h)&=		x\otimes \phi_{st}(y^*\otimes z)h\\
    &=    \lim_{ \beta} x \otimes \phi_{st}\big(y^*\otimes  \alpha_s(z_\beta y_\beta^{\prime *}) z_\beta^{\prime} \big)h\\
		&= \lim_{ \beta} x\otimes  \phi_{st}\big((y^*\otimes z_\beta)\cdot (y_{\beta}^{\prime *}\otimes z_{\beta}^{\prime})\big)h \\
		&= \lim_{ \beta} x\otimes  \phi_t(y^*\otimes z_{\beta}) \phi_{s}(y_{\beta}^{\prime *} \otimes z^{\prime}_{\beta})h \\
		&= \lim_{ \beta} V_t\big( \alpha_t(xy^*) z_{\beta} \otimes \phi_s(y^{\prime *}_{\beta} \otimes z^{\prime}_{\beta}) h\big) \\
		&=\lim_{ \beta} V_t\big( V_s\big (\alpha_{s}(\alpha_t(xy^*) z_{\beta} y^{\prime *}_{\beta}) z^{\prime}_{\beta}  \otimes h \big) \big) \\
		&= \lim_{\beta} V_tV_s \big( \alpha_{st} (xy^*) \alpha_s(z_{\beta}y^{\prime *}_{\beta}) z^{\prime}_{\beta} \otimes h \big) \\
		&= V_tV_s\big( \alpha_{st} (xy^*)z\otimes h\big).
	\end{align*}
	As  \(\{\alpha_{st}(xy^*)z\otimes h : x,y, z\in \mathcal{E}, h\in L, \textup{ and } s,t\in P \}\) is a total subset of \(H\), we can conclude that \(V_{st} = V_tV_s\) for \(s, t\in P\).
	
	Let \(T\in \mathcal{K}_B(\mathcal{E})\) and $t \in P$. Then 
	\begin{equation}\label{eq-cov-cond-1}
		\pi(T)V_t(\alpha_t(xy^*)z\otimes h) = \pi(T) (x\otimes \phi_t(y^*\otimes z)h) = Tx \otimes \phi_t(y^*\otimes z)h
	\end{equation}
	for \(x,y,z\in \mathcal{E}\) and \(h\in L\).
	Again,  for $x,y,z \in \cle$, $h \in L$, we have
	\begin{align}
		V_t\big( \pi(\alpha_t(T)) ( \alpha_t(xy^*)z\otimes h)\big)& = V_t\big( \alpha_t(T)  \alpha_t(xy^*)z\otimes h\big) \nonumber \\ 
		&= V_t\big( \alpha_t(Tx y^*)z\otimes h\big) \nonumber\\
		\label{eq-cov-cond-2}&= Tx\otimes \phi_t(y^*\otimes z)h.
	\end{align}
	As the set $\{\alpha_t(xy^*)z \otimes h: x,y,z \in \cle, h \in L\}$ is total, Eq.~\ref{eq-cov-cond-1} and Eq.~\ref{eq-cov-cond-2} give us the covariance condition \(\pi(T)V_t = V_t\pi(\alpha_t(T))\).
	Therefore, \((\pi, V)\) is a covariant representation of the semigroup dynamical system~\((\clk_B(\cle),P,\alpha)\).  This completes the proof.
\hfill $\Box$

    \begin{rmrk}
We have the following. 
\begin{enumerate}
    \item The `maps' $(\pi,V) \to \phi^{(\pi,V)}$ and $\phi \to (\pi^\phi,V^\phi)$ are inverses of each other under the natural identifications $\cle^{*}\otimes_{\clk_B(\cle)} (\cle \otimes_{\phi_e} L) \cong L$ and $\cle \otimes_{B} (\cle^{*}\otimes_{\pi} H) \cong H$.
    
   \item  Let $(\pi,V)$ be a covariant representation of $(\clk_B(\cle),P,\alpha)$, and let $\phi:=\phi^{(\pi,V)}$. Then, 
	the representation \(\phi = \{\phi_s\}_{s\in Q}\) of \(X\) is Cuntz-Pimsner covariant as defined by Fowler~(\cite{Fowler_2002}) if and only if  \(V_s\) is a unitary operator for  \(s\in P\).
    \end{enumerate}
    \end{rmrk}

    We denote \(\mathcal{K}_B(\mathcal{E})\rtimes_{red}P\) by \(\mathcal{W}\) for the rest of this section. We let $(\pi,V)$ be the regular representation of $(\clk_B(\cle),P,\alpha)$. Note that $\pi:\clk_B(\cle) \to \mathcal{W}$ is one-one, and $\overline{\pi(\clk_B(\cle))\mathcal{W}}=\mathcal{W}$ (Prop. \ref{non-degeneracy-faithful-reduced-basic}). Thus, we view $\clk_{B}(\cle)$ as a $*$-subalgebra of $\mathcal{W}$. Let \[M:=\cle^{*}\otimes_{\clk_B(\cle)}\mathcal{W}=\cle^{*}\otimes_{\pi}\mathcal{W}.\]
    Then, $M$ is a Hilbert $\mathcal{W}$-module. Since $\overline{\pi(\clk_B(\cle))\mathcal{W}}=\mathcal{W}$, $M$ is full. We show that $M$ carries a left action of $C_{red}^{*}(X)$, and then show that $M$ is an imprimitivity bimodule. 
   
\begin{lmma}
\label{different_representation}
Let $\widetilde{\pi}:\clk_B(\cle) \to \cll_B(\ell^2(P) \otimes \cle)$ be defined by 
\[
\widetilde{\pi}(x)(\delta_s \otimes \xi)=\delta_s \otimes \alpha_s(x)\xi.
\]
For $s \in P$, let $\widetilde{V}_s$ be the adjointable operator on $\ell^2(P) \otimes \cle$ defined by 
\[
\widetilde{V}_s(\delta_t \otimes \xi)=\delta_{ts} \otimes \xi.\] Let $(\pi,V)$ be the regular representation of the dynamical system $(\clk_B(\cle),P,\alpha)$. Then, there exists a faithful homomorphism $\lambda: \clk_B(\cle) \rtimes_{red} P \to \cll_B(\ell^2(P) \otimes \cle)$ such that 
\[
\lambda(\pi(x)V_s)=\widetilde{\pi}(x)\widetilde{V}_s
\]
for $x \in \clk_B(\cle)$ and $s \in P$. 
\end{lmma}
\textit{Proof.} Note that the map $U:(\ell^2(P) \otimes \clk_B(\cle))\otimes_{\clk_B(\cle)} \cle  \to \ell^2(P) \otimes \cle$ defined by \[
U(\delta_s \otimes T) \otimes \xi:=\delta_s \otimes T\xi\]
is a unitary operator. Here, the left action of $\clk_B(\cle)$ on $\cle$ in the interior tensor product $(\ell^2(P) \otimes \clk_B(\cle))\otimes_{\clk_B(\cle)}\cle$ is given by the identity representation.  Observe that, for $s \in P$ and $x \in \clk_B(\cle)$,   \[U(\widetilde{V}_s \otimes 1)U^{*}=V_s~; \quad U(\widetilde{\pi}(x) \otimes 1)U^{*}=\pi(x).\] Let $\lambda:\clk_B \rtimes_{red} P \to \cll_B(\ell^2(P)\otimes \cle)$ be the homomorphism defined by $\lambda(\cdot)=U(\cdot)U^*$. Then, $\lambda$ is the required map.  \hfill $\Box$

	Let $\lambda:\mathcal{W} \to \mathcal{L}_B(\ell^2(P) \otimes \mathcal{E})$ be the homomorphism given by Lemma \ref{different_representation}. Since $\lambda$
	  is faithful, we have a faithful homomorphism (see Remark \ref{lem:comp-faithful}) 
	\[\Phi_0:\cll_{\mathcal{W}}(\cle^* \otimes_{\clk_B(\cle)} \mathcal{W}) \to \cll_B((\cle^*\otimes_{\clk_B(\cle)} \mathcal{W})\otimes_{ \mathcal{W}}(\ell^2(P)\otimes \mathcal{W})=\cll_B((\cle^*\otimes_{\clk_B(\cle)}\mathcal{W})\otimes_{\lambda}(\ell^2(P)\otimes \mathcal{W})).\]
		defined by \(\Phi_0(T)=T \otimes 1\).
	Now, we have the following isomorphism of Hilbert modules
	\begin{align*}
		M \otimes_{\mathcal{W}}(\ell^2(P)\otimes \cle)& \cong (\mathcal{E}^*\otimes_{\mathcal{K}_B(\mathcal{E})}\mathcal{W})\otimes_{\mathcal{W}} (\ell^2(P)\otimes \mathcal{E}) \\
        & \cong \mathcal{E}^*\otimes_{\mathcal{K}_B(\mathcal{E})} (\ell^2(P)\otimes \mathcal{E}) \\
		&\cong \mathcal{E}^*\otimes_{\mathcal{K}_B(\mathcal{E})} \bigl( \bigoplus_{s\in P} \mathcal{E}\bigr) \\
        & \cong \bigoplus_{s\in P} \mathcal{E}^*\otimes_{\alpha_s} \mathcal{E}\\
        & \cong \bigoplus_{s\in P} X_s.
	\end{align*}
	The isomorphism of the second line above is given by \((x^*\otimes w)\otimes y \mapsto x^*\otimes \lambda(w)\cdot y\), the isomorphism of the fourth line is given by \(x^*\otimes (u \otimes \delta_s) \mapsto (x\otimes_{\alpha_s} u) \otimes \delta_s\). The resulting unitary from 
    $M \otimes_{\mathcal{W}}(\ell^2(P)\otimes \cle)=(\cle^* \otimes_{\clk_B(\cle)}\mathcal{W})\otimes_{\mathcal{W}}(\ell^2(P)\otimes \mathcal{W}) \to \bigoplus_{s \in P}X_s$ is denoted $J$. Define a $*$-homomorphism $\Phi:\cll_{\mathcal{W}}(\cle^* \otimes_{\clk_B(\cle)}\mathcal{W}) \to \cll_{B}(\bigoplus_{s \in P}X_s)$ by \[\Phi(T)=J(\Phi_0(T))J^*=J(T \otimes 1)J^*.\]
	\begin{ppsn}\label{lem:Mor-equ-left-act-2}
		 There exists a unique faithful $*$-homomorphism $\phi:C_{red}^{*}(X) \to \cll_{\mathcal{W}}(M)$ such that $\Phi \circ \phi=\mathcal{F}$, where $\mathcal{F}:C_{red}^{*}(X) \to \cll_B(\bigoplus_{s \in P}X_s)$ is the inclusion. With the left action given by $\phi$,  the Hilbert module $M=\mathcal{E}^*\otimes_{\mathcal{K}_B(\mathcal{E})} \mathcal{W}$ is a \(C^*\)-correspondence from \(C^*_{red}(X)\) to \(\mathcal{W}\).
	\end{ppsn}
	\textit{Proof.}
		Let \(C^*(X)\) denote the universal Toeplitz algebra of the product system \(X\), i.e. the universal $C^{*}$-algebra generated by $\{u: u \in X_s, s \in Q\}$ such that for $s \in P$, $t \in P$, $u \in X_s$ and $v \in X_t$,     $ uv=u\cdot v$ and $u^{*}v=\langle u|v \rangle$. Here, $\cdot$ denotes the multiplication in $X$. 
    
    Let $(\pi,V)$ be the regular representation of $(\clk_B(\cle), P, \alpha)$. For $x \in \clk_B(\cle)$, define the operator $\widehat{\pi}(x):\mathcal{W} \to \mathcal{W}$ by 
    \[
    \widehat{\pi}(x)y=xy.\]
 For $s \in P$, let  $\widehat{V}_s:\mathcal{W} \to \mathcal{W}$ be defined by \[\widehat{V}_{s}(y)=V_sy.\] As operators on $\mathcal{W}$ which is a Hilbert $\mathcal{W}$-module, $\widehat{\pi}(x)$ and $\widehat{V}_s$ are adjointable  for $x \in \clk_B(\cle)$ and $s \in P$. Moreover, $(\widehat{\pi},\widehat{V})$ is a non-degenerate covariant representation of $(\clk_B(\cle),P,\alpha)$. 
    We apply Prop.~\ref{equivalence of representation theory1} to \((\widehat{\pi},\widehat{V})\) to  obtain a representation \(\phi:=\phi^{(\widehat{\pi},\widehat{V})}=\{\phi_s\}_{s \in Q}\) of the product system $X$ on $\cle^{*}\otimes_{\clk_B(\cle)}\mathcal{W}$.

     Let \(\eta_1\colon C^*(X) \to \mathcal{L}_{\mathcal{W}}(\mathcal{E}^*\otimes_{\mathcal{K}_B(\mathcal{E})}\mathcal{W})\) be the homomorphism such that $\eta_1(x)=\phi_s(x)$ for $x \in X_s$ and for $s \in P$. The homomorphism $\eta_1$ exists by the universal property.  Now consider the following diagram, 
		\begin{figure*}[hbt]
		\begin{center}
			\begin{tikzpicture}[scale=0.5,cap=round]
				
				\draw [->] (-6.3, 3) -- (-1, 3);
				\draw [->] (-8, 2.2) -- (-8,-1);
				\draw [->] (-6.2,-2)--(-.2,-2);
				\draw [->] (2, 2)--(2,-1.2);
				
				\node (A) at (-8, 3) {\(C^*(X)\)};
				\node (B) at (3,3) {\(\mathcal{L}_{W}(\mathcal{E}^*\otimes_{\mathcal{K}_B(\mathcal{E})} W )\)};
				\node (C) at (-8, -2) {\(C^*_{red}(X)\)};
				\node (D) at (3, -2) {\(\mathcal{L}_{B}\big(\bigoplus_{s\in P}X_s\big)\)};
				\node (E) at  (-4, 3.8) {\(\eta_1\)};
				\node (F) at  (-8.8, .4) {\(q\)};
				\node (G) at  (-3, -2.5) {\(\mathcal{F}\)};
				\node (H) at  (3, .4) {\(\Phi\)};
			\end{tikzpicture}
		\end{center}
	\end{figure*}
	where~\(q\) is the quotient map and \(\mathcal{F}\) is the inclusion. 
    
\textit{Claim:} The above diagram is commutative.  
It is enough to check the equality 
\[
J\Phi_0(\phi_s(x^*_1\otimes y_1))J^* = \mathcal{F} q(x_1^*\otimes y_1)
\] 
for \(x^*_1\otimes y_1 \in X_s\) and \(s\in P\). Let~\((x^*\otimes_{\alpha_t}\alpha_t(z_1z_2^*) y)\otimes \delta_t\in \bigoplus_{s\in P} X_s \). The unitary \(J^*\) maps 
	\[
	\big(x^*\otimes_{\alpha_t}\alpha_t(z_1z_2^*)y\big)\otimes \delta_t \mapsto x^*\otimes \big(\alpha_t(z_1z_2^*)y\otimes \delta_t\big) \mapsto \big(x^*\otimes \pi(z_1z^*_2)\big)\otimes (y\otimes \delta_t).
	\]
	Now we have 
	\begin{align*}
		J\Phi_0(\phi_s(x^*_1\otimes y_1))J^* \bigl(	(x^*\otimes_{\alpha_t} \alpha_t(z_1z_2^*) y)\otimes \delta_t\bigr) &= 	J\Phi_0(\phi_s(x^*_1\otimes y_1))\bigl((x^*\otimes \pi(z_1z_2^*))\otimes (y\otimes \delta_t)\bigr)\\
		&= 	J\bigl(\phi_s(x_1^*\otimes y_1)(x^*\otimes \pi(z_1z_2^*))\otimes (y\otimes \delta_t)\bigr)\\
        &= J\bigl(x_1^*\otimes V_s\pi(y_1x^*)\pi(z_1z_2^*)) \otimes (y\otimes \delta_{t})\bigr)\\
        &=(x_1^* \otimes_{\alpha_t}\alpha_t(y_1x^*)\alpha_t(z_1z_2^*)y\otimes \delta_{ts}\\
        			&= \mathcal{F}q(x^*_1\otimes y_1) \bigl((x^*\otimes_{\alpha_t}\alpha_t(z_1z_2^*)y) \otimes \delta_t\bigr).
	\end{align*}
	This proves the claim.

Since the diagram is commutative, $\Phi$ is faithful and $q$ is onto, there exists a homomorphism \(\widetilde{\phi}\colon C^*_{red}(X) \to \mathcal{L}_{W}(\mathcal{E}^*\otimes_{\mathcal{K}_B(\mathcal{E})}W)\) such that \(\Phi \circ \widetilde{\phi} = \mathcal{F}\). The faithfulness of $\widetilde{\phi}$ follows from that of $\Phi$ and $\mathcal{F}$. We abuse notation, and we denote $\widetilde{\phi}$ by $\phi$.  Thus, $M=\cle^{*}\otimes_{\clk_B(\cle)}\mathcal{W}$ is a $C^{*}$-correspondence from $C_{red}^{*}(X)$ to $\mathcal{W}$, where the left action of $C_{red}^{*}(X)$ on $M$ is  given by $\phi$. 
    \hfill $\Box$

To show that $M$ is an imprimitivity bimodule, we need the following lemma.
    \begin{lmma}
\label{check only for dense}
    Let $\mathcal{A}$ be a $C^{*}$-algebra, and let $\mathcal{B}$ be a $C^{*}$-subalgebra of $\mathcal{A}$. Suppose that $F$ is a full Hilbert $\mathcal{B}$-module. Let $F_0:=F \otimes_{\mathcal{B}} \cla$. Assume that $\overline{\clb \cla}=\cla$. Suppose $\mathcal{D} \subset \cla$ is such that $\clb \subset \cld$ and the $C^{*}$-algebra generated by $\cld$ is $\cla$. Then, the $C^{*}$-algebra generated by $\{\theta_{e \otimes d,f \otimes b}: e,f \in F, d \in \cld, b \in \clb\}$ is $\clk_\cla(F_0)$. 
\end{lmma}
\textit{Proof.} For $e_1,e_2 \in F$ and $a_1,a_2 \in \cla$, the compact operator $(e_1 \otimes a_1)(e_2 \otimes a_2)^{*}$ on $F \otimes_{\clb}\cla$ will be denoted $\theta_{e_1 \otimes a_1,e_2 \otimes a_2}$. Note that since $\overline{\clb \cla}=\cla$, $\clb$ contains an approximate identity of~$\cla$. Let $\mathcal{C}$ be the $C^{*}$-algebra generated by $\{\theta_{e\otimes d,f \otimes b}: e,f \in F, d \in \cld, b \in \clb\}$. 

For $d_1,d_2 \in \cld$, $e_1,e_2,f_1,f_2 \in F$ and $b_1,b_2 \in \clb$, note that
\begin{equation}
\label{base for induction zero}
\mathcal{C} \ni \theta_{e_1 \otimes d_1,f_1 \otimes b_1}\theta_{e_2 \otimes d_2,f_2 \otimes b_2}^{*}=
\theta_{e_1 \otimes d_1b_1^*\langle f_1|f_2 \rangle b_2, e_2 \otimes d_2 }.
\end{equation}
The above computation along with the fact that $F$ is full and $\clb$ contains an approximate identity of $\cla$ imply that $\theta_{e_1 \otimes d_1,e_2 \otimes d_2} \in \mathcal{C}$ for every $e_1,e_2 \in F$ and $d_1,d_2 \in \mathcal{D}$. 

Let $e_1,e_2,e_3,e_4 \in F$, $d_1,d_2,d_3,d_4 \in \cld$, note that 
\begin{equation}
\label{base for induction}
\mathcal{C} \ni \theta_{e_1 \otimes d_1,e_2 \otimes d_2}\theta_{e_3 \otimes d_3,e_4 \otimes d_4}=\theta_{e_1 \otimes d_1d_2^*\langle e_2|e_3 \rangle d_3,e_4 \otimes d_4}.\end{equation}
Again using the fact that $\cld$ contains an approximate identity of $\cla$ and the fact that $F$ is full, we see that $\theta_{e_1 \otimes d_1^*,e_2 \otimes d_2} \in \mathcal{D}$ for $e_1,e_2 \in F$ and $d_1,d_2 \in \cld$. Similar arguments imply that $\theta_{e_1 \otimes d_1,e_2 \otimes d_2} \in \mathcal{C}$ for every $e_1,e_2 \in F$ and $d_1,d_2 \in \mathcal{D} \cup \mathcal{D}^{*}$. Thus, with no loss of generality, we can suppose that $\mathcal{D}=\mathcal{D}^{*}$

Using induction and by similar arguments that we used to arrive at Eq.~\ref{base for induction} and by repeatedly appealing to the fact that $\clb$ (and $\cld$) contains an approximate identity of $\cla$ and $F$ is full, we can prove that if $w$ is any word in $\mathcal{D}$ and for $e_1,e_2 \in F$ and $b \in \clb$, $\theta_{e_1 \otimes w,e_2 \otimes b} \in \mathcal{C}$. Since $\mathcal{C}$ is $*$-closed, if $w^{'}$ is a word in $\mathcal{D}$, then for $e_1,e_2 \in F$ and $b \in \clb$, $\theta_{e_1 \otimes b,e_2 \otimes w^{'}} \in \mathcal{C}$. Now, let $e_1,e_2 \in F$, and let $w_1,w_2$ be words in $\cld$. Let $f_1,f_2 \in F$ be given. Let $(u_i)_i$ be a net in $\clb$ which is an approximate identity for $\cla$. Then, 
\[
\theta_{e_1 \otimes w_1\langle f_1|f_2 \rangle,e_2 \otimes w_2}=\lim_{i}\theta_{e_1 \otimes w_1,f_1 \otimes u_i}\theta_{f_2 \otimes u_i,e_2 \otimes w_2} \in \mathcal{C}.
\]
Since $F$ is full, $\theta_{e_1 \otimes w_2,e_2 \otimes w_2} \in \mathcal{C}$. As $\mathcal{D}$ generates $\cla$, $\mathcal{C}=\mathcal{K}_\cla(F_0)$. This completes the proof. \hfill $\Box$

\medskip 

\textit{Proof of Thm.~\ref{main_Muhly}.} Let $M$ be the  \(C^*\)-correspondence from $C_{red}^{*}(X)$ to $\mathcal{W}$ considered in Prop.~\ref{lem:Mor-equ-left-act-2}. 
Let $s \in P$, $e_1,e_2 \in \cle$, and let $T,S \in \clk_B(\cle)$ be given. We claim that 
\begin{equation}
\label{compact and reduced coincide}
\theta_{e_1^* \otimes V_s\pi(T),e_2^{*}\otimes \pi(S)}=\phi_s(e_1^* \otimes_{\alpha_s} TS^*e_2).
\end{equation}
Let $e_3 \in \cle$ and $w \in \mathcal{W}$ be given. Calculate as follows to observe that 
\begin{align*}
\phi_s(e_1^* \otimes_{\alpha_s}TS^*e_2)(e_3^* \otimes w)&=e_1^{*}\otimes \widehat{V}_s\widehat{\pi}(TS^*e_2e_3^*)w=e_1^{*}\otimes V_s\pi(TS^*)\pi(e_2e_3^*)w\\
&=e_1^{*}\otimes V_s\pi(T) \langle e_2 ^* \otimes \pi(S)|e_3^* \otimes w\rangle_{\mathcal{W}}\\
&=\theta_{e_1^* \otimes V_s\pi(T),e_2^* \otimes \pi(S)}(e_3^* \otimes w).
\end{align*}
This proves the claim. Eq.~\ref{compact and reduced coincide} implies that $\phi(C_{red}(X)) \subset \clk_{\mathcal{W}}(M)$. 
From Eq.~\ref{compact and reduced coincide} and Lemma~\ref{check only for dense}, we get 
$\clk_{\mathcal{W}}(M)\subset \phi(C_{red}^{*}(X)$. Hence, $\phi(C_{red}^{*}(X))=\clk_\mathcal{W}(M)$. As $\phi$ is faithful, $M$ is a $C_{red}^{*}(X)$-$\mathcal{W}$ imprimitivity bimodule. \hfill $\Box$

\begin{rmrk}
The case when $P$ is a group deserves special mention. 
In this case, Thm.~\ref{main_Muhly}, applied to a group, is essentially the Packer-Raeburn stabilization theorem. This theorem states that a twisted crossed product can be written as an ordinary crossed product up to a  Morita equivalence.  
 For simplicity, let us consider the case of product systems of Hilbert spaces where the fibres are one-dimensional. 
 
In this case, a product system is determined by a $2$-cocycle. 
 Let $G$ be a group, and let $\omega:G \times G \to \bbt$ be a $2$-cocycle, i.e. 
\[
\omega(r,s)\omega(rs,t)=\omega(r,st)\omega(s,t)
\]
for every $r,s,t \in G$.
For $s \in G$, let $X_s:=\bbc$, and denote the unit vector $1$ in $X_s$ by $e_s$. Then, $\{X_s\}_{s \in G}$ is a product system of Hilbert spaces over $G$ with the multiplication given by 
\[
e_se_t=\omega(s,t)e_{ts}.
\]
Denote the resulting product system by $X^{\omega}$.  Note that $C_{red}^{*}(X^\omega)$ is the reduced twisted group $C^{*}$-algebra $C_{red}^{*}(G,\omega)$ as defined in \cite[Chapter 5]{Dana-Williams}.

Let $\pi:G \to B(H)$ be a $\omega$-projective unitary representation of $G^{op}$ on a Hilbert space $H$, i.e. for $s \in G$, $\pi(s)$ is a unitary operator, and for $s,t \in G$, 
\[
\pi(s)\pi(t)=\omega(s,t)\pi(ts).
\]
For $s \in G$, let $\alpha_s:\clk(H) \to \clk(H)$ be defined by $\alpha_s=Ad(\pi(s))$. Then, $\alpha_{\pi}:=\{\alpha_s\}_{s \in G}$ is an $E_0$-semigroup over $G^{op}$ on $H$. Then, the product system associated with $\alpha$ is $X^{\omega}$, and it follows from Thm.~\ref{main_Muhly} that $C_{red}^{*}(G,\omega)$ is Morita equivalent to $\clk(H) \rtimes_{\alpha, red} G^{op}$. 
\end{rmrk}

\section{Do product systems come from  $E_0$-semigroups?}
 \label{sec-exist-E0-semigroup} 

The main advantage of our Morita equivalence result (Thm.~\ref{main_Muhly}) is that while it is not known whether $C_{red}^{*}(X)$ has a groupoid crossed producd/Fell bundle presentation (except for the case of quasi-lattice ordered semigroups~(\cite{Rennie_Sims})), the reduced crossed product $\clk_B(\cle) \rtimes_{red} P$ has a groupoid crossed product presentation~(\cite{Sundar_Khoshkam}). However, it is necessary to settle the following fundamental question before we can apply Thm.~\ref{main_Muhly}. 

\textbf{(Q):} Is it true that every product system is isomorphic to the product system of an $E_0$-semigroup?
 
     Arveson first settled the above question in the setting of Hilbert spaces and for the topological semigroup $(0,\infty)$. It is interesting to note here that Arveson's original proof (\cite{Arv_Fock2}, \cite{Arv_Fock3}, \cite{Arv_Fock4}) relied on a very deep analysis of the reduced $C^{*}$-algebra of a product system over $(0,\infty)$, called the spectral $C^{*}$-algebra by Arveson, and we are trying to argue that settling \textbf{(Q)} in the affirmative helps us to better understand the reduced $C^{*}$-algebra.  Skeide gave an alternate proof via an induced construction~(\cite{skeide}), and also  
       settled the module version~(\cite{Skeide2}) through the same trick for $\bbr_{+}$. 
        The Hilbert space version beyond the case of $1$-parameter semigroups was considered in the works of the second author~(\cite{Murugan_Sundar_discrete_mathsci}, 
\cite{Murugan_Sundar_continuous}, \cite{Sundar_Existence}). In particular, the induced construction trick of Skeide was generalised in~\cite{Sundar_Existence}, and it was proved that for a large class of subsemigroups of discrete groups, which includes normal subsemigroups, \textbf{(Q)} has a positive answer if we consider  Hilbert spaces.      We explain below how to adapt these techniques to the Hilbert module setting, and we prove Thm.~\ref{main_existence}, which states that, under reasonable assumptions, every product system arises from an $E_0$-semigroup. 

Until further mention, let $Q$ be a semigroup, and let $B$ be a fixed separable $C^{*}$-algebra. 

\begin{dfn}[{\cite{Skeide_module1}, \cite{Skeide2}}]
Let $X$ be a product system of $B$-$B$-correspondences over $Q$. For $s,t \in Q$, let $U_{s,t}:X_s \otimes_B X_t \to X_{st}$ be the multiplication map given by 
\[
U_{s,t}(u \otimes v)=uv.\]  Let $\cle$ be a full Hilbert $B$-module. Suppose that  for every $s \in Q$, we have a unitary operator $\sigma_s:\cle \otimes_B X_s \to \cle$. Then, $(\cle,\sigma:=\{\sigma_s\}_{s \in Q})$ is called a left dilation  of $X$ if for $s,t \in Q$,
\begin{equation}
    \label{left dilation equation}
    \sigma_t(\sigma_s \otimes 1)=\sigma_{st}(1 \otimes U_{s,t}).
\end{equation}
    
\end{dfn}

Let $X$ be a product system of $B$-$B$-correspondences over $Q$. Suppose  $(\cle,\sigma=\{\sigma_s\}_{s \in Q})$ is a left dilation of $X$. For $s \in Q$, let $\alpha_s:\cll_{B}(\cle) \to \cll_B(\cle)$ be defined by \[\alpha_s(T)=\sigma_s(T \otimes 1)\sigma_s^{*}.\]
Then, Eq.~\ref{left dilation equation} translates to the fact that $\alpha:=\{\alpha_s\}_{s \in P}$ is an $E_0$-semigroup over $Q^{op}$ on $\cll_B(\cle)$, and  the product system associated with $\alpha$ is $X$. We call $\alpha$  \emph{the $E_0$-semigroup associated with the left dilation $(\cle,\sigma)$}.
Conversely, if $\alpha$ is an $E_0$-semigroup over $Q^{op}$ on $\cll_B(\cle)$, and $X$ is the associated product system over $Q$, then $(\cle,\sigma:=\{\sigma_s\}_{s \in Q})$ is a left dilation of the product system $X$, where $\sigma_s:\cle \otimes_B X_s \to \cle$ is given by the formula
\[
\sigma_s(x \otimes (y^* \otimes z))=\alpha_s(xy^*)z.\]
Verifications of the above assertions are routine. We also refer the reader to  \cite[Section 6]{Skeide_module1}. 

Summarising our discussion, we have the following.
\begin{ppsn}
    \label{left dilation versus}
    If $X$ is a product system over $Q$ of $B$-$B$-correspondences, then $X$ is isomorphic to the product system of an $E_0$-semigroup over $Q^{op}$ on $\cll_B(\cle)$ for some full Hilbert $B$-module $\cle$ if and only if $X$ has a left dilation. 
\end{ppsn}

Let $X$ be a product system of $B$-$B$-correspondences over $Q$. For a subsemigroup $R$ of $Q$, we denote the restriction of $X$ to $R$ by $X|_R$.

\begin{thm}
\label{concrete}
Let $Q$ be a subsemigroup of a group $G$ such that   $Q^{-1}Q=G$. Suppose $R \subset Q$ is a semigroup such that $R^{-1}R$ is a group, and $R^{-1}Q=G$. Let $X$ be a product system of $B$-$B$-correspondences over $Q$. If $X|_R$ has a left dilation, then $X$ has a left dilation.
\end{thm}

\begin{rmrk}
The `right dilation' version in the Hilbert space setting is Thm.~2.6 of~\cite{Sundar_Existence}. The `left dilation' version is analogous, and we give minimal details. 
 The trick is to induce a left dilation of $X|_R$ to a left dilation of $X$. This trick for the case $R=\bbn$ and $Q=(0,\infty)$ is originally due to Skeide. 
\end{rmrk}

\textit{Proof of Thm.~\ref{concrete}}. 
	The idea of the proof is the same as in ~\cite[Theorem 2.6]{Sundar_Existence}.
	Suppose \((\mathcal{E}, \sigma:=\{\sigma_s\}_{s\in R})\) is a left dilation for the product system \(X|_{R}\). Let \(K=R^{-1}R\) which is a subgroup of \(G\). The set of right cosets of $K$ is denoted $K\backslash G$. For \([x]\in K\backslash G\), we define 
	\[
	\Delta ([x]) := \coprod_{y \in Q, [y]=[x]} \mathcal{E}\otimes_BX_y.
	\]
	Since \(R^{-1}Q=G\), for \(x\in G\), there exists \(y\in Q\) such that \([y]=[x]\). Thus, \(\Delta ([x])\) is non-empty. For \(z,y\in [x]\cap Q\), \(u\in \mathcal{E}\otimes_{B}X_y\) and \(v\in \mathcal{E}\otimes_{B}X_z\), we define \(u\sim v\) if and only if there are \(a,b\in R\) such that \(ay=bz\) and
	\[
	(1\otimes U_{a,y})(\sigma^*_{a}\otimes 1) u = (1\otimes U_{b,z})(\sigma^*_{b} \otimes 1)v.
	\]
	A similar argument as in \cite[Prop.~2.3]{Sundar_Existence} ensures that `\(\sim\)' is an equivalence relation on \(\Delta([x])\) for \([x]\in K\backslash G\).
	Let \(\mathcal{E}_{[x]} \) be the set of equivalence classes of \(\Delta([x])\). Then, \(\mathcal{E}_{[x]}\) has a natural Hilbert \(B\)-module structure after it is identified with ~\(\mathcal{E}\otimes_B X_y\) for some \(y\in [x]\cap Q\) (the choice of $y$ does not affect the Hilbert module structure).
	
	For $y \in Q$, $t \in Q$ and $\eta \in X_t$, denote the map \[\cle \otimes_{B} X_y \ni u \otimes \xi \mapsto  u \otimes \xi \eta \in \cle \otimes_{B} X_{yt}\] by $\theta(y,t,\eta)$.
	
	For $[x] \in K\backslash G$, $t \in Q$ and $\eta \in X_t$, denote the map $\displaystyle \coprod_{y \in Q, [y]=[x]}\theta(y,t,\eta):\Delta([x]) \to \Delta([xt])$ by $\theta([x],t,\eta)$. Then, $\theta([x],t,\eta)$ preserves the equivalence relation $\sim$ and descends to a map from $\cle_{[x]} \to \cle_{[xt]}$ which we still denote by $\theta([x],t,\eta)$. 
	
        Now, define a new Hilbert \(B\)-module
	\[
	\widetilde{\mathcal{E}} := \bigoplus_{[x]\in K/G} \mathcal{E}_{[x]}.
	\]
    Let $t \in Q$. We claim that there exists an isometry \(\theta_t \colon \widetilde{\mathcal{E}}\otimes_B X_t \to \widetilde{\mathcal{E}}\) such that  \[
	\theta_t(([u\otimes \xi] \otimes \delta_{[x]})\otimes \eta)= [u\otimes \xi\eta] \otimes \delta_{[xt]}=\theta([x],t,\eta)([u\otimes \xi])\otimes \delta_{[xt]} .\]
	For $s \in Q$, $u_1,u_2 \in \cle$, $\xi_1,\xi_2 \in X_{s}$, and $\eta_1,\eta_2$ in $X_t$, we have
	\begin{align*}
				 \bigl \langle ([u_1\otimes \xi_1 \eta_1]\otimes \delta_{[st]} ) | ([u_2\otimes \xi_2\eta_2]\otimes \delta_{[st]}) \bigr \rangle 
		&= \bigl \langle \xi_1\eta_1 | \langle u_1|u_2\rangle \xi_2\eta_2  \bigr \rangle \\
		&= \bigl \langle \eta_1 | \big\langle \xi_1| \langle u_1|u_2\rangle \xi_2\big\rangle \eta_2  \bigr \rangle \\
		&= \bigl \langle \eta_1| \big\langle u_1\otimes \xi_1|u_2\otimes \xi_2\big\rangle \eta_2\bigr \rangle \\ 
		&= \bigl \langle ([u_1\otimes \xi_1]\otimes \delta_{[s]})\otimes \eta_1| ([u_2\otimes \xi_2]\otimes \delta_{[s]})\otimes \eta_2 \bigr \rangle.
	\end{align*}
		Therefore, the map \(\theta_t\) is a well defined isometry.
	
	To prove \(\theta_t\) is an unitary, consider an element in $\widetilde{\cle}$ of the form \([u\otimes x] \otimes \delta_{[g]}\). Since \(R^{-1}Q =G\), we can choose \(r\in R\) such that \([rg] =[g]\) and \(rgt^{-1}\in Q\). We set \(rg = \tilde{g}\). Since the map \(U_{\tilde{g}t^{-1}, t} \colon X_{\tilde{g}t^{-1}} \otimes X_t \to X_{\tilde{g}}\) is a unitary operator, we can choose a net $(z_\beta)_{\beta}$ in $X_{\tilde{g}t^{-1}}\otimes X_t$ of the form \(\displaystyle z_\beta: = \sum_{i=1}^{n_{\beta}} x_i^{\beta}\otimes y_{i}^{\beta}\) such that  \(U_{\tilde{g}t^{-1}, t}(z_\beta) \to x\) as \(\beta \to \infty\). 
	Then, we have 
	\[
	\theta_t\biggl( \lim_{\beta} \sum_{i=1}^{n_{\beta}} \big([u\otimes x^{\beta}_i])\otimes \delta_{[\tilde{g}t^{-1}]} \big) \otimes y^{\beta}_i\biggr) = \lim_{\beta} ([u\otimes U_{\tilde{g}t^{-1}, t}(z_\beta)]\otimes \delta_{ [\tilde{g}]}) = [u\otimes x]\otimes \delta_{[g]}.
	\] 
	Therefore, \(\theta_t\) is unitary for \(t\in Q\). The verification that for $s,t \in P$, 
    \[\theta_t(\theta_s \otimes 1)=\theta_{ts}(1 \otimes U_{s,t})\] is left to the reader as it follows directly from the associativity of the multiplication rule of the product system $X$. 
\hfill $\Box$

\medskip

\textit{Proof of Thm.~\ref{main_existence}}. With the notation of Thm.~\ref{main_existence}, set $Q:=P^{op}$ and $R:=\{a^n: n \geq 0\}$. Then, $Q^{-1}Q=G$ and $R^{-1}Q=G$. If $a$ is of infinite order, then $R\cong \bbn_0$. In this case, we can apply  \cite[Thm. 3.1]{Skeide_non_unital}, Thm.~\ref{concrete} and Prop.~\ref{left dilation versus} to conclude the result. 
    If $a$ is of finite order, then $Q$ is a group. In this case, we apply Thm.~\ref{concrete}  with $R=\{e\}$ and  Prop.~\ref{left dilation equation} to conclude the proof. \hfill $\Box$

\begin{rmrk}
Thm.~\ref{main_existence} is applicable to the following examples.
\begin{enumerate}
    \item Finitely generated subsemigroups of abelian groups.
    \item $\bbn^k \cup \{0\}$. Note that this is not finitely generated but has $(1,1,\cdots,1)$ as an order unit. Recall that $\bbn=\{1,2,\cdots\}$ and $\bbn_0=\{0,1,2,\cdots\}$ as per our notation. 
    \item Let $R$ be an integral domain, and let $\Gamma$ be a finitely generated subsemigroup of $R^{*}$. Define 
    \[
    R \rtimes \Gamma:=\Big\{\begin{bmatrix}
     a & b \\
     0 & 1
        \end{bmatrix}: a \in \Gamma, b \in R\Big\}.\]
    Then, $R \rtimes \Gamma$ is a semigroup with composition given by matrix multiplication. Note that $(R \rtimes \Gamma)^{op}$ is right Ore and has an order unit if $\Gamma$ is finitely generated. So, a product system over $R \rtimes \Gamma$ comes from an $E_0$-semigroup over $(R \rtimes \Gamma)^{op}$.
    \end{enumerate}
There are also many examples of semigroups that satisfy the Ore condition but fail to have an order unit; the simplest example being $\bbn$ with the semigroup law given by multiplication.  We do not know whether Thm.~\ref{main_existence} is true for this example. (However,  see~\cite[Thm. 4.11]{Sundar_Existence} for the Hilbert space version.) Also, the induced construction trick certainly does not adapt well beyond the Ore case, and we do not know whether the conclusion of Thm.~\ref{main_existence} remains valid in the non-Ore case. 
For the case of semigroups that do not embed in groups, an example of a product system of Hilbert spaces that is not isomorphic to the product system of an $E_0$-semigroup was constructed in~\cite[Section 5]{Sundar_Existence}. 
\end{rmrk}

\section{A groupoid presentation}
\label{sec-gpd-presentation}

With Thm.~\ref{main_Muhly} and Thm.~\ref{main_existence} established, in the case of Ore semigroups with an order unit, the study of the reduced $C^{*}$-algebra of product systems of full proper $C^{*}$-correspondences is essentially the same as the study of the reduced crossed product of semigroup dynamical systems. In this section, we recall from \cite{Sundar_Khoshkam} the groupoid crossed product presentation of the reduced $C^{*}$-algebra of a semigroup dynamical system. We use this picture in the next two sections when we discuss the nuclearity and exactness of the reduced crossed product and the invariance of $K$-theory under homotopy.

The following proposition that generalises \cite[Lemma 3.5]{Sundar_Khoshkam}  allows us to restrict ourselves to the unital case. A semigroup dynamical system $(A,P,\alpha)$ is said to be \emph{unital} if $A$ is unital, and $\alpha_s$ is unital for every $s \in P$.  We need the generalised version to study the invariance of the $K$-theory of the reduced crossed product under homotopy.

\begin{ppsn}
\label{split_exact_homotopy}
Let $(A,P,\alpha)$ and $(B,P,\beta)$ be two unital semigroup dynamical  systems, and let $\phi:A \to B$ be a unital homomorphism which is $P$-equivariant, i.e. \[\phi(\alpha_s(a))=\beta_s(\phi(a))\] for $s \in P$ and $a \in A$. Suppose that the homomorphism $\phi$ admits a $P$-equivariant splitting~$\epsilon$, i.e. $\epsilon:B \to A$ is a $P$-equivariant homomorphism such that $\phi \circ \epsilon=id$. Let $I:=Ker( \phi)$. Suppose that $\overline{\alpha_s(I)I}=I$ for each $s \in P$. Then, we have a short exact sequence 
\[
0 \longrightarrow I \rtimes_{red} P \longrightarrow A \rtimes_{red} P \longrightarrow B \rtimes_{red} P \longrightarrow 0
\]
  which is also split-exact.
\end{ppsn}
\textit{Proof.} Let $(\pi,V)$ be the regular representation of $(A,P,\alpha)$. Thanks to Corollary \ref{functorial_ideal}, we can identify $I \rtimes_{red} P$ with the $C^{*}$-algebra generated by $\{V_s\pi(x):s \in P, x\in I\}$. Let $(\widetilde{\pi},\widetilde{V})$ be the regular representation of $(B,P,\beta)$. By Prop. \ref{equivariance implies homomorphism}, there exists a $*$-homomorphism $\widetilde{\phi}: A \rtimes_{red}P \to B \rtimes_{red} P$ such that 
\[
\widetilde{\phi}(V_s\pi(x))=\widetilde{V}_s\widetilde{\pi}(\phi(x))
\] for $s \in P$ and $x \in A$. By Corollary \ref{functorial_one}, there exists a $*$-homomorphism $\widetilde{\epsilon}:B \rtimes_{red} P \to A \rtimes_{red} P$ such that 
\[
\widetilde{\epsilon}(\widetilde{V}_s\widetilde{\pi}(b))=V_s\pi(\epsilon(b))
\]
for $b \in B$ and $s \in P$. It is clear that $\widetilde{\phi} \circ \widetilde{\epsilon}$ is the identity map. Also, $I \rtimes_{red} P \subset Ker(\widetilde{\phi})$.

\textit{Claim:} For every $d \in A \rtimes_{red} P$, $d-\widetilde{\epsilon} \circ\widetilde{\phi}(d) \in I \rtimes_{red} P$.

Let \[\cla_0:=\{d \in A \rtimes_{red} P: d-\widetilde{\epsilon}(\widetilde{\phi}(d)) \in I \rtimes_{red} P\}.\]
Since $I \rtimes_{red} P$ is an ideal, it follows that $\cla_0$ is a $*$-subalgebra of $A \rtimes_{red} P$. Also, $\cla_0$ is norm closed. Let $s \in P$ and $x \in A$ be given. Since $x-\epsilon \circ \phi(x) \in I$, it follows that 
\[
V_s\pi(x)-\widetilde{\epsilon}(\widetilde{\phi}(V_{s}\pi(x)))=V_s\pi(x-\epsilon \circ \phi(x)) \in I \rtimes_{red} P.
\]
Hence, $V_s\pi(x) \in \cla_0$. Since $\{V_s\pi(x):s \in P, x \in P\}$ generates $A \rtimes_{red} P$, it follows that $\cla_0=A \rtimes_{red} P$. The proof of the claim is over. 

Let $d \in Ker(\widetilde{\phi})$. Then, $d=d-\widetilde{\epsilon}(\widetilde{\phi}(d)) \in I \rtimes_{red} P$. Hence, $Ker(\widetilde{\phi}) \subset I \rtimes_{red} P$. The other inclusion is obvious. 
\hfill $\Box$

Let $(A,P,\alpha)$ be a semigroup dynamical system. Let $A^{+}:=A \oplus \bbc$ be the unitisation. For $s \in P$, let $\alpha_{s}^{+}:A^{+} \to A^{+}$ be defined by \[\alpha_s^{+}((a,\lambda)):=(\alpha_s(a),\lambda).\] Then, $(A^{+},P,\alpha^{+})$ is a unital semigroup dynamical system called the unitization of $(A,P,\alpha)$.

\begin{crlre}[{\cite[Lemma 3.5]{Sundar_Khoshkam}}]
\label{useful_split_exact}
    With the above notation, we have a short exact sequence 
    \[
    0 \longrightarrow A \rtimes_{red} P \longrightarrow A^{+} \rtimes_{red} P \longrightarrow C_{red}^{*}(P)\longrightarrow 0\]
        which is also split-exact.
    \end{crlre}

We review here the groupoid crossed product presentation of the Wiener-Hopf algebra $\mathcal{W}(A,P,G,\alpha)$ (Defn.~\ref{first_Wiener}) of a semigroup dynamical system $(A,P,\alpha)$. 
 Let us first recall the definition of the \emph{Wiener-Hopf groupoid} which first appeared in the work of Muhly and Renault (\cite{Renault_Muhly}) in their study of Wiener-Hopf algebras associated with cones in an Euclidean space. 

Let $P$ be a subsemigroup of a discrete, countable group $G$. Denote the power set of $G$ by $\mathcal{P}(G)$. Then, $\mathcal{P}(G)$ is a compact metric space where the convergence is as follows: 
for a sequence $(F_n) \in \mathcal{P}(G)$ and $F \in \mathcal{P}(G)$, $F_n \to F$ iff $1_{F_n}(g) \to 1_{F}(g)$ for every $g \in G$. 
Note that $G$ acts on the right on $\mathcal{P}(G)$, where the action is given by the map 
\[
\mathcal{P}(G) \times G \ni (F,g) \mapsto Fg \in \mathcal{P}(G).
\]
Define 
\begin{align*}
    \Omega:&=\overline{\{P^{-1}a: a \in P\}}\\
    \widetilde{\Omega}:&=\bigcup_{g \in G}\Omega g.
\end{align*}
Then, $\Omega$ is compact, and $\widetilde{\Omega}$ is locally compact. For $F \in \widetilde{\Omega}$, $F \in \Omega$ if and only if $e \in F$.
Moreover, $G$ leaves $\widetilde{\Omega}$ invariant. The \emph{Wiener-Hopf groupoid} $\mathcal{G}$ associated with the semigroup $P$ is defined to be the reduction of the transformation groupoid $\widetilde{\Omega}\rtimes G$ to the clopen set $\Omega$, i.e. as a set
\[
\mathcal{G}:=\widetilde{\Omega}\rtimes G|_{\Omega}=\{(F,g): F \in \Omega, Fg \in \Omega\},
\]
 and the groupoid multiplication and inversion are as follows: \[
(F,g)(G,h)=(F,gh) \textrm{~~if $Fg=G$}; \quad (F,g)^{-1}=(Fg,g^{-1}).
\]
If we wish to stress the dependence of $\Omega$, $\widetilde{\Omega}$, and $\mathcal{G}$ on $P$, we denote them by $\Omega_P$, $\widetilde{\Omega}_P$ and $\mathcal{G}_P$, respectively. 

\begin{rmrk}
\label{Equivalence of groupoids amenability}
Note that since $\widetilde{\Omega}=\bigcup_{g \in G}\Omega g$ and $\Omega$ is clopen in $\widetilde{\Omega}$, $\mathcal{G}$ and $\widetilde{\Omega} \rtimes G$ are equivalent as groupoids (see \cite[Example 2.7]{MRW}). Thus, $\mathcal{G}$ is amenable if and only if $\widetilde{\Omega}\rtimes G$ is amenable (see \cite[Thm. 2.2.17]{Anantharaman}). 
\end{rmrk}

Let $(A,P,\alpha)$ be a unital semigroup dynamical system. We assume that $A$ is separable. We recall from \cite[Section 4]{Sundar_Khoshkam} the construction of the groupoid dynamical system associated with $(A,P,\alpha)$. Let $\ell^{\infty}(G,A)$ denote the $C^{*}$-algebra of bounded $A$-valued functions.  For $s \in G$, let $\beta_s$ be the automorphism of $\ell^{\infty}(G,A)$ given by  $\beta_s(f)(h):=f(hs)$. Then, $\beta:=\{\beta_s\}_{s \in G}$ is a group of automorphisms of $\ell^{\infty}(G,A)$.  For $g \in G$ and for $x \in A$, let $j_g \in \ell^{\infty}(G,A)$ be defined by 
\begin{equation}
\label{jg}
j_g(x)(h):=\begin{cases}
    \alpha_{hg^{-1}}(x)& \mbox{ if
} hg^{-1} \in P,\cr
   &\cr
    0 &  \mbox{ if } hg^{-1} \notin P.
         \end{cases}
\end{equation}
Let $\widetilde{D}$ be the $C^{*}$-algebra generated by $\{j_g(x):g \in G, x \in A\}$. Note that $\widetilde{D}$ is $G$-invariant.  For, $\beta_s(j_g(x))=j_{gs^{-1}}(x)$ for $s,g \in G$ and $x \in A$.
Also, there exists a unique injective $G$-equivariant $*$-homomorphism $\phi:C_0(\widetilde{\Omega}) \to \widetilde{D}$ such that 
\begin{equation}\label{eq-embed-inside-D}
  \phi(1_{\Omega g})=j_g(1).  
\end{equation}
Note that for $A \in \widetilde{\Omega}$ and $g \in G$, $1_{\Omega g}(A)=1_{A}(g)$. Hence, $\{1_{\Omega g}: g \in G\}$ generates $C_{0}(\widetilde{\Omega})$. Then, it follows from Eq.~\ref{eq-embed-inside-D} that 
\begin{equation}
    \label{eq-embed-inside-D-one}
    \phi(f)(g)=f(P^{-1}g)
\end{equation}
for $f \in C_{0}(\widetilde{\Omega})$ and $g \in G$. 

 This way, we identify $C_0(\widetilde{\Omega})$ as a $*$-subalgebra of $\widetilde{D}$ (see \cite[Section 4]{Sundar_Khoshkam} and \cite[Remark 4.2]{Sundar_Khoshkam}). Moreover, $C_0(\widetilde{\Omega}) \subset Z(\widetilde{D})$. Hence, $\widetilde{D}$ is a $C_{0}(\widetilde{\Omega})$-algebra which can now be realised as the section algebra of an upper semicontinuous bundle 
$\widetilde{\mathcal{D}}:= \displaystyle \coprod_{F \in \widetilde{\Omega}}\mathcal{D}_F$. 
For $F \in \widetilde{\Omega}$, fibre $\widetilde{\mathcal{D}}_F$ is defined by $\widetilde{\mathcal{D}}_F:=\widetilde{D}/{I_F}$, where $I_F:=C_{0}(\widetilde{\Omega}\backslash\{F\})\widetilde{D}$. 
Thanks to the equivariance of the homomorphism $\phi:C_{0}(\widetilde{\Omega}) \to \widetilde{D}$, we have an action $\beta:=\{\beta_{(F,g)}\}$ of the  transformation groupoid $\widetilde{\Omega} \rtimes G$  on the bundle $\widetilde{\mathcal{D}}$ given by the following formula: \begin{equation}
    \label{groupoid action}
    \beta_{(F,g)}(d+I_{Fg})=\beta_g(d)+I_{F}
\end{equation}
for $(F,g) \in \widetilde{\Omega} \rtimes G$ and $d \in \widetilde{\mathcal{D}}$. 
Denote the bundle $\widetilde{\mathcal{D}}$ restricted to the clopen set $\Omega$ by $\mathcal{D}$. Since the Wiener-Hopf groupoid $\mathcal{G}$ is the restriction $\widetilde{\Omega} \rtimes G|_{\Omega}$, the groupoid $\mathcal{G}$ acts on $\mathcal{D}$ whose action we again denote by $\beta$. Moreover, for $(F,g) \in \mathcal{G}$, $\beta_{(F,g)}$ is given by the same formula as in Eq.~\ref{groupoid action}.

The following theorem is the main result of \cite{Sundar_Khoshkam}. 
\begin{thm}[{\cite[Thm. 4.3]{Sundar_Khoshkam}}]
\label{unital_case_KS}
Keep the above notation. Suppose that $(P,G)$ satisfies the Toeplitz condition. Then, $A \rtimes_{red} P$ and  the reduced crossed product $\mathcal{D} \rtimes_{red} \mathcal{G}$ of the groupoid dynamical system $(\mathcal{D},\mathcal{G})$ are isomorphic.  Moreover, $A \rtimes_{red} P$ and $\widetilde{\mathcal{D}}\rtimes_{red} (\widetilde{\Omega} \rtimes G)\cong \widetilde{D} \rtimes_{red} G$ are Morita equivalent.
\end{thm}
\textit{Proof.} By Remark \ref{Toeplitz condition}, it follows that $A \rtimes_{red} P$ coincides with the Wiener-Hopf algebra $\mathcal{W}(A,P,G,\alpha)$. Now, the first statement is exactly~\cite[Thm. 4.3]{Sundar_Khoshkam}. 
As far as the second statement is concerned, note that the groupoids $\mathcal{G}$ and $\widetilde{\Omega} \rtimes G$ are equivalent as $\Omega$ is clopen and $\widetilde{\Omega}=\bigcup_{g \in G}\Omega g$. By construction, $(\mathcal{D},\mathcal{G})$ and $(\widetilde{\mathcal{D}},\widetilde{\Omega} \rtimes G)$ are equivalent in the sense of \cite[Defn. 5.1]{MW08}.  The second conclusion follows. \hfill $\Box$

\begin{rmrk}
The non-unital version of Thm.~\ref{unital_case_KS} with a slight modification in the definition of the bundles $\mathcal{D}$ and $\mathcal{\widetilde{D}}$ also holds. The non-unital version together with Thm.~\ref{main_existence} and Thm.~\ref{main_Muhly} applied to the reduced $C^{*}$-algebra of a proper product system $X$ says that, under `good conditions', $C_{red}^{*}(X)$ is Mortia equivalent to a groupoid crossed product.
However, to derive structural results concerning  $C_{red}^{*}(X)$, thanks to Corollary~\ref{useful_split_exact}, we can pass to the unitization, and it then suffices to prove the desired results only for unital semigroup dynamical systems. For this reason, we have not included the details of the non-unital version of Thm.~\ref{unital_case_KS}.   For more details, we refer the reader to~\cite{Sundar_Khoshkam}.
Note that Thm.~\ref{unital_case_KS} for the trivial dynamical system $(\bbc,P,\alpha)$, where $\alpha$ is the trivial action, says that $C_{red}^{*}(P) \cong C^{*}_{red}(\clg)$ and the Morita equvalence $C_{red}^{*}(\clg)\sim_{M} C_{0}(\widetilde{\Omega})\rtimes_{red} G$ forms the first step in the computation of the $K$-theory of semigroup $C^{*}$-algebras (\cite{Li_Oberwolfach}).
\end{rmrk}

\section{Nuclearity and exactness}
\label{sec-nuc-exact}

In this section, we discuss the exactness and nuclearity of the reduced crossed product $A \rtimes_{red} P$ of a semigroup dynamical system $(A,P,\alpha)$, and then apply the results to the reduced $C^{*}$-algebra of a proper product system. We will keep the notation of the previous section for the remainder of this paper. For the rest of this paper, $P$ stands for a subsemigroup of a group~$G$. 

We need a better understanding of the fibre $\mathcal{D}_F$ to proceed further. 
First, we introduce some notation. For $x,y \in G$, we say $x \leq y$ if $yx^{-1} \in P$. A subset $F \subset G$ is said to be directed if given $x,y \in F$, there exists $z \in F$ such that $x,y \leq z$.
A subset $F \subset G$ is hereditary if  $x \in F$ and $y \leq x$, then $y \in F$.  For $F \subset G$, if $\displaystyle F=\lim_{n \to \infty}P^{-1}g_n$ for some sequence $(g_n)_n$ in $G$, then $P^{-1}F \subset F$. Hence, $F$ is  hereditary for every $F \in \widetilde{\Omega}$. 

\begin{lmma}
\label{crucial_repeat}
Let $s \in G$, $a \in A$ and $F \in \widetilde{\Omega}$. Then, $j_s(a)+I_F=0$ if $s \notin F$. The $C^{*}$-algebra $\widetilde{\mathcal{D}}_F$ is generated by $\{j_s(a)+I_F: s \in F, a \in A\}$
\end{lmma}
\textit{Proof.} Suppose that $s \notin F$. Note that $j_s(a)=j_s(1)j_s(a)$. By Eq. \ref{eq-embed-inside-D}, \[j_s(1)(F)=1_{\Omega s}(F)=1_{\Omega}(Fs^{-1})=1_{F}(s)=0.\] Hence, $j_s(a)+I_F=0$. The second assertion is now clear as $\widetilde{D}$, by definition, is generated by $\{j_s(a):a \in A,s \in G\}$.\hfill $\Box$

\begin{ppsn}
\label{inductive limit}
Let $(A,P,\alpha)$ be a unital semigroup dynamical system.   Suppose that every element of $\widetilde{\Omega}$ is directed. For $F \in \widetilde{\Omega}$, consider the following directed system of $C^{*}$-algebras: for $s \in F$, let $B_s:=A$, and for $s \leq t$, let $\phi_{t,s}:B_s \to B_t$ be the connecting map defined by $\phi_{t,s}=\alpha_{ts^{-1}}$. Let \ $\displaystyle B_F:= \underset{s\in F}{\varinjlim} B_s$ be the inductive limit. Then, $\widetilde{\mathcal{D}}_F$ and $B_F$ are isomorphic.
\end{ppsn}
\textit{Proof.}	Clearly, \(\big(B_s, \{\phi_{t,s}\}_{s\leq t\in F}\big)\) is an inductive system of \(C^*\)-algebras. Suppose  \((B_F,\{\lambda_s\}_{s \in F})\) is the inductive limit, where $\lambda_s:B_s \to B_F$ is the canonical map. 
	For $s \in F$, we define \(u_s\colon B_s \to \widetilde{\mathcal{D}}_F\) by \[u_s(a): = j_s(a) + I_F.\] For \(s,t\in F\) and \(s\leq t\), we claim that  the following diagram 
	\begin{center}
	\begin{tikzcd}
		 B_s \arrow[rd, "u_s"'] \arrow[r, "\phi_{t,s}"] & B_t \arrow[d, "u_t"] \\
		& \widetilde{\mathcal{D}}_F
	\end{tikzcd}
		\end{center}
		commutes. To see that, let \(a\in B_s\) and $s,t \in F$ be such that $s \leq t$. We have 
		\[
		u_t\circ \phi_{t,s}(a)= j_t(\alpha_{ts^{-1}}(a)) + I_F.
		\]
		Now for \(h\in G\), we have
		\begin{align*}
			j_t(\alpha_{ts^{-1}}(a))(h) &=  \begin{cases}
			 \alpha_{ht^{-1}}(\alpha_{ts^{-1}}(a)) & \textup{if } h\in Pt,\\
				0 & \textup{otherwise};
			\end{cases}\\
			& = \begin{cases}
				\alpha_{hs^{-1}}(a) & \textup{if } h\in Pt,\\
				0 & \textup{otherwise}.
			\end{cases}
		\end{align*}
	Using the above computation, we  get the following equality
	\begin{equation}\label{eq-ind-lim-fibre-1}
	j_s(a) - j_t(\alpha_{ts^{-1}}(a)) = \bigl(j_s(1) - j_t(1) \bigr)j_s(a).
	\end{equation}
	By Eq.~\ref{eq-embed-inside-D}, we have 
	\[
	 j_s(1)(F)-j_t(1)(F) = 1_{\Omega s}(F) - 1_{\Omega t}(F) = 1_F(s) -1_F(t) =0.
	\]
	Therefore, Eq.~\ref{eq-ind-lim-fibre-1} ensures that \(u_t\circ \phi_{s,t} = u_s\) for \(s,t\in F\).
	Thus, there exists a \(^*\)-homomorphism \(\mu \colon B_F\to \widetilde{\mathcal{D}}_F\) such that \(\mu \circ \lambda_s=u_s\) for \(s\in F\).
	
	For \(d\in \widetilde{D}\), we claim that the net \((\lambda_g(d(g)))_{g\in F}\) is a Cauchy net in \(B_F\). It suffices to prove when $d=j_s(a)$ for some $s \in P$ and $a \in A$. Thus, let $s \in G$, let $a \in A$ and let $d=j_s(a)$. Since $F$ is hereditary, if $s \notin F$, then $s \nleq g \implies j_s(a)(g)=0$ for every $g \in F$. Thus, in this case, $(\lambda_g(d(g)))_{g \in F}$ is convergent in $B_F$.
    Suppose $s \in F$. Then,		\begin{align*}
		\lambda_h\big(j_s(a)(h)\big) &=  \begin{cases}
			\lambda_{h}(\alpha_{hs^{-1}}(a)) & \textup{if } h\in Ps,\\
			0 & \textup{otherwise};
		\end{cases}\\
		& = \begin{cases}
			\lambda_s(a) & \textup{if } h\in Ps,\\
			0 & \textup{otherwise}.
		\end{cases}
	\end{align*}
    Hence, \begin{equation}
    \label{limitexists}
   \displaystyle \lim_{g \in F}\lambda_g(j_s(a)(g))=\displaystyle \lim_{h \in F, h \geq s}\lambda_h(j_s(a)(h))=\lambda_s(a).\end{equation} This proves the claim. Define $\lambda:\widetilde{D} \to B_F$ by setting
   \[
   \lambda(d):=\lim_{g \in F}\lambda_g(d(g)).
   \]
\textit{Claim:} The map $\lambda$ vanishes on $I_F$. 

Let $f \in C_{0}(\widetilde{\Omega}\backslash \{F\})$ and $d \in \widetilde{D}$ be given. Let $d_1=f d$. Since $F$ is directed, the net $(P^{-1}g)_{g \in F} \to F$ in $\widetilde{\Omega}$. Therefore, $\lim_{g \in F}f(P^{-1}g)= f(F)=0$.  Thanks to Eq. \ref{eq-embed-inside-D-one},
\[
\lambda(d_1)=\lim_{g \in F}\lambda_g(d_1(g))=\lim_{g \in F}f(P^{-1}g)\lambda_g(d(g)) = f(F)\lim_{g \in F}\lambda_g(d(g))=0.
\]
This proves the claim. 

Hence, the map $\lambda$ descends to a homomorphism $\lambda_1:\widetilde{\mathcal{D}}_F \to B_F$. 
We next check that $\lambda_1$ and  $\mu$ are inverses of each other. 
        Let $s \in F$ and $a \in A$ be given.  We have 
    \begin{align*}
    \mu \circ \lambda_1 (j_s(a)+I_F)& = \mu \big( \lim_{g\in F} \lambda_g(j_s(a)(g))\big)\\
    &= \mu (\lambda_s(a)) \quad \quad (\textrm{by Eq.~\ref{limitexists}})\\
       &= j_s(a)+I_F.
    \end{align*}
    Since  \(\{j_s(a)+I_F : s\in F, a\in A\}\) generates \(\mathcal{\widetilde{D}}_F\) (by Lemma \ref{crucial_repeat}),  we conclude that \(\mu \circ \lambda_1 = Id_{\mathcal{\widetilde{D}}_F}.\)
	On the other hand, for $s \in F$ and $a \in A$, thanks to Eq.~\ref{limitexists}, 
	\[
	 \lambda_1 \circ \mu (\lambda_s(a)) = \lambda(u_s(a)) = \lambda_1(j_s(a)+ I_F) = \lambda(j_s(a))= \lambda_s(a).
 	\]
 Hence, \(\lambda_1 \circ \mu = Id_{B_F}\). 
	Therefore, \(\mathcal{\widetilde{D}}_F\) and \(B_F\) are isomorphic. The proof is over. \hfill $\Box$ 

To prove the exactness part of Thm.~\ref{main_nuclearity}, we need the following proposition. If $P$ is left Ore, i.e. $P^{-1}P=G$, and when the action is injective, then the fibre $\widetilde{\mathcal{D}}_F$ can be described as in the following proposition. This case was considered in \cite{SundarJFA} where a groupoid dynamical system was constructed.  We verify in the next proposition that the bundle constructed in \cite{SundarJFA} coincides with the one described in this paper.

\begin{ppsn}
\label{picture for car}
    Suppose $P$ is left Ore in $G$, i.e. $P^{-1}P=G$. Let $(A,P,\alpha)$ be a unital semigroup dynamical system. Suppose that $\alpha_s$ is injective for every $s \in P$.  Let $(B,G,\beta,\phi)$ be  the Laca dilation of $(A,P,\alpha)$, i.e. $(B,G,\beta)$ is a $C^{*}$-dynamical system, $\phi:A \to B$ is an injective map that is $P$-equivariant such that $\displaystyle \overline{\bigcup_{s \in P}\beta_{s^{-1}}(\phi(A))}=B$. 
Let \[
\Gamma:=\{f:C_{0}(\widetilde{\Omega}, B): \textrm{$f(F) \in \cla_F$ for every $F \in \widetilde{\Omega}$} \}.
\]
Then, $\widetilde{\mathcal{D}}$ is isomorphic to $ \Gamma$. Moreover, the fibre  $\widetilde{\mathcal{D}}_F$ is isomorphic to the $C^{*}$-algebra generated by $\{\beta_{g}^{-1}(\phi(a)):a \in A, g \in F\}$.
\end{ppsn}
\textit{Proof.} We abuse notation and consider $A$ as a subalgebra of $B$, and suppress $\phi$ from the notation. Moreover, we consider $\ell^{\infty}(G,A)$ as a subalgebra of $\ell^{\infty}(G,B)$. We do not distinguish between $\beta$ and $\alpha$, and we use the same letter $\alpha$ to denote both of them.

For $F \in \widetilde{\Omega}$, let $\cla_F:=C^{*}\{\alpha_g^{-1}(x):x \in A, g \in F\}$. Then, $\{\cla_F\}_{F \in \widetilde{\Omega}}$ is an upper semi-continuous bundle with the space of continuous sections $\Gamma$ given by 
\[
\Gamma:=\{f:C_{0}(\widetilde{\Omega}, B): \textrm{$f(F) \in \cla_F$ for every $F \in \widetilde{\Omega}$} \}.
\]
We claim that the map $R:\Gamma \to \ell^{\infty}(G,B)$ defined by 
\[
R(f)(g):=\alpha_{g}(f(P^{-1}g)).\]
is a $C_0(\widetilde{\Omega})$ isomorphism between $\Gamma$ and $\widetilde{D}$. 
Note that $R$ is one-one as $\{P^{-1}g:g \in G\}$ is dense in $\widetilde{\Omega}$. 
Observe that for $s \in G$ and $x \in A$, \[
R(1_{\Omega s}\otimes \alpha_{s}^{-1}(x))=j_s(x).
\]
Now the claim follows by applying (\cite{Dana-Williams}, Prop. C.25). The proof is over. \hfill $\Box$

We now prove Thm.~\ref{main_nuclearity}.
\medskip

    \textit{Proof of Thm.~\ref{main_nuclearity}.} 
    \noindent (1). Suppose that $P^{-1}P=G$ and $\alpha_s$ is injective for every $s \in P$. Let $(A^{+},P,\alpha^{+})$ be the unitisation of $(A,P,\alpha)$. It is clear that, for $s \in P$, $\alpha_{s}^{+}$ is injective. Since a $C^{*}$-subalgebra of an exact $C^{*}$-algebra is exact, thanks to Corollary \ref{useful_split_exact}, it suffices to prove the exactness result when  $(A,P,\alpha)$ is unital.  Suppose $A$ is exact. For $g \in G$, let $A_g:=\beta_{g}^{-1}(\phi(A))$. Then, $A_g$ is exact.  Note that $A_g \subset A_h$ if $g \leq h$. Since $B$ is the closure of an increasing net $(A_g)_{g \in G}$ of $C^{*}$-algebras isomorphic to $A$, it follows that $B$ is exact. It follows from Thm.~\ref{picture for car} that $\widetilde{D} \subset C_{0}(\widetilde{\Omega},B)$ which implies that $\widetilde{D}$ is exact. Thanks to~\cite[Thm. 6.6]{Lalonde}, $\widetilde{D} \rtimes_{red} (\widetilde{\Omega} \rtimes G)$ is exact. As $A \rtimes_{red} P$ and $\widetilde{D}\rtimes_{red}(\widetilde{\Omega} \rtimes G)$ are Morita equivalent, it follows that $A \rtimes_{red} P$ is exact.
      Conversely, assume that \(A\rtimes_{red}P\) is exact. Since $A$ is a $C^{*}$-subalgebra of $A \rtimes_{red} P$, $A$ is exact.

    Next, we prove the nuclearity result. Assume that for every $F \in \widetilde{\Omega}$, $F$ is directed. 
    As ideals of nuclear $C^{*}$-algebras are nuclear, thanks to Corollary ~\ref{useful_split_exact}, it suffices to consider the case when  $(A,P,\alpha)$ is unital.  Let the notation be as in Prop.~\ref{inductive limit}. 
	Since \((P,G)\) satisfies the Toeplitz condition, we have that \(A \rtimes_{red} P\) is isomorphic to the groupoid crossed product \(\mathcal{D} \rtimes_{red} \mathcal{G}\) (see Thm.~\ref{unital_case_KS}). As we have assumed that $\mathcal{G}$ is amenable, the reduced crossed product $\mathcal{D}\rtimes_{red} \mathcal{G}$ and the full crossed product $\mathcal{D} \rtimes \mathcal{G}$ are isomorphic.  
	 By ~\cite[Thm. 4.1]{Takeishi-2014-Nuclearity-of-Fell-bund-alg-over-etale-gpd}, \(\mathcal{D} \rtimes_{red} \mathcal{G}\) is nuclear if and only if \(\mathcal{D}_F\) is nuclear for every \(F \in \Omega\). Suppose that $A$ is nuclear. Since, for $F \in \Omega$, $\mathcal{D}_F=\widetilde{\mathcal{D}}_F$ is the inductive limit of \(\{B_s=A, (\phi_{t,s})_{s \leq t \in F}\}\), it follows that \(\mathcal{D}_F\) is nuclear. 
     
     Conversely, suppose $A \rtimes_{red} P=\mathcal{D}\rtimes_{red} \mathcal{G}$ is nuclear. By  \cite[Thm. 4.1]{Takeishi-2014-Nuclearity-of-Fell-bund-alg-over-etale-gpd}, it follows that $\mathcal{D}_F$ is nuclear for every $F$. Note that for $F=P^{-1}$, $\mathcal{D}_F=\widetilde{\mathcal{D}}_F \cong B_F=A$ as the identity element $e$ is an upper bound for $F$. Hence, $A$ is nuclear.

    \noindent (2). As nuclearity and exactness are preserved under Morita equivalence, the product system version follows from Thm.~\ref{main_Muhly}, Thm.~\ref{main_existence}, Remark \ref{lem:comp-faithful} and the fact that for a full Hilbert $B$-module $\cle$, $B$ and $\clk_B(\cle)$ are Morita equivalent. 
\hfill $\Box$

We next show that Thm.~\ref{main_nuclearity} is applicable for finitely generated subsemigroups of abelian groups. Note that if $P$ is a finitely generated subsemigroup of an abelian group $G$ and if $\{a_1,a_2,\cdots,a_n\}$ generates $P$ as a semigroup, then $\sum_{i=1}^{n}a_i$ is an order unit. Thus, we only need to check the directedness hypothesis. We need the following lemma for such a verification. 
\begin{lmma}
\label{finitely generated}
Let \(P\) be a subsemigroup of a countable group \(G\).
Suppose that given $g,h \in G$, either $Pg \cap Ph=\emptyset$ or there exists $n \in \bbn$ and $g_1,g_2,\cdots,g_n \in G$ such that 
\[
Pg \cap Ph=\bigcup_{i=1}^{n}Pg_i.
\]
Then, every element of $\widetilde{\Omega}$ is directed. 
\end{lmma}
\textit{Proof.}
	Let \(A\in \widetilde{\Omega}\) and \(g,h\in A\). Since right multiplication by an element of $G$ preserves the order $\leq$ and $\Omega=\displaystyle \bigcup_{s \in G}\Omega s$, we can assume that $A \in \Omega$.  Then, there is a sequence \((a_n)_n\) in $P$ such that $P^{-1}a_n \to A$. Set \(A_n = P^{-1}a_n\). Then, each \(A_n\) is directed as \(a_n\) is an upper bound of any two points in \(A_n\). 
	Since \(A_n \to A\), we have \(1_{A_n}(g) \to 1_{A}(g) =1\) and \(1_{A_n}(h) \to 1_{A}(h) =1\). Thus, \(g,h\in A_n\) eventually. Since \(A_n\) is directed, for large $n$,  there exist  $x_n \in A_n$ such that \(x_n\geq g,h\). This gives us \(x_n \in Pg\cap Ph\) for all \(n\in \mathbb{N}\).  Hence, $Pg \cap Ph\neq \emptyset$, and let $g_1,g_2,\cdots,g_r$ be such that $\displaystyle Pg \cap Ph=\bigcup_{k=1}^{r}Pg_k$. Since $ \displaystyle x_n \in Pg \cap Ph=\bigcup_{k=1}^{r}Pg_k$,  there exist \(g_l\) and a subsequence \((x_{n_k})\) such that \(x_{n_k} \in Pg_{l}\) for all \(k\).
    The hereditariness of \(A_{n}\) ensures that 
	\(g_{l}\in A_{n_k}\) for all \(k\). Therefore, 
	\[
	1=1_{A_{n_k}}(g_{l}) \to 1_{A}(g_{l}).
	\]
	which implies \(g_l\in A\). Since the identity \(e\in P\), \(g_l\in \bigcup_{i=1}^{r}Pg_i =	Pg \cap Ph\). Thus, \(g_l\geq g,h\), and hence, \(A\) is directed. 
\hfill $\Box$

\begin{rmrk}
If the hypothesis of the above lemma is satisfied with $n=1$, the pair $(P,G)$ is said to be quasi-lattice ordered (\cite{Nica92}).
\end{rmrk}

 The first statement in the following proposition is known. A proof is given for completeness. 
\begin{ppsn}
\label{finitely_generated_Noetherian}
    Let $G$ be an abelian group, and let $P \subset G$ be a finitely generated subsemigroup of $G$ such that $P-P=G$. Then, we have the following. 
    \begin{enumerate}
        \item Let $I \subset P$ be an ideal,  i.e. $I$ is non-empty and $I+P \subset P$. Then, there exists $x_1,x_2,\cdots,x_n \in I$ such that 
    \[
    I:=\bigcup_{i=1}^{n}(x_i+P).\]
\item Every element of $\widetilde{\Omega}$ is directed.
    \end{enumerate}
\end{ppsn}
\textit{Proof.} For $s \in P$, let $v_s:\ell^2(P) \to \ell^2(P)$ be the operator defined by \[v_s(\delta_t):=\delta_{s+t}.\] Let $R$ be the linear span  of $\{v_s:s \in P\}$. Then, $R$ is a unital commutative algebra that is finitely generated. Hence, $R$ is Noetherian. 

For an ideal $I \subset P$, set $
R_I:=span\{v_s: s \in I\}$. Let $I$ and $J$ be ideals of $P$ such that $I$ is a proper subset of $J$. We claim that $R_I$ is a proper subset of $R_J$.  Let $s \in J$ be such that $s \notin I$. Suppose $v_s \in R_I$. Then, 
there exists $t_1,t_2,\cdots,t_n \in I$ and complex numbers $a_1,a_2,\cdots,a_n$ such that $v_s=\sum_{i=1}^{n}a_iv_{t_i}$. Hence, 
\[
\delta_s=v_s(\delta_0)=\sum_{i=1}^{n}a_iv_{t_i}(\delta_0)=\sum_{i=1}^{n}a_i\delta_{t_i}.
\]
The above equality implies that there exists $i_0$ such that $s=t_{i_0} \in I$, which is a contradiction. Hence, $v_s \notin R_I$. This proves the claim. 

Let $I=\{x_1,x_2,\cdots\} \subset P$ be an ideal. Set $I_n:=\bigcup_{i=1}^{n}(P+x_i)$. Note that $I_n$ is an ideal and $I_n \subset I_{n+1}$ for each $n$. Hence, $(R_{I_n})_n$ is an increasing chain of ideals in $R$. Since $R$ is Noetherian, there exists $n$ such that $R_{I_k}=R_{I_n}$ for all $k \geq n$. Hence, $I_k=I_n$ for $k \geq n$. Since $\displaystyle I=\bigcup_{j=1}^{\infty}I_j$, it follows that 
$I=I_n$. This proves $(1)$. 

Let $g,h \in G$ be given. Since $P-P=G$, there exists $a,b,c,d$ such that $g=a-b$ and $h=c-d$. Replacing $a,b,c,d$ by $a+d,b+d,c+b,d+b$ respectively, we can assume that $g=a-c$ and $h=b-c$ for some $a,b,c \in P$. Note that $I:=(P+a) \cap (P+b)$ is an ideal of~$P$. Thus, there exists $x_1,x_2,\cdots,x_n \in I$ such that $I=\bigcup_{i=1}^{n}(P+x_i)$. For $i \in \{1,2,\cdots,n\}$, let $y_i=x_i-c$. Then, 
\[
(P+g) \cap (P+h)=I-c=\bigcup_{i=1}^{n}(P+y_i).
\]
We can now apply Lemma~\ref{finitely generated} to conclude that every element of $\widetilde{\Omega}$ is directed. \hfill $\Box$

\begin{rmrk}
 
     It is not true in general that, for a subsemigroup $P$ of a group $G$, every element of $\widetilde{\Omega}$ is directed. For example, consider  $P:=\bbn \times \bbn \cup \{(0,0)\}$.  Recall that as per our notation, $\bbn=\{1,2,\cdots\}$. Let $\bbz_\infty=\bbz \cup \{\infty\}$, $\bbn_\infty = \bbn \cup \{\infty\}$. We let $\bbz^2$ acts on $\bbz_\infty^{2}$ by translations.
    
     For $(a,b) \in \bbz_\infty \times \bbz_\infty$, define 
   \[
    F_{(a,b)}:=\{(x,y) \in \bbz^2: \textrm{$(x,y)=(a,b)$ or 
    ($x<a$ and $y<b$})\}.\]
     Then, the map $\bbz_\infty \times \bbz_\infty \ni (a,b) \mapsto F_{(a,b)} \in \widetilde{\Omega}$ is a $\bbz^2$-equivariant homeomorphism. Under this homeomorphism, we can identify $\Omega$ with $(\bbn_\infty \times \bbn_\infty) \cup \{(0,0)\}$.  Note that $F_{(a,b)}$ is not directed if $a<\infty$ and $b<\infty$.

\end{rmrk}

Even if we do not have the directedness hypothesis of Thm.~\ref {inductive limit} (for example, for the semigroup $P=(\bbn \times \bbn) \cup \{0\}$), the picture of the fibre $\widetilde{\mathcal{D}}_F$ given by Prop. \ref{picture for car} comes handy in certain situations. We illustrate with a semigroup dynamical system arising out of the canonical anti-commutation relation, where we can apply Prop.~\ref{picture for car} to conclude nuclearity.

\begin{xmpl}
\label{CAR_example_one}
Let $K$ be a separable Hilbert space. For a subspace $L \subset K$, let $\cla(L)$ be the CAR algebra associated with $L$, i.e. $\cla(L)$ is the universal unital $C^{*}$-algebra generated by $\{a(\xi):\xi \in L\}$ that satisfy the canonical anti-commutation relation, i.e. for $\xi,\eta \in L$, 
\begin{align*}
    a(\xi)a(\eta)+a(\eta)a(\xi)&=0,\\
    a(\xi)^{*}a(\eta)+a(\eta)^{*}a(\xi)&=\langle \xi|\eta \rangle.
\end{align*}
Since the CAR algebra is simple, if $L_1 \subset L_2$, $\cla(L_1) \subset \cla(L_2)$.

Suppose that $P^{-1}P=G$ and $PP^{-1}=G$. Let $U:=\{U_g\}_{g \in P}$ be a group of unitary operators on a separable Hilbert space $K$, and let $H \subset K$ be a subspace such that for $s \in P$, $U_sH \subset H$. For $s \in P$, set $V_{s}:=U_{s}|_{H}$. Then, $V:=\{V_s\}_{s \in P}$ is a semigroup of isometries on $H$. Thanks to the universal property, for every $s \in P$, there exists a unique unital endomorphism $\alpha_s$ of $\cla(H)$ such that 
\[
\alpha_s(a(\xi)):=a(V_s\xi).\]
Then, $\alpha:=\{\alpha_s\}_{s \in P}$ is a semigroup of endomorphisms on $\cla(H)$ which gives rise to the dynamical system $(\cla(H),P,\alpha)$. Since $\cla(H)$ is simple, $\alpha_s$ is injective for every $s \in P$. 

We can apply Prop. \ref{picture for car} to deduce that $\cla(H) \rtimes_{red} P$ is nuclear. We can, without loss of generality, assume that $K=\overline{\bigcup_{g \in G}U_{g}H}=\overline{\bigcup_{a \in P}U_a^{*}H}$. For $g \in G$, let $\beta_g$ be the automorphism of $\cla(K)$ such that
\[
\beta_g(a(\xi)):=a(U_g\xi)
\]
for every $\xi \in H$. 
The automorphism $\beta_g$ exists because of the universal property of $\cla(K)$.
Then, $(\cla(K),G,\beta,\iota)$ is the Laca dilation of $(\cla(H),P,\alpha)$, where $\iota:\cla(H) \to \cla(K)$ is the natural inclusion given by the map $\cla(H) \ni a(\xi) \mapsto a(\xi) \in \cla(K)$. 

Let $F \in \widetilde{\Omega}$. By Prop. \ref{picture for car}, \[\widetilde{\mathcal{D}}_F=C^{*}\{\beta_{g}^{-1}(a(\xi)):\xi \in H, g \in F\}=C^{*}\{a(U_g^*\xi): \xi \in H, g \in F\}=\cla(L_F),\] where $L_F$ is the closure of the linear span of $\overline{\bigcup_{g \in F}U_{g}^{*}H}$. Since $\cla(L_F)$ is nuclear, $\widetilde{\mathcal{D}}_F$ is nuclear for every $F \in \widetilde{\Omega}$. Thus, if $\widetilde{\Omega} \rtimes G$ is amenable , then  $\cla(H) \rtimes_{red}
P$ is nuclear.

\end{xmpl}

\section{Homotopy and $K$-theory}
\label{sec-homo-K-theo}

In this section, we discuss the invariance of the $K$-theory of the reduced $C^*$-algebra of a product system under homotopy.  We use the `descent principle' (see \cite{Phillips_Ecterhoff}, \cite{Li-Cuntz-Echterhoff}) in $K$-theory to deduce invariance. These $K$-theoretic techniques, in the context of semigroup $C^{*}$-algebras, appeared in the work of Cuntz, Echterhoff and Li (\cite{Li-Cuntz-Echterhoff}), where they calculated the $K$-theory of semigroup $C^{*}$-algebras (i.e. the reduced $C^{*}$-algebra of the trivial product system) for a class of semigroups that satisfy a technical condition called \emph{independence}. For more on the computation of $K$-groups of semigroup $C^{*}$-algebras,  we refer the reader to~\cite{Li_Oberwolfach}.

 We begin by considering the $C(Z)$-version of $E_0$-semigroups, dynamical systems and product systems.  In this section, we only consider the case where the coefficient algebra is unital.  
Let $\cle$ be a full Hilbert $B$-module. We assume that $B$ is a unital \(C^*\)-algebra, whose center is denoted $Z(B)$. For $b \in Z(B)$, the map $R_b:\cle \to \cle$ defined by $R_b(x):=xb$ is adjointable, and for $b \in Z(B)$, $R_b \in Z(\cll_B(\cle))=Z(M(\clk_B(\cle)))$. 
Let $Z$ be a compact metric space. Suppose $B$ is a $C(Z)$-algebra with the $C(Z)$-structure given by the homomorphism $\rho:C(Z) \to Z(B)$. Observe that $\clk_B(\cle)$ is a $C(Z)$-algebra, where the $C(Z)$-structure is given by the homomorphism $C(Z) \ni f \mapsto R_{\rho(f)} \in Z(\cll_{B}(\cle))=Z\bigl(M(\clk_B(\cle))\bigr)$. 
We suppress notation and we write $fb$ or $f\cdot b$ instead of  $\rho(f)(b)$. Similarly, we write $ef$ or $e\cdot f$ in place of $R_{\rho(f)}(e)$. 

\begin{dfn}
    Let $Z$ be a compact metric space, and let  $B$ be a $C(Z)$-algebra.
    \begin{enumerate}
  \item  Let $\alpha:=\{\alpha_s\}_{s \in P}$ be an $E_0$-semigroup  on $\cll_B(\cle)$, where $\cle$ is a full Hilbert $B$-module. We say that $\alpha$ is a $C(Z)$-$E_0$-semigroup if for $s \in P$, $T \in \cll_B(\cle)$ and $f \in C(Z)$, $\alpha_s(fT)=f\alpha_s(T)$, i.e. $\alpha_s(R_{\rho(f)}T)=R_{\rho(f)}\alpha_s(T)$.

    \item Let $X:=\{X_s\}_{s \in P}$ be a product system over $B$. We say $X$ is a $C(Z)$-product system if for every $s \in P$, $f \in C(Z)$ and $x \in X_s$, $f\cdot x=x\cdot f$. \footnote{Note that $f\cdot x$ makes sense as $X_s$ carries a left action of $B$. More precisely, $f\cdot x=\rho(f)\cdot x$ where $\rho:C(Z) \to Z(B)$ is the homomorphism that determines the $C(Z)$-structure on $B$. Similarly, $x\cdot f=x\cdot \rho(f)$.}

    \item Let $(A,P,\alpha)$ be a semigroup dynamical system, where $A$ is a $C(Z)$-algebra. We call $(A,P,\alpha)$ a $C(Z)$-semigroup dynamical system if $\alpha_s(f\cdot a)=f\cdot \alpha_s
    (a)$ for $f \in C(Z)$, $a \in A$ and $s \in P$. 
\end{enumerate}
    \end{dfn}

We make a list of assertions whose proofs we omit as they are not difficult. Let $Z$ be a compact metric space. Let $B$ be a $C(Z)$-algebra. We assume that $B$ is separable.  For $z \in Z$, the fibre over $z$ is denoted $B^z$, and we let $ev_{z}: B \to B^z$ be the evaluation map. Let $z \in Z$ be given. 

\begin{enumerate}
\item[(1)] Suppose that $X:=\{X_s\}_{s \in P}$ is a product system of $B$-$B$-correspondences over $P$. Assume that $X$ is also a $C(Z)$-product system. For $s \in P$, set 
\[
X_{s}^{z}:=X_{s}\otimes_B B^z=X_{s}\otimes_{ev_z}B^z.
\]
For $s \in P$, $X_{s}^{z}$ is a $B^z$-$B^z$-correspondence, where the left action of $B_z$ is given by the formula
\[
\widetilde{b}\cdot (u \otimes \widetilde{c})=bu \otimes \widetilde{c}\]
for $\widetilde{b} \in B^z$ and $u \otimes \widetilde{c} \in \cle \otimes_{ev_z} B^z$. 
The condition that $f\cdot u=u\cdot f$ for every $f \in C(Z)$ and $u \in X_s$ ensures that the above left action is well defined. 
Then, $X^z:=\{X^z_{s}\}_{s \in P}$ is a product system of $B^z$-$B^z$ correspondences with the product rule given by 
\[
(u \otimes b)\cdot (v \otimes c):=u(bv)\otimes c
\]
for $u \otimes b \in X_s^{z}$ and $v \otimes c \in X_{t}^z$. If $X$ is proper, then $X^z$ is proper.
The reader is also referred to~\cite{Gillaspy} for more details. 

\item[(2)] Let $X$ be a $C(Z)$-product system over $P^{op}$, and suppose $(\cle,\sigma=\{\sigma_s\}_{s \in P})$ is a left dilation of $X$. Let $\alpha$ be the $E_0$-semigroup associated with $(\cle,\sigma)$, i.e. for $s \in P$, the endomorphism $\alpha_s:\cll_B(\cle) \to \cll_B(\cle)$ is defined by $\alpha_s(T):=\sigma_s(T \otimes 1)\sigma_s^*$. Then, $\alpha$ is a $C(Z)$-$E_0$-semigroup. Let $\cle^z:=\cle \otimes_{ev_z}B^z$. For $s \in P$, there exists a unique unitary operator $\sigma_s^z:\cle^z \otimes_{B^z} X^z_s \to \cle^z$ such that 
\[
\sigma_s^z((e \otimes ev_z(b))\otimes (u \otimes ev_z(c))=\sigma_s(e \otimes bu)\otimes ev_z(c).\]
Then, $(\cle^z,\sigma^z:=\{\sigma^z_s\}_{s \in P})$ is a left dilation of $X^z$. Let $\alpha^z$ be the $E_0$-semigroup associated with $\cle^z$. Then, for $T \in \cll_B(\cle)$ and $s \in P$, 
\[
\alpha_{s}^z(T \otimes 1)=\alpha_s(T) \otimes 1.
\]
If $\alpha$ is of compact type, then $\alpha^z$ is of compact type. 
Also, $\alpha^z$ is an $E_0$-semigroup over $P$ on $\cll_{B^z}(\cle^z)$, and the product system associated with $\alpha^z$ is $X^z$.  

\item[(3)] Suppose that $\alpha$ is a $C(Z)$-$E_0$-semigroup over $P$ on $\cll_B(\cle)$, where $\cle$ is a full Hilbert $B$-module. Let $X$ be the product system associated with $\alpha$. Then, $X$ is a $C(Z)$-product system.  For $s \in P$, let $\sigma_s:\cle \otimes_B X_s \to \cle$ be defined by 
\[
\sigma_s(x \otimes (y^* \otimes z))=\alpha_s(xy^*)z.
\]
Then, $(\cle,\sigma)$ is a left dilation of $X$. 

Define $\cle^z:=\cle \otimes_B  B^z=\cle \otimes_{ev_z} B^z$. Then, there exists a unique $E_0$-semigroup $\alpha^z:=\{\alpha^z_{s}\}_{s \in P}$ on $\cll_{B^z}(\cle^z)$ such that 
\[
\alpha^z_{s}(T \otimes 1)=\alpha_s(T) \otimes 1
\]
for $T \in \cll_B(\cle)$. Also, $X^z$ is the product system associated with $\alpha^z$.  
Moreover, thanks to Remark \ref{lem:comp-faithful}, if $\alpha$ is of compact type, then $\alpha^z$ is of compact type. 
 Note that $\clk_B(\cle)$ carries a  $C(Z)$-structure with fibre $\clk_B(\cle)^z =\clk_{B^z}(\cle^z)$ (see~\cite[1.7]{Kumjian-1998-Fell-bundle-gpd}).

\end{enumerate}  
Thus, under good conditions, to compare the $K$-theory of $C_{red}^{*}(X)$ and $C_{red}^{*}(X^z)$, it suffices to consider the case of semigroup dynamical systems. 

Let $(A,P,\alpha)$ be a $C(Z)$-dynamical system. We can write $A$ as the section algebra $C(Z,\cla)$ of an upper semicontinuous bundle of \(C^*\)-algebras $\cla$ over $Z$, where the fibre $A^z$, for $z \in Z$, is given by 
$A^z:=A/I^z$. Here, $I^z:=C_0(Z\backslash\{z\})A$. Since, for $s \in P$, $a \in A$ and $f \in C(Z)$, $\alpha_s(f\cdot a)=f\cdot\alpha_s(a)$, for every $z \in Z$ and $s \in P$, there exists a unique endomorphism $\alpha^z_{s}$ of $A^z$ such that 
\[
\alpha^z_{s}(a+I^z)=\alpha_s(a)+I^z
\]
for every $a \in A$. Then, $(A^z,P,\alpha^z)$ is a semigroup dynamical system.

Let $\cla^{+}$ be the unitization of $\cla$, i.e. for every $z \in Z$, set $(\cla^z)^{+}:=\cla^z \oplus \bbc$ be the unitization of $\cla^z$. Let $\Gamma^{+}$ be the set of all sections $s=(s_0,f)$, where $s_0$ is a continuous section of $\cla$, and $f$ is a continuous function on $Z$. Then, $(\cla^{+},\Gamma^{+})$ is an upper semicontinuous bundle over $Z$ (see \cite{Blanchard}). For $z \in Z$ and $s \in P$, define $(\alpha^
{z}_{s})^{+}:(\mathcal{A}^z)^{+} \to (\cla^z)^{+}$ by 
\[
(\alpha^{z}_{s})^{+}(x,\lambda)=(\alpha^z_{s}(x),\lambda).
\]
For $s \in P$, let $\widetilde{\alpha}_{s}:C(Z,\cla^{+}) \to C(Z,\cla^{+})$ be the endomorphism defined by
\[
\widetilde{\alpha}_{s}(f,\phi)(z):=(\alpha^z_{s}(f(z)),\phi(z)).
\]
Let $\widetilde{A}:=C(Z,\cla^{+})$. 
Then, $(\widetilde{A},P,\widetilde{\alpha})$ is a unital semigroup dynamical system.  Moreover, thanks to the split-exact sequence 
\[
0 \longrightarrow \cla^z \longrightarrow (\cla^z)^{+} \longrightarrow \bbc \longrightarrow 0\]
for every $z \in Z$, we have the following ``natural'' short exact sequence 
\[
0 \longrightarrow A \longrightarrow \widetilde{A} \longrightarrow C(Z) \longrightarrow 0.
\]
which is also split-exact.  Note that the above short exact sequence is $P$-equivariant, where the action of $P$ on $C(Z)$ is trivial. The following proposition is now a consequence of Prop. \ref{split_exact_homotopy}.

\begin{ppsn}
\label{short_exact}
    With the above notation, the  sequence 
    \[
    0 \longrightarrow A \rtimes_{red} P \longrightarrow \widetilde{A} \rtimes_{red} P \longrightarrow C(Z)\rtimes_{red} P \cong C(Z)\otimes C_{red}^{*}(P) \longrightarrow 0
    \]
    is a split-exact sequence. 
\end{ppsn}

We need the following result, which is probably well-known to the experts, and could be considered the baby version of a result of Dadarlat (\cite[Thm. 1.1]{Dadarlat}).

\begin{ppsn}\label{prop-fibre-K-theory-isomor}
Let $Z$ be a locally compact metric space which has a basis of compact open sets. Let $\cla$ and $\clb$ be two upper semicontinuous bundles of separable $C^{*}$-algebras over $Z$, and let $A:=C_0(Z,\cla)$ and $B:=C_0(Z,\clb)$. We assume that the fibres of $A$ and $B$ are unital. For $z \in Z$, the fibre of $\cla$ and the fibre of $\clb$ over $z$ are denoted $A(z)$ and $B(z)$ respectively. 
Let $\displaystyle \phi:=\coprod_{z \in Z}\phi(z):\coprod_{z \in Z}A(z)\to \coprod_{z \in Z}B(z)$ be a $C(Z)$-homomorphism.    Suppose that for $i \in \{0,1\}$, $K_i(\phi(z)): K_i(A(z)) \to K_i(B(z))$ is an isomorphism for every $z \in Z$. Then, $K_i(\phi):K_i(A) \to K_i(B)$ is an isomorphism. 
\end{ppsn}
\textit{Proof.} 
We first consider the case when $Z$ is compact.
 To prove the injectivity of $K_0(\phi)$, let
 \(K_0(\phi)([p]-[q]) = 0\) for some \([p]-[q] \in K_0(A)\). Without loss of generality, we can assume that  $p,q \in \mathcal{P}_n(A)$, where $\mathcal{P}_n(A)$ denotes the set of projections in $M_n(A)$.  Then we have
 \[
 0 = K_0(\phi)([p]-[q]) = K_0(\phi)([p]) - K_0(\phi)([q]) = [\phi(p)] - [\phi(q)].
 \]
 For a point \(x_0 \in Z\), we have
 \[
 0 = [\phi(p)(x_0)] - [\phi(q)(x_0)] = K_0(\phi(x_0)([p(x_0)] - [q(x_0)]).
 \]
 The injectivity of \(K_0(\phi((x_0))\) ensures that \([p(x_0)] = [q(x_0)]\).
 Then, there exist $m_{x_0} \in \bbn$ and a  partial isometry \(v_{x_0}\in M_{n+m_{x_0}}(A(x_0))\) such that
 \[
 v^*_{x_0}v_{x_0} =\begin{bmatrix}
  	p(x_0) & 0 \\
  	0 & 1_{m_{x_0}}
  \end{bmatrix} \quad \textup{and}\quad  v_{x_0}v_{x_0}^* = \begin{bmatrix}
  q({x_0}) & 0 \\
  0 & 1_{m_{x_0}}
  \end{bmatrix}.
  \]
   By Lemma~2.10 of~\cite{Blanchard}, there exist a neighbourhood  \(U_{x_0}\) of \(x_0\) and a partial isometry \(V_{x_0} \in M_{n+m_{x_0}}(A)\) such that \(V_{x_0}(x_0) =v_{x_0}\) and for every $x \in U_{x_0}$, 
  \[
  V_{x_0}^*(x)V_{x_0}(x) = \begin{bmatrix}
  	p(x) & 0 \\
  	0 & 1_{m_{x_0}}
  \end{bmatrix} \quad \textup{and}\quad  V_{x_0}(x)V_{x_0}^*(x) = \begin{bmatrix}
  	q({x}) & 0 \\
  	0 & 1_{m_{x_0}}
  	\end{bmatrix}.
  \]
  As $Z$ has a basis consisting of clopen sets, we can assume that $U_{x_0}$ is clopen. Then, $\{U_{x}:x \in Z\}$ covers $X$. 
    Choose a finite sub-cover \(\{U_{x_i}\}_{i=0}^{k}\) such that \(Z = \bigcup_{i=0}^{k}U_{x_i}\).  We can assume that \(U_{x_i}\)'s are mutually disjoint. Otherwise, replace $U_{x_i}$ by $\widetilde{U}_{x_i}=U_{x_i}\backslash (\bigcup_{k=0}^{i-1}U_{x_k})$. 
 
Set \(m = \sum_{i=0}^{k}m_{x_i}\).  We consider the partial isometry \(V_i = V_{x_i} \oplus 1_{m-m_{x_i}}\) for \(i=0,1,\cdots, k\). Define 
 \[
 V:=\sum_{i=0}^{k}1_{U_{x_i}} V_{i}.
 \] 
 Then, \(V\in M_{n+m}(A)\) and
  \[
 V^*V = \begin{bmatrix}
 	p & 0 \\
 	0 & 1_{m}
 \end{bmatrix} \quad \textup{and}\quad  VV^* = \begin{bmatrix}
 	q & 0 \\
 	0 & 1_{m}
 \end{bmatrix}.
 \]
 Therefore, \([p] =[q]\) and \(K_0(\phi)\) is injective.

Next, we claim that \(K_0(\phi)\) is surjective. Let \([p] - [q] \in K_0(B)\). We can assume that there exists $n$ such that $p,q \in \mathcal{P}_n(B)$. 
For a point \(x_0 \in Z\), we have \([p(x_0)] - [q(x_0)] \in K_0(B(x_0))\). Since \(K_0(\phi(x_0))\) is surjective, there exist an element \([r_{x_0}] - [s_{x_0}] \in K_0(A(x_0))\) such that 
 \[
 [p(x_0)] - [q(x_0)] = K_0(\phi(x_0))([r_{x_0}] - [s_{x_0}]) = [\phi(x_0)(r_{x_0})] - [\phi(x_0)(s_{x_0})].
 \]
 The above equation can be written as 
 \[
 [\phi(x_0)(r_{x_0})] + [q(x_0)] = [\phi(x_0)(s_{x_0})] + [p(x_0)].
 \]
 We can assume that there exists $m_{x_0}$ such that $r_{x_0},s_{x_0} \in \mathcal{P}_{m_{x_0}}(A(x_0))$. 
   Then, there exist \(t_{x_0}\in \mathbb{N}\) and a  partial isometry \(v_{x_0} \in M_{n+ t_{x_0}+m_{x_0}}(B(x_0))\) such that 
 \[
 v^*_{x_0}v_{x_0} =\begin{bmatrix}
 	p(x_0) & 0 & 0\\
 	0 & \phi(x_0)(s_{x_0}) & 0\\
 	0 & 0 &  1_{t_{x_0}}
 \end{bmatrix} \quad \textup{and}\quad  v_{x_0}v_{x_0}^* = \begin{bmatrix}
 	q({x_0}) & 0 & 0\\
 	0 & \phi(x_0)(r_{x_0}) & 0\\
 	0 & 0 & 1_{t_{x_0}}
 \end{bmatrix}.
 \]
 
  Then, by Lemma~2.10 of~\cite{Blanchard}, there exist a neighbourhood \(U_{x_0}\) of \(x_0\), a partial isometry \(V_{x_0} \in M_{n+ t_{x_0}+ m_{x_0}}(A)\) and sections \(R_{x_0}, S_{x_0} \in \mathcal{P}_{m_{x_0}}(A)\) such that \(V_{x_0}(x_0) =v_{x_0}\), \(R_{x_0}(x_0) = r_{x_0}, S_{x_0}(x_0) =s_{x_0}\) and   
 \[
 V_{x_0}^*(x)V_{x_0}(x) = \begin{bmatrix}
 	p(x) & 0 & 0 \\
 	0 & \phi(x)(S_{x_0}(x)) & 0\\
 	0 & 0 & 1_{t_{x_0}}
 \end{bmatrix} \textup{ and } V_{x_0}(x)V_{x_0}^*(x) = \begin{bmatrix}
 	q({x}) & 0  & 0\\
 	0 & \phi(x)(R_{x_0}(x)) & 0\\  
 	0 & 0 &1_{t_{x_0}}
 \end{bmatrix}
 \]
 for \(x\in U_{x_0}\).
 We can assume $U_{x_0}$ is clopen. As \(Z\) is compact, we can choose a finite collection \(\{U_{x_i}\}_{i=0}^{k}\)  such that \(Z = \bigcup_{i=0}^{k}U_{x_i}\). Moreover, as earlier, we can assume  \(U_{x_i}\)'s are mutually disjoint.
 
 Now set \(m = \sum_{i=0}^{k}m_{x_i},\) and let \(t = \sum_{i=0}^{k} t_{x_i}\).
   For $i=0,1,2,\cdots,k$, consider the partial isometry \(V_i = V_{x_i} \oplus 1_{m +t-m_{x_i} -t_{x_i}}\). Define \(V=\sum_{i=0}^{k}1_{U_{x_i}} V_{i} \in  M_{n+m+t}(B)\).
   We also define the projections \(S_{i} = S_{x_i} \oplus 1_{t-t_{x_i}}\) and \(R_{i} = R_{x_i} \oplus 1_{t-t_{x_i}}\) for \(i=0,1,\cdots, k\).  Set 
  \(S = \sum_{i=0}^{k}1_{U_{x_i}} S_{i}\) and \(R = \sum_{i=0}^{k}1_{U_{x_i}} R_{i}\). Then we have
 \[
 V^*V = \begin{bmatrix}
 	p & 0 & 0 \\
 	0 & \phi(S) & 0\\
 	0 & 0 & 1_{m}
 \end{bmatrix} \quad \textup{and}\quad  VV^* = \begin{bmatrix}
 	q & 0 & 0 \\
 	0 & \phi(R) & 0 \\
 	0 & 0 & 1_{m}
 \end{bmatrix}.
 \]
 Therefore, \([p] + [\phi(S)] =[q] + [\phi(R)]\) and hence \(\phi([R] -[S]) = [p] -[q]\). Therefore, \(K_0(\phi)\) is surjective.
 
 The proof of \(K_1\) is similar to \(K_0\), so we omit the proof.

	 Now suppose \(Z\) is a locally compact space. We choose an increasing sequence of compact open sets \((Z_i)_{i}\) such that \(Z = \bigcup_{i}Z_i\). Let \(\mathcal{A}|_{Z_i}\) be the restriction of the bundle \(\mathcal{A}\) on \(Z_i\) and \(A_i = C_0(Z_i; \mathcal{A}|_{Z_i})\).  We can view $A_i$ as a subalgebra of $A$ as $Z$ is clopen. Since, \(Z_i\subseteq Z_{i+1}\), we can view $A_i$ as a subalgebra of $A_{i+1}$, where the connecting map  \(A_i \ni f \mapsto \hat{f} \in  A_{i+1}\) is given by 
	 \[
	\hat{f}(x)= \begin{cases}
	 	  f(x) \quad \textup{if } x\in Z_i\\
	 	           0 \quad \textup{if } x\in Z_{i+1}\setminus Z_i.
	 \end{cases}
	 \]
Similarly, define $B_i$. Note that $\phi(A_i) \subset B_i$, and let $\phi_i$ be the restriction of $\phi$. 
     
	  Then,   
	 \[\varinjlim A_i \cong A, ~~\varinjlim B_i \cong B, \textrm{~~and~~}\phi=\varinjlim \phi_i.\] 
	 Since \(K_*(\phi(z)): K_*(A(z)) \to K_*(B(z))\) is an isomorphism for \(z\in Z\), the first part ensures that \(K_*(\phi_i): K_*(A_i) \to K_*(B_i)\) is an isomorphism for all \(i\). As \(K\)-theory respects inductive limits, we have that \(K_*(\phi): K_*(A) \to K_*(B)\) is an isomorphism.
\hfill $\Box$

We now prove Thm.~\ref{main_invariance}.

\textit{Proof of Thm.~\ref{main_invariance}.} Let us recall the notation and the hypothesis. Let $(A,P,\alpha)$ be a $C[0,1]$-semigroup dynamical system. For $z \in [0,1]$, let $A^z$ be the fibre, and let $ev_z:A \to A^z$ be the evaluation map. We have assumed that the following conditions are satisfied. \begin{enumerate}
     \item The map $K_*(ev_z):K_*(A) \to K_*(A^z)$ is an isomorphism.
     \item The group $G$ satisfies the Baum-Connes conjecture with coefficients and is torsion-free. 
     \item The pair $(P,G)$ satisfies the Toeplitz condition.
     \item Every element of $\widetilde{\Omega}$ is directed.
\end{enumerate}
Let $z \in [0,1]$ be given. 
As a consequence of  Corollary \ref{functorial_ideal}, Corollary \ref{equivariance implies homomorphism}, Prop.~\ref{split_exact_homotopy} and Prop.~\ref{useful_split_exact}, we have the following commutative diagram 
	\begin{center}
	\begin{tikzcd}
	0 \arrow[r] \arrow[dr, phantom,]
		& A\rtimes_{red} P \arrow[r] \arrow[d] & \widetilde{A}\rtimes_{red} P  \arrow[r] \arrow[d] & C([0,1])\otimes C^*_{red}(P) \arrow[r] \arrow[d] & 0\\
		0 \arrow[r] & A^z\rtimes_{red} P \arrow[r] & \widetilde{A}^z\rtimes_{red} P \arrow[r]&C^*_{red}(P) \arrow[r] & 0
	\end{tikzcd}
	\end{center}	
whose top and bottom rows are split-exact.

Applying the functor $K_*$ and noting that the functor $K_*$ preserves split-exactness, we obtain the following commutative diagram whose top and bottom rows are also split-exact. 
	\begin{center}
		\begin{tikzcd}
			0 \arrow[r] \arrow[dr, phantom,]
			& K_*(A\rtimes_{red} P) \arrow[r] \arrow[d] & K_*(\widetilde{A}\rtimes_{red} P)  \arrow[r] \arrow[d] & K_*(C([0,1])\otimes C^*_{red}(P)) \arrow[r] \arrow[d] & 0\\
			0 \arrow[r] & K_*(A^z\rtimes_{red} P) \arrow[r] & K_*(\widetilde{A}^z \rtimes_{red} P) \arrow[r] &K_*(C^*_{red}(P)) \arrow[r] & 0
		\end{tikzcd}
	\end{center}	
	Thanks to the five lemma, to conclude the result,  we can assume $A$ is unital and $\alpha_s$ is unital for every $s \in P$.

Let $(\widetilde{\mathcal{D}},\widetilde{\Omega}\rtimes G)$ be the groupoid dynamical system considered in Section~\ref{sec-gpd-presentation} that corresponds to $(A,P,\alpha)$. Similarly, let $(\widetilde{\mathcal{D}}^{z},\widetilde{\Omega}\rtimes G)$ be the groupoid dynamical system  that corresponds to $(A^z,P,\alpha^z)$. It follows from Thm.~\ref{unital_case_KS} that $A \rtimes_{red} P$ and $C_{0}(\widetilde{\Omega},\widetilde{\mathcal{D}}) \rtimes_{red} G$ are Morita-equivalent, and $A^z \rtimes_{red} P$ is Morita-equivalent to $C_{0}(\widetilde{\Omega},\widetilde{\mathcal{D}}^z)\rtimes_{red} G$. Thus, it suffices to prove that $C_{0}(\widetilde{\Omega},\widetilde{\mathcal{D}}) \rtimes_{red} G$ and $C_{0}(\widetilde{\Omega},\widetilde{\mathcal{D}}^{z})\rtimes_{red} G$ have the same $K$-theory. 

Let $\widetilde{D}=C_0(\widetilde{\Omega},\widetilde{\mathcal{D}})$ and $\widetilde{D}^{z}=C_{0}(\widetilde{\Omega},\widetilde{\mathcal{D}}^z)$. Recall that $\widetilde{D} \subset \ell^{\infty}(G,A)$ and $\widetilde{D}^z \subset \ell^{\infty}(G,A^z)$. Let $\widetilde{ev}_z: \ell^{\infty}(G,A) \to \ell^{\infty}(G,A^z)$ be defined by 
\[
\widetilde{ev}_z(f)(s)=ev_z(f(s)).
\]
For $g \in G$, $a \in A$, recall that $j_g(a) \in \ell^{\infty}(G,A)$ is defined as 
\begin{equation}
j_g(x)(h):=\begin{cases}
    \alpha_{hg^{-1}}(x)& \mbox{ if
} hg^{-1} \in P,\cr
   &\cr
    0 &  \mbox{ if } hg^{-1} \notin P.
         \end{cases}
\end{equation}
For $g \in G$, $a \in A$, $j_g(a+I^z) \in \ell^{\infty}(G,A^z)$ is similarly defined. To avoid confusion, we denote $j_g(a+I^z)$ by $j_g^{z}(a+I^z)$. Observe that 
\[
\widetilde {ev}_z(j_g(a))=j_g^{z}(a+I^z).
\]
for $g \in G$ and $a \in A$. 
Thus, $\widetilde{ev}_z$ maps $\widetilde{D}$ to $\widetilde{D}^z$. Moreover, $\widetilde{ev}_z$ is a $C_0(\widetilde{\Omega})$ homomorphism, and  is $G$-equivariant. Thus, $\widetilde{ev}_z$ induces a map $\epsilon^z: C_{0}(\widetilde{\Omega},\widetilde{\mathcal{D}})\rtimes_{red}  G=\widetilde{D}\rtimes_{red} G \to \widetilde{D}^z \rtimes_{red} G=C_{0}(\widetilde{\Omega},\widetilde{\mathcal{D}}^z)\rtimes_{red} G$. By definition, $\epsilon^z|_{\widetilde{D}}=\widetilde{ev}_z$. 

We claim that for $i=0,1$, $K_i(\epsilon^z):C_{0}(\widetilde{\Omega},\widetilde{\mathcal{D}}^z) \rtimes_{red} G \to C_0(\widetilde{\Omega},\widetilde{\mathcal{D}}^z) \rtimes_{red} G$ is an isomorphism. 
For $F \in \widetilde{\Omega}$, we denote the fibre of $\widetilde{\mathcal{D}}$ over $F$ by $\widetilde{\mathcal{D}}_F$. Similarly, the fibre of $\widetilde{\mathcal{D}}^z$ over $F$ is denoted $\widetilde{\mathcal{D}}^z_F$. 
 Since the group \(G\) is torsion-free and satisfies the Baum-Connes conjecture with coefficients,  the descent principle (see~\cite[Prop. 2.1] {Echterhoff_Nest_Oyono}) ensures that the above claim reduces to the claim: \(K_i(\widetilde{ev}_z) \colon K_i(\widetilde{\mathcal{D}}) \to  K_i(\widetilde{\mathcal{D}}^z)\) is an isomorphism for \(i=0,1\).
 
 Let \(F\in \widetilde{\Omega}\). Again Prop.~\ref{inductive limit} says that 
 \[
 \widetilde{\mathcal{D}}_F \cong \underset{s\in F}{\varinjlim}  B_s \quad \textup{and} \quad  \widetilde{\mathcal{D}}^z_F \cong \underset{s\in F}{\varinjlim} B^z_s,
 \]
 where \(B_s =A\) and \(B^z_s =A^z\) and the connecting maps are as in Prop. \ref{inductive limit}. Under this identification, thanks to Eq. \ref{limitexists}, $\displaystyle \widetilde{ev}_{z}(F)= \underset{s\in F}{\varinjlim} ev_z$.  Since $K_i(ev_z):K_i(A) \to K_i(A^z)$ is an isomorphism and $K$-theory commutes with inductive limits, we have
	\(K_i(\widetilde{ev}_z(F)) \colon K_i(\widetilde{\mathcal{D}}_F) \to  K_i(\widetilde{\mathcal{D}}^z_F)\) is an isomorphism for \(i=0,1\). Therefore, by Prop.~\ref{prop-fibre-K-theory-isomor}, 
	 \(K_i(\widetilde{ev}_z) \colon K_i(\widetilde{D}) \to  K_i(\widetilde{D}^z)\) is an isomorphism for \(i=0,1\). This completes the proof.
\hfill $\Box$

\begin{crlre}
\label{homotopy_invariance_prod_version}
Let $P$ be a subsemigroup of a group $G$ such that $\bigcup_{n=1}^{\infty}Pa^{-n}=G$ for some  $a \in P$.    Let $X$ be a proper $C[0,1]$-product system over $P^{op}$ with coefficient algebra $B$. Suppose that $X_s$ is full for every $s \in P$. Assume that $B$ is unital and separable.  For $z \in [0,1]$, denote the fibre of $B$ over $z $ by $B^z$. Suppose that the following conditions are satisfied.
\begin{enumerate}
    \item[(1)] For every $z \in [0,1]$, the map $K_*(ev_z):K_*(B) \to K_*(B^z)$ is an isomorphism.
    \item[(2)] The group $G$ is torsion-free and satisfies the Baum-Connes conjecture with coefficients. 
    \item[(3)] Every element of $\widetilde{\Omega}$ is directed.
    \end{enumerate}
Then, $K_*(C_{red}^{*}(X))$ and $K_*(C_{red}^{*}(X^z))$ are isomorphic for every $z \in [0,1]$. 
\end{crlre}
\textit{Proof.} Let $(\cle,\sigma:=\{\sigma_s\}_{s \in P})$ be a left dilation of $X$. By Thm.~\ref{main_existence}, such a left dilation exists.  Let $\alpha$ be the $E_0$-semigroup associated with  $(\cle,\sigma)$, i.e. for $s \in P$, the endomorphism $\alpha_s:\cll_B(\cle) \to \cll_B(\cle)$ is given by 
\[
\alpha_s(T)=\sigma_s(T \otimes 1)\sigma_s^{*}.\]
Then, $X$ is the product system associated with $\alpha$. As discussed in the beginning of this section, $\alpha$ is a $C([0,1]$-semigroup, and $X^z$ is the product system associated with $\alpha^z$. Moreover, $\alpha^z$ is of compact type. Since $C_{red}^{*}(X)$ is Morita equivalent to $\clk_{B}(\cle) \rtimes P$ and $C_{red}^{*}(X^z)$ is Morita equivalent to $\clk_{B^t}(\cle^z)$, it suffices to prove that $\clk_B(\cle) \rtimes_{red} P$ and $\clk_{B^z}(\cle^z) \rtimes_{red} P$ have the same $K$-theory. 

Observe that $(\clk_B(\cle),P,\alpha)$ is a $C[0,1]$-semigroup dynamical system, and, for every $z \in [0,1]$, $(\clk_B(\cle)^z,P,\alpha^z)=(\clk_{B^z}(\cle^z),P,\alpha^z)$. We can now apply Thm.~\ref{main_invariance}. The only thing that requires justification is that, for every $z \in [0,1]$,  the evaluation map $\widetilde{ev}_z:\clk_B(\cle) \to \clk_{B}(\cle)^z= \clk_{B^z}(\cle^z)$ induces an isomorphism at the $K$-theory level. 

Let $z \in [0,1]$ be given. Let $[\cle]$ denote the element in $KK(\clk_B(\cle),B)$  given by the imprimitivity bimodule $\cle$, and similarly, let $[\cle^z]$ denote the element in $KK(\clk_{B^z}(\cle^t),B^z)$. Then,
\[
[\cle] \cdot [ev_z] = [\widetilde{ev}_z] \cdot [\cle^z].
\]
In the above, $\cdot$ denotes the $KK$-product. Note that  $[\cle],[\cle^z]$ are invertible and $[ev_z]$ is invertible by assumption. Hence, $[\widetilde{ev}_z]$ is invertible. Thus, $\widetilde{ev}_z$ induces an isomorphism at the $K$-theory level. \hfill $\Box$

\textit{Notation:} For a $C^{*}$-algebra $A$, $C([0,1],A)=C[0,1]\otimes A$ denotes the $C^{*}$-algebra of continuous $A$-valued functions.  The algebra $C([0,1],A)$ is a $C([0,1])$-algebra with the $C([0,1])$-structure given by 
\[
f\cdot (g \otimes a)=fg \otimes a\]
for $f \in C([0,1])$ and $a \in A$. 

\begin{dfn}
\label{homotopy_definition_Gillaspy}
\begin{enumerate}
    \item  Let $X_0,X_1$ be product systems over $P$ with the same coefficient algebra~$B$. We say that $X_0$ and $X_1$ are homotopic if there exists a $C[0,1]$-product system $\mathcal{X}$ with coefficient algebra $\mathcal{B}=C[0,1]\otimes B$  such that $\mathcal{X}^{0} \cong X_0$ and $\mathcal{X}^{1}\cong X_1$. If $X_0$ and $X_1$ are proper, we demand that $\mathcal{X}$ is proper. We call such a $\mathcal{X}$ a homotopy between $X_0$ and $X_1$. If the fibres of $X_0$ and $X_1$ are full, we demand that the fibres of $\mathcal{X}$ are also full. 

   \item  Let $(A,P,\alpha)$ and $(A,P,\beta)$ be semigroup dynamical systems. Then, $(A,P,\alpha)$ and $(A,P,\beta)$ are said to be homotopic if there exists a $C[0,1]$-semigroup dynamical system $(C,P,\gamma)$ with $C=C[0,1]\otimes A$, $(C^0,P,\gamma^0) \cong (A,P,\alpha)$ and $(C^1,P,\gamma^1)=(A,P,\beta)$. The semigroup dynamical system $(C,P,\gamma)$ is called a homotopy between $(A,P,\alpha^0)$ and $(A,P,\alpha^1)$.
\end{enumerate}
    \end{dfn}

    \begin{rmrk}
        The notion of homotopy of product systems as defined in Defn.~\ref{homotopy_definition_Gillaspy} was earlier considered in~\cite[Defn. 3.3]{Gillaspy}. 
    \end{rmrk}

The following are immediate corollaries of Thm.~\ref{main_invariance} and Thm.~\ref{homotopy_invariance_prod_version}. 

\begin{crlre} Let $P$ be a subsemigroup of a torsion-free group $G$. 
    Let $(A,P,\alpha)$ and $(A,P,\beta)$ be semigroup dynamical systems that are homotopic. Suppose that $(P,G)$ satisfies the Toeplitz condition, $G$ satisfies the Baum-Connes conjecture with coefficients, and every element of $\widetilde{\Omega}$ is directed. 
    Then, $K_i(A \rtimes_{red,\alpha} P)$ and $K_i(A \rtimes_{red,\beta} P)$ are isomorphic for \(i\in \{0,1\}.\)
\end{crlre}

\begin{crlre}
\label{homotopoy_invariance_product_two}
    Let $P$ be a subsemigroup of a torsion-free group $G$ that satisfies the Baum-Connes conjecture with coefficients. Assume that $PP^{-1}=G$, $P$ has an order unit, and every element of $\widetilde{\Omega}$ is directed. Let $X_0$ and $X_1$ be proper product systems over $P$ with coefficient algebra $B$ that is separable and unital. Suppose that the fibres of $X_0$ and $X_1$ are full. If $X_0$ and $X_1$ are homotopic, then  $K_*(C_{red}^{*}(X_0))$ and $K_*(C_{red}^{*}(X_1))$ are isomorphic. 
\end{crlre}

\subsection{Examples} In this subsection, we give a few examples of homotopic semigroup dynamical systems and product systems. 

\begin{xmpl}
\label{abstract_example_1}
Let $A$ be a separable $C^{*}$-algebra. Suppose for every $z\in [0,1]$, we have an $E_0$-semigroup $\alpha^z:=\{\alpha^z_s\}_{s \in P}$ on $M(A)$ which is of compact type. Suppose that for $a \in A$ and $s \in P$, the map 
\[
[0,1] \ni z \mapsto \alpha^z_s(a) \in A
\]
is norm continuous. Then, $(A,P,\alpha^0)$ and $(A,P,\alpha^1)$ are homotopic. To see this, consider  $C:=C[0,1] \otimes A=C([0,1],A)$.  For $f \in C([0,1],A)$, $s \in P$, define $\gamma_s:C \to C$ by 
\[
\gamma_s(f)(z)=\alpha^z_s(f(z)).\]
Then, $\overline{\alpha_s(C)C}=C$ for every $s \in P$. Also, $(C,P,\gamma)$ is a $C[0,1]$-semigroup dynamical system and is a homotopy between $(A,P,\alpha^0)$ and $(A,P,\alpha^1)$.
\end{xmpl}

A concrete example that comes under the framework of Example \ref{abstract_example_1} is given below.
\begin{xmpl}
    \label{concrete_example_1}
    Let $K$ be a separable Hilbert space, and let $U:=\{U_g\}_{g \in G}$ be a group of unitary operators on $K$. Let $H$ be a closed subspace of $K$ such that $U_s(H) \subset H$ for every $s \in P$. For $s \in P$, let $V_s:=U_{s}\big|_{H}$. Let $c:G \to \bbr$ be a homomorphism. For $z \in [0,1]$, set $V_{s}^{z}:=e^{izc(s)}V_s$. Then,  there exists, for every $z \in [0,1]$, a semigroup $\alpha^z:=\{\alpha^z_{s}\}_{s \in P}$ of unital endomorphisms of $\cla(H)$ (see Example~\ref{CAR_example_one}) such that 
    \[
    \alpha_{s}^{z}(a(\xi))=a(V_{s}^z\xi)=e^{izc(a)}a(V_s\xi)
    \]
    for every $\xi \in H$. Then, $(\cla(H),P,\alpha^0)$ and $(\cla(H),P,\alpha^1)$ are homotopic. 
    As an example of $(U,K,H)$, we can take $K=\ell^2(G)$, $U$ the left regular representation, and $H=\ell^2(P)$.
    \end{xmpl}

\begin{xmpl}
\label{Deforming multiplication}
Let $B$ be a separable $C^{*}$-algebra. For each $s \in P$, let $X_s$ be a $C^{*}$-correspondence from $B$ to $B$. 
Assume that for each $z \in [0,1]$ and $s,t \in P$, we have a unitary bimodule map $U_{s,t}^{z}:X_s \otimes_B X_t \to X_{st}$. Let $z \in [0,1]$. Define a multiplication $\odot_z$ on $X:=\coprod_{s \in P}X_s$ by setting 
\[
u\odot_{z} v=U_{s,t}^{z}(u \otimes v).
\]
We suppose that $X$ with the multiplication rule $\odot_z$ is a product system. We denote the resulting product system by $X^z$. 

Suppose that  for $u \in X_s$ and $v \in X_t$, the map 
\[
[0,1] \ni z \mapsto U_{s,t}^{z}(u \otimes v) \in X_{st}
\]
is norm continuous. 

For $s \in P$, let $\mathcal{X}_s:=C([0,1],X_s)=C[0,1]\otimes X_s$ be the external tensor product which is a Hilbert $IB:=C([0,1],B)=C[0,1]\otimes B$-module. It carries a left $IB$ action given by 
\[
f\cdot \phi(z)=f(z)\phi(z)\]
for $f \in IB$ and $\phi \in \mathcal{X}_s$. For $s,t \in P$, $\phi \in \mathcal{X}_s$ and $\psi \in \mathcal{X}_t$, define 
\[
(\phi \odot \psi)(z)=U_{s,t}^{z}(\phi(z) \otimes \psi(z)).
\]
Then, $\mathcal{X}:=\{\mathcal{X}_s\}_{s \in P}$ with the product $\odot$ is a $C[0,1]$ product system. Moreover, the product system $\mathcal{X}$ is a homotopy between $X^0$ and $X^1$. 
The Hilbert space version of this example for $P=\bbn^k$ was also discussed in \cite[Section 4]{Gillaspy}. 
\end{xmpl}

Two concrete examples that fit within the setup of Example~\ref{Deforming multiplication} are given below. In the first example, we deform the product rule by a \(2\)-cocycle. 

\begin{xmpl}
    Let $X:=\{X_s\}_{s \in P}$ be a product system over $P$. The product rule on $X$ is denoted by $\odot$. Let $\{U_{s,t}:X_{s} \otimes_B X_t \to X_{st} \}_{s,t \in P}$ be the unitaries given by the multiplication rule. Let $c:G \times G \to \bbr$ be a \(2\)-cocycle, i.e. for $r,s,t \in G$, 
    \[
    c(r,s)+c(r+s,t)=c(r,s+t)+c(s,t).\]
    Define a new product rule $\odot$ on $X=\coprod_{s \in P}X_s$ by setting 
    \[
    u \odot_c v:=e^{ic(s,t)}u\odot v
    \]
    for $u \in X_s$ and $v \in X_t$.  Denote the resulting product system by $X_c$.
    
    For $z \in [0,1]$ and for $s,t \in P$, define $U_{s,t}^{z}:=e^{izc(s,t)}U_{s,t}$. Then, the conditions of Example~\ref{Deforming multiplication} are satisfied. Hence, the product system $X$
    and $X_c$ are homotopic whose reduced $C^{*}$-algebras have the same $K$-theory if for every $s \in P$, $X_s$ is proper, full and if  $(P,G)$ satisfies the conditions mentioned in Corollary~\ref{homotopoy_invariance_product_two}. In particular, if $X$ is proper and $X_s$ is full for every $s$, and if $P$ is a finitely generated subsemigroup of an abelian group, then $C_{red}^{*}(X)$ and $C_{red}^{*}(X_c)$ have the same $K$-theory. 
\end{xmpl}

\begin{xmpl}
Keep the notation of Example~\ref{concrete_example_1}. For $s \in P$, set $X_s:=\Gamma_a(Ker(V_s^*))$. Here, for a Hilbert space $L$, $\Gamma_a(L)$ denotes the anti-symmetric Fock space.    Define a product $\odot$ on the disjoint union $X:=\coprod_{s \in P}X_s$ as follows: 
\[
(\xi_1 \wedge \xi_2 \wedge \cdots \wedge \xi_m) \odot (\eta_1 \wedge \eta_2 \wedge \cdots \wedge \eta_n):=V_{t}\eta_1 \wedge V_t\eta_2 \wedge \cdots \wedge V_t\eta_n \wedge \xi_1 \wedge \xi_2 \wedge \cdots \wedge \xi_m\]
for $\xi_1,\xi_2,\cdots, \xi_m \in Ker(V_s^*)$ and $\eta_1,\eta_2,\cdots,\eta_n \in Ker(V_t^*)$. Then, $X$ is a product system of  Hilbert spaces. In the theory of $E_0$-semigroups, the product system $X$ is called the product system of the CAR flow associated with $V$, and CAR flows are one of the well-studied examples. Observe that $X$ is proper if and only if $\dim Ker(V_s^*)<\infty$ for every $s \in P$.

Let $z \in [0,1]$. For $g \in P$, let $U^z_{g}:=e^{izc(g)}U_g$, and for $s \in P$, let $V^z_{s}:=U^z_{s}|_{H}$. Let $X^z$ be the product system of the CAR flow associated with $V^z$. Note that the fibres of $X^z$ and $X$ remain the same; only the multiplication rule is changed. For $z \in [0,1]$ and $s,t \in P$, define 
\[
U_{s,t}^{z}((\xi_1 \wedge \xi_2 \wedge \cdots \wedge \xi_m)\otimes (\eta_1 \wedge \eta_2 \wedge \cdots \wedge \eta_n))=e^{inzc(a)}V_{t}\eta_1 \wedge V_t\eta_2 \wedge \cdots \wedge V_t\eta_n \wedge \xi_1 \wedge \xi_2 \wedge \cdots \wedge \xi_m.
\]
The conditions of Example \ref{abstract_example_1} are satisfied, and we can conclude that $X$ and $X^z$ are homotopic for every $z$. Thus, $C_{red}^{*}(X)$ and $C_{red}^{*}(X^z)$ have isomorphic $K$-groups if the hypotheses of Corollary \ref{homotopoy_invariance_product_two} are satisfied. 

For a concrete example that satisfies the conditions mentioned above,  let $P$ be a numerical subsemigroup, i.e. $P$ is a semigroup of $\bbn_0$ such that $\bbn_0\backslash P$ is finite. Let $U$ be the regular representation of $\bbz$ on $\ell^2(\bbz)$, and let $H:=\ell^2(P)$. Then, $Ker(V_s^*)$ is finite dimensional where $V_s=U_{s}|_H$. 

\end{xmpl}

 We end our paper with a final remark. 
\begin{rmrk}
     We believe that many of our results should have a Cuntz-Pimsner version, and our techniques can be applied. It would be interesting to see whether our methods can be extended to the non-proper case and beyond Ore semigroups. Does  the conclusion of Thm.~\ref{main_existence} hold in the quasi-lattice ordered case? 
\end{rmrk}


\end{document}